\documentclass{amsart}
\usepackage{amsmath,amssymb,amscd}
\usepackage[dpi=600,PostScript=dvips]{diagrams}
\usepackage{bbold}
\usepackage{psfrag}
\usepackage[dvips]{graphicx}
\usepackage{subfigure}
\newtheorem{thm}{Theorem}[section]
\newtheorem{cor}[thm]{Corrolary}
\newtheorem{lem}[thm]{Lemma}

\newtheorem{prop}[thm]{Proposition}
\theoremstyle{definition}
\newtheorem{defn}[thm]{Definition}
\newtheorem{remsbold}[thm]{Remarks}
\newtheorem{rembold}[thm]{Remark}
\theoremstyle{remark}

\numberwithin{equation}{section}

\newcommand{\kt}{$\Bbbk$\nobreakdash-\hspace{0pt}}
\newcommand{\id}{\mathrm{id}}
\newcommand{\bd}{\mathrm{Bd}}
\newcommand{\Ob}{\mathrm{Ob}}

\newcommand{\ptsum}{{\textstyle \sum }}

\newcommand{\al}{\alpha}
\newcommand{\be}{\beta}

\newcommand{\kk}{\Bbbk}

\newcommand{\opp}{\mathrm{op}}

\newcommand{\cc}{\mathcal{C}}
\newcommand{\aaa}{\mathcal{A}}
\newcommand{\vv}{\mathcal{V}}
\newcommand{\ww}{\mathcal{W}}
\newcommand{\bb}{\mathcal{B}}
\newcommand{\dd}{\mathcal{D}}

\newcommand{\Z}{\mathbb{Z}}

\newcommand{\ds}{\displaystyle}
\newcommand{\un}{\mathbb{1}}

\newcommand{\Endo}{\mathrm{End}}
\newcommand{\Homo}{\mathrm{Hom}}

\newcommand{\reph}{\mathrm{rep}_H}
\newcommand{\rephn}{\mathrm{rep}_{H_n}}
\newcommand{\Negl}{\mathrm{Negl}}
\newcommand{\Tr}{\mathrm{Tr}}
\newcommand{\tr}{\mathrm{tr}_q}
\newcommand{\qdim}{\mathrm{dim}_q}
\newcommand{\kdim}{\dim_\kk}
\newcommand{\R}{$R$\nobreakdash-\hspace{0pt}}
\newcommand{\norm}{\mathrm{norm}}

\newcommand{\Ki}{\mathcal{I}}
\newcommand{\RT}{\mathrm{RT}}
\newcommand{\sss}{s}

\newcommand{\ev}{\mathrm{ev}}
\newcommand{\tev}{\widetilde{\mathrm{ev}}}
\newcommand{\coev}{\mathrm{coev}}
\newcommand{\tcoev}{\widetilde{\mathrm{coev}}}

\begin{document}

\title{Kirby elements and quantum invariants}
\author{Alexis Virelizier}
\address{D\'epartement des Sciences Math\'ematiques\\
Universit\'e Montpellier II \\ Case Courrier 051 \\Place Eug\`ene Bataillon\\
34095 Montpellier Cedex 5\\ France}

\email{virelizi@math.univ-montp2.fr}


\date{\today}

\begin{abstract}
We define the notion of a Kirby element of a ribbon category~$\cc$ (not necessarily semisimple). Kirby elements lead to
3-manifolds invariants. We characterize a set of Kirby elements (in terms of the structural morphisms of a Hopf algebra
in $\cc$) which is sufficiently large to recover the quantum invariants of 3-manifolds of Reshetikhin-Turaev, of
Hennings-Kauffman-Radford, and of Lyubashenko. The cases of a semisimple ribbon category and of a category of
representations are explored in details.
\end{abstract}
\maketitle

\setcounter{tocdepth}{1} \tableofcontents

\section*{Introduction}

During the last decade, deep connections between low-dimensional topology and the purely algebraic theory of quantum
groups (or, more generally, of braided categories) were highlighted. In particular, this led to a new class of
3-manifolds invariants, called {\it quantum invariants}, defined in several ways.

The aim of the present paper is to give a as general as possible method of constructing quantum invariants of 3-manifolds
starting from a ribbon category or a ribbon Hopf algebra. With this formalism, we recover the 3-manifolds invariants of
Reshetikhin-Turaev~\cite{RT2,Tur2}, of Hennings-Kauffman-Radford~\cite{KR1,He2}, and of Lyubashenko~\cite{Lyu2} when
these are well-defined.

Let $\kk$ be a field and $\cc$ be a \kt linear ribbon category (not necessarily semisimple). Under some technical
assumption, namely the existence of a coend $A\in \Ob(\cc)$ of the functor $(X,Y)\in \cc^\opp \times \cc \mapsto X^*
\otimes Y \in \cc$, a scalar $\tau_\cc(L;\al)$ can be associated to any framed link $L$ in $S^3$ and any morphism
$\al\in\Homo_\cc(\un,A)$, see \cite{Lyu2}. Recall, see \cite{Lyu1}, that the object $A$ of $\cc$ is then a Hopf algebra
in the category $\cc$.

By a \emph{Kirby element} of $\cc$, we shall mean a morphism $\al \in\Homo_\cc(\un,A)$ such that $\tau_\cc(L;\al)$ is
invariant under isotopies of $L$ and under 2-handle slides. By using the Kirby theorem \cite{Ki}, we have that if $\al$
is a Kirby element of $\cc$ such that $\tau_\cc(\bigcirc^{\pm 1};\al)\neq 0$, then $\tau_\cc(L;\al)$ can be normalized to
an invariant $\tau_\cc(M_L;\al)$ of 3-manifolds, where $\bigcirc^{\pm 1}$ is the unknot with framing $\pm 1$ and $M_L$
denotes the 3-manifold obtained from $S^3$ by surgery along $L$.

In general, determining all the Kirby elements of $\cc$ is a quite difficult problem. In this paper we characterize, in
terms of the structure maps of the categorical Hopf algebra~$A$, a set $\Ki(\cc)$ made of Kirby elements of $\cc$ which
is sufficiently large to contain the Kirby elements corresponding to the known quantum invariants.

If the categorical Hopf algebra $A$ admits a two-sided integral $\lambda: \un \to A$, then $\lambda$ belongs to
$\Ki(\cc)$ and the corresponding invariant $\tau_\cc(M;\lambda)$ is the Lyubashenko invariant~\cite{Lyu2}.

When $\cc$ is semisimple, we give necessary conditions for being in $\Ki(\cc)$ and we show that there exist (even in the
non-modular case) elements of $\Ki(\cc)$ corresponding to the Reshetikhin-Turaev invariants \cite{RT2,Tur2} computed from
finitely semisimple full ribbon subcategories of $\cc$. Note that these elements are not in general two-sided integrals.

More generally, we show that $\Ki(\cc)$ contains the Kirby elements leading to the invariants obtained from $\Ki(\bb)$,
where $\bb$ is a finitely semisimple full ribbon subcategory of the semisimple quotient of $\cc$. We verify that, in
general, there exist Kirby elements in $\Ki(\cc)$ which are not of this last form. This means that the semisimplification
process ``lacks'' some invariants.

Let $H$ be a finite-dimensional ribbon Hopf algebra. Suppose that $\cc$ is the category $\reph$ of finite-dimensional
left $H$-modules. We describe $\Ki(H)=\Ki(\reph)$ in purely algebraic terms. One of the interest of such a description is
to avoid the representation theory of $H$ (which may be of wild type, see \cite{Ben}).

If $H$ is unimodular, then $1 \in \Ki(H)$ and the corresponding invariant $\tau_{(H,1)}$ is the Hennings-Kauffman-Radford
invariant \cite{KR1,He2}. More generally, and even if $H$ is not unimodular, we show that the invariant $\tau_{(H,z)}$ of
3-manifolds corresponding to $z \in \Ki(H)$ can be computed by using the Kauffman-Radford algorithm.

If $\vv$ is a set of simple left $H$-modules which makes $(H,\vv)$ a premodular Hopf algebra, then there exists $z_\vv
\in \Ki(H)$ such that $\tau_{(H,z_\vv)}$ is the Reshetikhin-Turaev invariant computed from $(H,\vv)$, which can then be
computed by using the Kauffman-Radford algorithm.

When $H$ is semisimple and $\kk$ is of characteristic $0$, we show that the Hennings-Kauffman-Radford invariant (computed
from $H$) and the Reshetikhin-Turaev invariant (computed from $\reph$) are simultaneously well-defined and coincide (even
in the non-modular case). In the modular case, this was first shown in \cite{Ker1}.

We explicitly determine $\Ki(H)$ for a family of non-unimodular ribbon Hopf algebras which contains Sweedler's Hopf
algebra.

As an algebraic application, the operators involved in the description of $\Ki(H)=\Ki(\reph )$ in algebraic terms allows
us to parameterize all the traces on a finite-dimensional ribbon Hopf algebra $H$. When $H$ is unimodular, we
recover the parameterization given in \cite{Rad1,He2}.\\

The paper is organized as follows. In Section~\ref{sect-1}, we review ribbon categories and coends. In
Section~\ref{sect-2}, we define and study Kirby elements. We focus, in Section~\ref{RTinv}, on the case of semisimple
ribbon categories and, in Section~\ref{HKRsect}, on the case of categories of representations of ribbon Hopf algebras. In
Section~\ref{Sect-examples}, we treat an example in detail. Finally, in Appendix~\ref{sect-app},
we study traces on ribbon Hopf algebras.\\

{\it Acknowledgements.} I thank A. Brugui\`eres for helpful discussions, in particular for questions concerning category
theory.

\section{Ribbon categories and coends}\label{sect-1}

In this section, we review some basic definitions concerning ribbon categories and coends. Throughout this paper, we let
$\Bbbk$ be a field.

\subsection{Ribbon categories}\label{sect-ribcat}

Let $\mathcal{C}$ be a strict monoidal category with unit object $\un$ (note that every monoidal category is equivalent
to a strict monoidal category in a canonical way, see, \cite{ML1}). A \emph{left duality} in $\mathcal{C}$ associates to
any object $U\in \mathcal{C}$ an object $U^*\in \mathcal{C}$ and two morphisms $\ev_U: U^* \otimes U \to \un$ and
$\coev_U:\un \to U \otimes U^*$ such that
\begin{align}
 & \label{catdual1} (\id_U \otimes \ev_U) (\coev_U \otimes \id_U) =\id_U \\
 & \label{catdual2} (\ev_U\otimes \id_{U^*}) (\id_{U^*}\otimes\coev_U) =\id_{U^*}.
\end{align}
We can (and we always do) impose that $\un^*=\un$, $\ev_\un=\id_\un$ and $\coev_\un=\id_\un$.

By a \emph{braided category} we shall mean a monoidal category $\mathcal{C}$ with left duality and endowed with a system
$\{c_{U,V}:U\otimes V \to V \otimes U\}_{U,V\in \mathcal{C}} $ of invertible morphisms (the \emph{braiding}) satisfying
the following three conditions:
\begin{align}
  & c_{U',V'} (f\otimes g)= (g \otimes f) \,c_{U,V},\label{catbraid1} \\
  & c_{U\otimes V,W} =(c_{U,W} \otimes \id_V) (\id_U \otimes c_{V,W}), \label{catbraid2}\\
  & c_{U,V\otimes W}= (\id_{V} \otimes c_{U,W})(c_{U,V} \otimes \id_W), \label{catbraid3}
\end{align}
for all objects $U,V,W\in \mathcal{C}$ and all morphisms $f:U\to U', g:V\to V'$ in $\mathcal{C}$. Note that, by applying
\eqref{catbraid2} to $U=V=\un$ and \eqref{catbraid3} to $V=W=\un$ and using the invertibility of $ c_{U,\un}$ and
$c_{\un,U}$, we obtain that $c_{U,\un}=c_{\un,U}=\id_U$ for any object $U\in \mathcal{C}$.

A \emph{ribbon category} is a braided category $\mathcal{C}$ endowed with a family of invertible morphisms
$\{\theta_U:U\to U\}_{U\in \mathcal{C}}$ (the \emph{twist}) satisfying the following conditions:
\begin{align}
  & \label{cattwist1}
           \theta_V f=f\, \theta_U;\\
  & \label{cattwist2}
           (\theta_U\otimes \id_{U^*})\coev_U =(\id_U \otimes
           \theta_{U^*}) \coev_{U};\\
  & \label{cattwist3}
          \theta_{U\otimes V} =c_{V,U}\,
          c_{U,V}\, (\theta_{U}\otimes \theta_V);
\end{align}
for all objects $U,V\in \mathcal{C}$ and all morphism $f:U\to V$ in $\mathcal{C}$. It follows from \eqref{cattwist3} that
$\theta_{\un}=\id_{\un}$.

A ribbon category $\mathcal{C}$ canonically has a right duality by associating to any object $U\in \mathcal{C}$ its left
dual $U^*\in \mathcal{C}$ and two morphisms $\tev_U: U \otimes U^* \to \un$ and $\tcoev_U:\un \to U^* \otimes U$ defined
by
\begin{align}
 & \label{catrightdual1} \tev_U=\ev_U c_{U,U^*}
       (\theta_U \otimes \id_{U^*}); \\
 & \label{catrightdual2} \tcoev_U=(\id_{U^*} \otimes \theta_U^{-1})
       (c_{U^*,U})^{-1}\coev_U.
\end{align}
Note that we have $\tev_\un=\id_\un$ and $\tcoev_\un=\id_\un$.

The \emph{dual} morphism $f^*:V^*\to U^*$ of a morphism $f:U\to V$ in a ribbon category $\cc$ is defined by
\begin{align}
 f^* & =(\ev_V\otimes \id_{U^*})(\id_{V^*}\otimes f\otimes \id_{U^*}) (\id_{V^*}\otimes \coev_U) \label{cattranspose} \\
     & =(\id_{U^*} \otimes \tev_V)(\id_{U^*}\otimes f\otimes \id_{V^*}) (\tcoev_U \otimes \id_{V^*}). \nonumber
\end{align}
It is well-known that $(\id_U)^*=\id_{U^*}$ and $(f g)^*=g^* f^*$ for composable morphisms $f,g$. Axiom \eqref{cattwist2}
can be shown to be equivalent to $\theta_{U}^*= \theta_{U^*}$.

Let $\mathcal{C}$ be a ribbon category. Note that $\Endo_\cc(\un)$ is a semigroup, with composition as multiplication,
which is commutative (since $\cc$ is monoidal). The \emph{quantum trace} of an endomorphism $f:U\to U$ of an object $U\in
\mathcal{C}$ is defined by
\begin{equation}\label{cattrace}
  \tr(f)= \tev_U (f\otimes \id_{U^*}) \,
        \coev_{U} =\ev_U (\id_{U^*}\otimes f) \,
        \tcoev_{U}\in \Endo_\mathcal{C}(\un).
\end{equation}
For any morphisms $u:U\to V$, $v:V\to U$ and for any endomorphisms $f,g$ in $\mathcal{C}$, we have
\begin{equation}\label{tracepptstens}
\tr(uv)=\tr(vu), \quad \tr(f^*)=\tr(f), \quad \text{and} \quad \tr(f\otimes g)=\tr(f)\, \tr(g).
\end{equation}
The \emph{quantum dimension} of an object $U\in \mathcal{C}$ is defined by
\begin{equation}\label{catdim}
 \qdim (U)= \tr (\id_U) = \tev_U \, \coev_{U}= \ev_U \, \tcoev_{U} \in \Endo_\mathcal{C}(\un).
\end{equation}
Isomorphic objects have equal dimensions and $\qdim(U\otimes V)=\qdim (U)\, \qdim (V)$ for any objects $U,V \in \cc$.
Note that $\qdim (\un)=\id_\un$.

\subsection{$\Bbbk$-categories} \label{ribkkcats}
Let $\kk$ be a field. By a \emph{\kt category}, we shall mean a category for which the sets of morphisms are \kt spaces
and the composition is \kt bilinear. By a \emph{monoidal \kt category}, we shall mean \kt category endowed with a
monoidal structure whose tensor product is \kt bilinear. Note that if $\cc$ is a monoidal \kt category, then
$\Endo_\cc(\un)$ is a commutative \kt algebra (with composition as multiplication). A monoidal \kt category is said to be
\emph{pure} if $\Endo_\cc(\un)=\kk$.

By a \emph{ribbon \kt category}, we shall mean a pure monoidal \kt category endowed with a ribbon structure.

\subsection{Graphical calculus}\label{graphcalkul}
Let $\cc$ be ribbon category. Any morphism in $\cc$ can be graphically represented by a plane diagram (we use the
conventions of \cite{Tur2}). This pictorial calculus will allow us to replace algebraic arguments involving commutative
diagrams by simple geometric reasoning. This is justified, e.g., in \cite{Tur2}.

\begin{figure}[t]
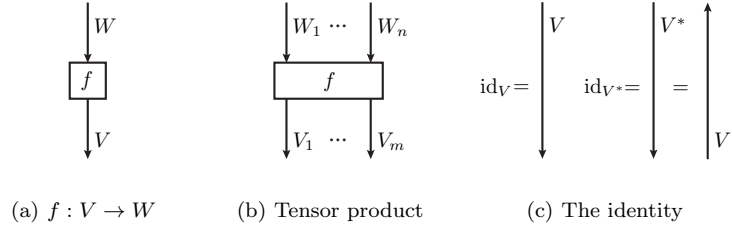

     \subfigure[$f:V \to W$]{\label{gc-mor}
                               \psfrag{f}[][]{$f$}
                               \psfrag{V}{$V$}
                               \psfrag{W}{$W$}
                               \scalebox{.8}{\includegraphics{gc-mor.eps}}} 
     \subfigure[Tensor product]{\label{gc-mortens}
                               \psfrag{f}{$f$}
                               \psfrag{V}{$V_1$}
                               \psfrag{W}{$W_1$}
                               \psfrag{R}{$V_m$}
                               \psfrag{T}{$W_n$}
                               \scalebox{.8}{\hspace*{.2cm}\includegraphics{gc-mortens.eps}}\hspace*{.4cm}}
     \subfigure[The identity]{\label{gc-id}
                               \psfrag{I}[Br][Br]{$\id_{V^*}$}
                               \psfrag{V}{$V$}
                               \psfrag{J}[Br][Br]{$\id_V$}
                               \psfrag{A}{$V^*$}
                               \psfrag{=}{$=$}
                               \scalebox{.8}{\includegraphics{gc-id.eps}}}
     \caption{Plane diagrams of morphisms}
     \label{gc-diagmor}
\end{figure}

A morphism $f:V \to W$ in $\cc$ is represented by a box with two vertical arrows oriented downwards, as in
Figure~\ref{gc-mor}. Here $V,W$ should be regarded as ``colors'' of the arrows and $f$ should be regarded as a ``color''
of the box. More generally, a morphism $f: V_1 \otimes  \cdots \otimes V_m \to W_1 \otimes \cdots \otimes W_n$ may be
represented as in Figure~\ref{gc-mortens}.

We also use vertical arrows oriented upwards under the convention that the morphism sitting in a box attached to such an
arrow involves not the color of the arrow but rather the dual object.

The identity endomorphism of an object $V \in \cc$  or of its dual $V^*$ will be represented by a vertical arrow as
depicted in Figure~\ref{gc-id}. Note that a vertical arrow colored with $\un$ may be deleted from any picture without
changing the morphism represented by this picture. The symbol ``$=$'' displayed in the figures denotes equality of the
corresponding morphisms in $\cc$.

The tensor product $f \otimes g$ of two morphisms $f$ and $g$ in $\cc$ is represented by placing a picture of $f$ to the
left of a picture of $g$. A picture for the composition $g \circ f$ of two (composable) morphisms $g$ and $f$ is obtained
by putting a picture of $g$ on the top of a picture of $f$ and by gluing the corresponding free ends of arrows.

The braiding $c_{V,W}: V \otimes W \to W \otimes V$ and its inverse $c_{V,W}^{-1}: W \otimes V \to V \otimes W$, the
twist $\theta_V: V \to V$ and its inverse $\theta_V^{-1}:V \to V$, and the duality morphisms $\ev_V: V^* \otimes V \to
\un$, $\coev_V: \un \to V \otimes V^*$, $\tev_V: V \otimes V^* \to \un$, and $\tcoev_V: \un \to V^* \otimes V$ are
represented as in Figures~\ref{gc-braid}, \ref{gc-twist}, and \ref{gc-dual} respectively. The quantum trace of an
endomorphism $f : V \to V$ in $\cc$ and the quantum dimension of an object $V \in \cc$ may be depicted as in
Figure~\ref{gc-trace}.
\begin{figure}[h]
        \subfigure[Braiding]{ \label{gc-braid}
                               \psfrag{c}[][]{$c_{V,W}$}
                               \psfrag{V}{$V$}
                               \psfrag{S}[Br][Br]{$V$}
                               \psfrag{T}[Br][Br]{$W$}
                               \psfrag{W}{$W$}
                               \psfrag{R}{$W$}
                               \psfrag{a}[][]{$c_{V,W}^{-1}$}
                               \psfrag{=}[][]{$=$}
                               \scalebox{.8}{\includegraphics{gc-braid.eps}}}
        \subfigure[Twist]{ \label{gc-twist}
                               \psfrag{c}[Br][Br]{$\theta_V^{-1}$}
                               \psfrag{V}[Br][Br]{$V$}
                               \psfrag{W}[Br][Br]{$V$}
                               \psfrag{a}[Br][Br]{$\theta_V$}
                               \psfrag{=}[B][B]{$=$}
                               \scalebox{.8}{\includegraphics{gc-twist.eps}}}
        \subfigure[Duality morphisms]{ \label{gc-dual}
                               \psfrag{c}[Br][Br]{$\tev_V$}
                               \psfrag{V}{$V$}
                               \psfrag{e}[Br][Br]{$\ev_V$}
                               \psfrag{u}[Br][Br]{$\tcoev_V$}
                               \psfrag{a}[Br][Br]{$\coev_V$}
                               \psfrag{=}[B][B]{$=$}
                               \scalebox{.8}{\includegraphics{gc-dual.eps}}}
        \subfigure[Trace and dimension]{ \label{gc-trace}
                               \psfrag{V}{$V$}
                               \psfrag{c}[Br][Br]{$\qdim V$}
                               \psfrag{a}[Br][Br]{$\tr(f)$}
                               \psfrag{f}[cr][cr]{$f$}
                               \psfrag{=}[B][B]{$=$}
                               \scalebox{.8}{\includegraphics{gc-trace.eps}}}
        \caption{}
        \label{appareilmoustic}
\end{figure}

\subsection{Negligible morphisms}\label{sect-negl}
Let $\cc$ be a ribbon \kt category. A morphism $f \in \Homo_\cc(X,Y)$ is said to be \emph{negligible} if $\tr(gf)=0$ for
all $g \in \Homo_\cc(Y,X)$. Denote by $\Negl_\cc(X,Y)$ the \kt subspace of $\Homo_\cc(X,Y)$ formed by the negligible
morphisms.

It is important to note that $\Negl_\cc$ is a two-sided $\otimes$-ideal of $\cc$. This means that the composition
(resp.\@ the tensor product) of a negligible morphism with any other morphism is negligible.

Note that since $\Endo(\un)=\kk$, a morphism $f: \un \to X$ is negligible if and only if $gf=0$ for all morphism $g: X
\to \un$.

\subsection{Dinatural transformations and coends}\label{sect-coend}
Recall that to each category $\cc$ we associate the \emph{opposite category} $\cc^\opp$ in the following way: the objects
of $\cc^\opp$ are the objects of $\cc$ and the morphisms of $\cc^\opp$ are morphisms $f^\opp$, in one-one correspondence
$f \mapsto f^\opp$ with the morphisms in $\cc$. For each morphism $f:U \to V$ of $\cc$, the domain and codomain of the
corresponding $f^\opp$ are as in $f^\opp: V \to U$ (the direction is reversed). The composite $f^\opp g^\opp=(gf)^\opp$
is defined in $\cc^\opp$ exactly when the composite $gf$ is defined in $\cc$. This makes $\cc^\opp$ a category.

Let $\cc$ and $\bb$ be two categories. A \emph{dinatural transformation} between a functor $F: \cc^\opp \times \cc \to
\bb$ and an object $B \in \bb$ is a function $d$ which assigns to each object $X \in \cc$ a morphism $d_X : F(X,X) \to B$
of $\bb$ in such a way that the diagram
\begin{equation*}
\begin{CD}
F(Y,X) @>{F(\id_Y,f)}>> F(Y,Y) \\
@V{F(f,\id_X)}VV @VV{d_{Y}}V \\
F(X,X) @>{d_X}>> B
\end{CD}
\end{equation*}
commutes for every morphism $f:X \to Y$ in $\cc$.

A \emph{coend} of the functor $F$ is a pair $( A, i )$ consisting of an object $A$ of $\bb$ and a dinatural
transformation $i$ from $F$ to $A$ which is \emph{universal} among the dinatural transformation from $F$ to a constant,
that is, with the property that, to every dinatural transformation $d$ from $F$ to $B$, there exists a unique morphism
$r:A \to B$ such that, for all object $X \in \cc$,
\begin{equation}\label{coendfacto}
  d_X=r \circ i_X.
\end{equation}

By using the factorization property~\eqref{coendfacto}, it is easy to verify that if $ (A, i)$ and $(A', i')$  are two
coends of $F$, then they are isomorphic in the sense that there exists an isomorphism $I:A \to A'$ in $\bb$ such that
$i'_X=I \circ i_X$ for all object $X \in \cc$.

\subsection{Categorical Hopf algebras from coends}\label{coendHA}

Let $\cc$ be a ribbon category. Consider the functor $F: \cc^\opp \times \cc \to \cc$ defined by
\begin{equation}\label{Ffunctor}
  F(X,Y)=X^* \otimes Y \text{\quad  and \quad }
  F(f^\opp,g)=f^* \otimes g
\end{equation}
for all objects $X \in \cc^\opp$, $Y \in \cc$ and all morphisms $f,g$ in $\cc$.

\begin{figure}[h]
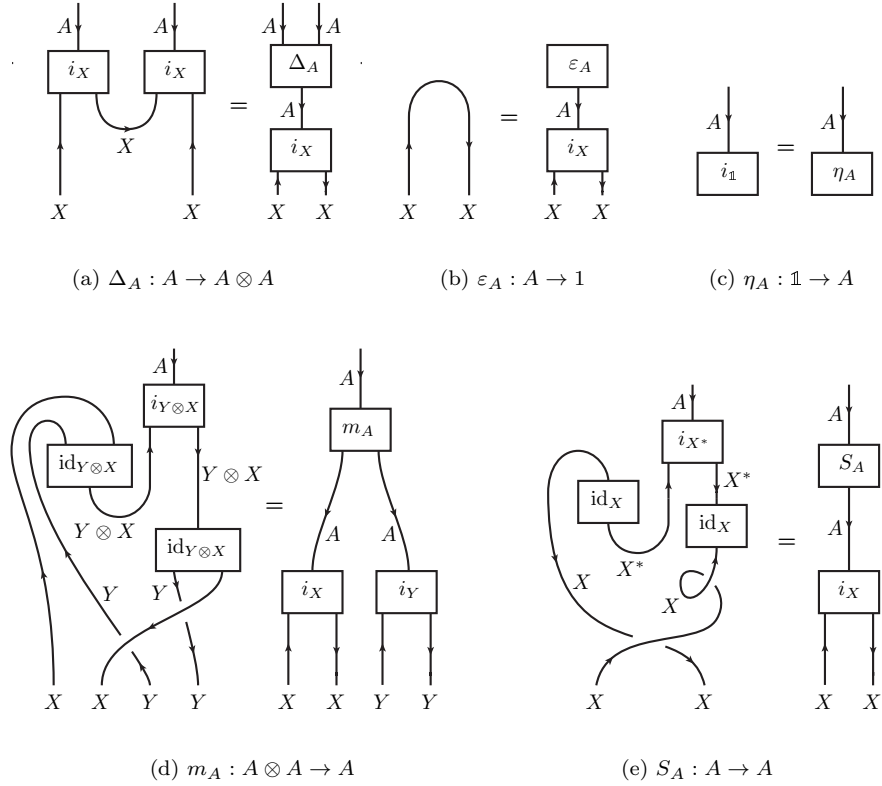

   \begin{center}
     \subfigure[$\Delta_A : A \to A \otimes A$]{
       \psfrag{A}[Br][Br]{$A$}
       \psfrag{B}{$A$}
       \psfrag{X}[B][B]{$X$}
       \psfrag{E}[B][B]{$i_X$}
       \psfrag{H}[B][B]{$\Delta_A$}
       \psfrag{=}{\textbf{=}}
       \scalebox{.8}{\includegraphics{coend-th1.eps}}}
       \subfigure[$\varepsilon_A: A \to 1$]{
       \psfrag{A}[Br][Br]{$A$}
       \psfrag{X}[B][B]{$X$}
       \psfrag{E}[B][B]{$i_X$}
       \psfrag{H}[B][B]{$\varepsilon_A$}
       \psfrag{=}{\textbf{=}}
       \scalebox{.8}{\includegraphics{coend-th1c.eps}} \qquad}
       \subfigure[$\eta_A: \un \to A$]{
       \psfrag{A}[Br][Br]{$A$}
       \psfrag{B}{$A$}
       \psfrag{X}[B][B]{$X$}
       \psfrag{E}[B][B]{$i_\un$}
       \psfrag{H}[B][B]{$\eta_A$}
       \psfrag{=}{\textbf{=}}
       \scalebox{.8}{\;\includegraphics{coend-th1u.eps}\;}}\\
     \subfigure[$m_A:A \otimes A \to A$]{
      \psfrag{A}[Br][Br]{$A$}
      \psfrag{B}[Bl][Bl]{$A$}
      \psfrag{I}[B][B]{$\id_{Y \otimes X}$}
      \psfrag{X}[B][B]{$X$}
      \psfrag{Y}[B][B]{$Y$}
      \psfrag{U}[Br][Br]{$Y$}
      \psfrag{=}{\textbf{=}}
      \psfrag{Z}{$Y$}
      \psfrag{V}[Br][Br]{$Y \otimes X$}
      \psfrag{T}{$Y \otimes X$}
      \psfrag{E}[B][B]{$i_{Y \otimes X}$}
      \psfrag{S}[B][B]{$i_X$}
      \psfrag{H}[B][B]{$i_Y$}
      \psfrag{m}[B][B]{$m_A$}
      \scalebox{.8}{\includegraphics{coend-th2.eps}} \qquad}
     \subfigure[$S_A: A \to A$]{
      \psfrag{A}[Br][Br]{$A$}
      \psfrag{I}[B][B]{$\id_{X}$}
      \psfrag{E}[B][B]{$i_{X^*}$}
      \psfrag{S}[B][B]{$i_X$}
      \psfrag{H}[B][B]{$S_A$}
      \psfrag{X}[B][B]{$X$}
      \psfrag{U}[Br][Br]{$X$}
      \psfrag{Z}{$X$}
      \psfrag{V}[Br][Br]{$X^*$}
      \psfrag{=}{\textbf{=}}
      \psfrag{T}{$X^*$}
      \quad \scalebox{.8}{\includegraphics{coend-th3.eps}}}
     \end{center}
     \caption{Structural morphisms of $A$}
     \label{figcathopfalg}
\end{figure}

Suppose that the functor $F$ admits a coend $( A, i)$. Then the object $A$ has a structure of a Hopf algebra in the
category $\cc$ (see \cite{Lyu1}). This means that there exist morphisms $m_A : A \otimes A \to \un$, $\eta_A: \un \to A$,
$\Delta_A: A \to A \otimes A$, $\varepsilon_A : A \to \un$, and $S_A : A \to A$, which verify the same axioms as those of
a Hopf algebra except the usual flip is replaced by the braiding $c_{A,A}: A \otimes A \to A \otimes A$. By using the
factorization property \eqref{coendfacto} and the naturality of the braiding and twist of $\cc$, these structural
morphisms are defined as follows:
\begin{align}
 &\Delta_A \circ i_X = \Bigl ( X^* \otimes X  \rTo^{\id_{X^*} \otimes \coev_X \otimes \id_X}
     X^* \otimes X \otimes X^* \otimes X
 \rTo^{i_X \otimes i_X} A \otimes A \Bigr ), \label{defrel1}\\
 &\varepsilon_A \circ i_X=\ev_X,\\
 &\eta_A: \un = \un^* \otimes \un \rTo^{i_\un} A,\\
 &m_A \circ (i_X \otimes i_Y) = \Bigl ( X^* \otimes X \otimes Y^* \otimes Y \rTo^{c_{X,Y^* \otimes Y}}
 X^* \otimes Y^* \otimes Y \otimes X\\
 & \phantom{XXXxxxxxxxXXXXXXXxxxXX} \simeq (Y \otimes X)^* \otimes (Y \otimes X) \rTo^{i_{Y \otimes X}} A \Bigr ),\nonumber \\
 &S_A \circ i_X = \Bigl ( X^* \otimes X \rTo^{c_{X^*,X}} X \otimes X^* \rTo^{\id_X \otimes \theta_{X^*}}
   X \otimes X^* \label{defrel2}\\
 & \phantom{XXXXXXXXXXXXxxxxxxxx;xxxxxxxxxx} \simeq  X^{**} \otimes X^* \rTo^{i_{X^*}} A \Bigr ). \nonumber
\end{align}
Here $X$ and $Y$ are any objects of $\cc$. The defining relations \eqref{defrel1}-\eqref{defrel2} are graphically
represented in Figure~\ref{figcathopfalg}.

It can be shown that the antipode $S_A:A \to A$ is an isomorphism and that $S_A^2=\theta_A$ (see \cite{Lyu1}). Moreover,
as for Hopf algebras, the antipode is anti-(co)multiplicative, that is,
\begin{align}
  & S_A m_A = m_A (S_A \otimes S_A) c_{A,A}, & & S_A \eta_A=\eta_A, \label{podantimult} \\
  & \Delta_A S_A = c_{A,A} (S_A \otimes S_A) \Delta_A, & & \varepsilon_A S_A=\varepsilon_A.\label{podanticomult}
\end{align}

The Hopf algebra $A$ posses a Hopf pairing $\omega_A:A \otimes A \to \un$ (see \cite{Lyu2}) defined by
\begin{equation}\label{parBB}
 \omega_A \circ (i_X \otimes i_Y)=(\ev_X \otimes \ev_Y)(\id_{X^*} \otimes c_{Y^*,X}c_{X,Y^*} \otimes \id_Y)
\end{equation}
for any objects $X,Y \in \cc$. This pairing is said to be \emph{non-degenerate} if $\omega_A(\id_A \otimes \coev_A) :A
\to A^*$ and
$\omega_A(\tcoev_A \otimes \id_A) :A \to A^*$ are isomorphisms.\\

Set:
\begin{align}
 & \Gamma_l=(\id_A \otimes m_A)(\Delta_A \otimes \id_A):A \otimes A \to A \otimes A,\label{defgaml}\\
 & \Gamma_r=(m_A \otimes \id_A)(\id_A \otimes \Delta_A):A \otimes A \to A \otimes A. \label{defgamr}
\end{align}

\begin{lem}\label{lemgarl}
$\Gamma_l(S_A \otimes S_A) c_{A,A}=c_{A,A}(S_A \otimes S_A) \Gamma_r$.
\end{lem}

\begin{proof}
By using \eqref{catbraid1}, \eqref{podantimult} and \eqref{podanticomult}, we have
\begin{align*}
\Gamma_l(S_A \otimes S_A) c_{A,A}
  & = (\id_A \otimes m_A)(\Delta_A S_A \otimes S_A) c_{A,A} \\
  & = (\id_A \otimes m_A)(c_{A,A} (S_A \otimes S_A) \Delta_A \otimes S_A) c_{A,A} \\
  & = (S_A \otimes m_A(S_A \otimes S_A))(c_{A,A}\Delta_A \otimes \id_A) c_{A,A} \\
  & = (S_A \otimes S_A m_A)(\id_A \otimes c_{A,A}^{-1})(c_{A,A} \otimes \id_A)(\Delta_A \otimes \id_A) c_{A,A}.
\end{align*}
Then, using \eqref{catbraid1} and \eqref{catbraid2}-\eqref{catbraid3}, we get that
\begin{eqnarray*}
\lefteqn{\Gamma_l(S_A \otimes S_A) c_{A,A}}\\
  &  & =(S_A \otimes S_A m_A)(\id_A \otimes c_{A,A}^{-1})(c_{A,A} \otimes \id_A) c_{A,A\otimes A} (\id_A \otimes \Delta_A)\\
  &  & =(S_A \otimes S_A m_A)(\id_A \otimes c_{A,A}^{-1})c_{A,A\otimes A}(\id_A \otimes c_{A,A})(\id_A \otimes \Delta_A)\\
  &  & =(S_A \otimes S_A m_A)(\id_A \otimes c_{A,A}^{-1})(\id_A \otimes c_{A,A}) (c_{A,A} \otimes \id_A)(\id_A \otimes c_{A,A})
        (\id_A \otimes \Delta_A)\\
  &  & =(S_A \otimes S_A m_A)c_{A\otimes A,A}(\id_A \otimes \Delta_A)\\
  &  & =c_{A,A}(S_A \otimes S_A)(m_A \otimes \id_A)(\id_A \otimes \Delta_A)\\
  &  & =c_{A,A}(S_A \otimes S_A)\Gamma_r.\\
\end{eqnarray*}
\end{proof}

\begin{cor}\label{corgamm}
Suppose that $\cc$ is moreover a (monoidal) \kt category. Let $\al \in \Homo_\cc(\un,A)$. If $S_A \al - \al  \in
\Negl_\cc(\un,A)$, then the following assertions are equivalent:
\begin{enumerate}
  \renewcommand{\labelenumi}{{\rm (\alph{enumi})}}
  \item $\Gamma_l(\al \otimes \al)-\al \otimes \al : \un \to A\otimes A$ is negligible;
  \item $\Gamma_r(\al \otimes \al)-\al \otimes \al:\un \to A\otimes A$ is negligible.
\end{enumerate}
Moreover, if  $S_A \circ \al=\al$, then $\Gamma_l(\al \otimes \al)=\al \otimes \al$ if and only if $\Gamma_r(\al \otimes
\al)=\al \otimes \al$.
\end{cor}

\begin{proof}
This is an immediate consequence of Lemma~\ref{lemgarl} since $\Negl_\cc$ is a two-sided $\otimes$-ideal of $\cc$.
\end{proof}

\section{Kirby elements of a ribbon category}\label{sect-2}
In this section, we generalize the Lyubashenko's method \cite{Lyu2} of constructing 3-manifolds invariants from ribbon
categories.

\subsection{Special tangles}
Let $n$ be a positive integer. By a \emph{$n$-special tangle} we shall mean an oriented framed tangle $T \subset
\mathbb{R}^2 \times [0,1]$ with $2n$ bottom endpoints and no top endpoints, consisting of $n$ arc components such that
the $k$th arc ($1\leq k \leq n$) joins the $(2k-1)$th and $(2k)$th bottom endpoints and is oriented out of $\mathbb{R}^2
\times [0,1]$ near the $(2k)$th bottom endpoint (and so inside $\mathbb{R}^2 \times [0,1]$ near the $(2k-1)$th bottom
endpoint).

Diagrams of special tangles are drawn with blackboard framing. An example of a $3$-special tangle is depicted in
Figure~\ref{exspectang1}.

Let $\cc$ be a ribbon category. Suppose that the functor \eqref{Ffunctor} admits a coend $( A, i)$. Let $T$ be a
$n$-special tangle. For objets $X_1, \dots , X_n \in \cc$, let $T_{(X_1, \cdots , X_n)}$ be the morphism $X_1^*\otimes
X_1 \otimes \cdots \otimes X_n^* \otimes X_n \to \un$ in $\cc$ graphically represented by a diagram of $T$ where the
$k$th component of $T$ has been colored with the object $X_k$. Note that $T_{(X_1, \cdots , X_n)}$ does not depend on the
choice of the diagram of $T$, that is, only depends on the isotopy class of $T$ (see \cite{Tur2}). Since the braiding and
twist of $\cc$ are natural and by using the Fubini theorem for coends (see \cite{ML1}), there exists a (unique) morphism
$\phi_T: A^{\otimes n} \to \un$ such that
\begin{equation}\label{phiTL}
T_{(X_1, \dots,X_n)}=\phi_T \circ (i_{X_1} \otimes \cdots \otimes i_{X_n})
\end{equation}
for all objects $X_1, \dots , X_n\in \cc$ (see Figure~\ref{exspectang2} for $n=3$).

\begin{figure}[h]
   \begin{center}
     \subfigure[A 3-special tangle $T$]{\label{exspectang1}
       \scalebox{.8}{\;\, \includegraphics{coend-special1.eps}\;\,}
       } \quad
       \subfigure[$T_{(X_1,X_2,X_3)}=\phi_T \circ (i_{X_1} \otimes i_{X_2} \otimes i_{X_2})$]{\label{exspectang2}
       \psfrag{a}[B][B]{$i_{X_1}$}
       \psfrag{u}[B][B]{$i_{X_2}$}
       \psfrag{c}[B][B]{$i_{X_3}$}
       \psfrag{X}[B][B]{$X_1$}
       \psfrag{A}[B][B]{$A$}
       \psfrag{Y}[B][B]{$X_2$}
       \psfrag{T}[B][B]{$X_3$}
       \psfrag{E}[B][B]{$T_{(X_1,X_2,X_3)}$}
       \psfrag{H}[B][B]{$\phi_T$}
       \psfrag{=}{\textbf{=}}
       \scalebox{.8}{\includegraphics{coend-special2.eps}}
       }
     \end{center}
     \caption{}
\end{figure}

\subsection{Kirby elements}
Let $\cc$ be a ribbon category such that the functor \eqref{Ffunctor} admits a coend $( A, i)$.

Let $L$ be a framed link in $S^3$ with $n$ components. Fix an orientation for $L$. There always exists a (non-unique)
$n$-special tangle $T_L$ such that $L$ is isotopic $T_L \circ (\cup_- \otimes \cdots \otimes \cup_-)$, where $\cup_-$
denote the cup with clockwise orientation, see Figure~\ref{figdem3man3}. For $\alpha \in \Homo_\cc(\un,A)$, set
\begin{equation*}
\tau_\cc(L;\alpha)=\phi_{T_L} \circ \alpha^{\otimes n} \in \Endo_\cc(\un),
\end{equation*}
where $\phi_{T_L}:A^{\otimes n} \to \un$ is defined as in \eqref{phiTL}.

\begin{defn}
By a \emph{Kirby element} of $\cc$, we shall mean a morphism $\alpha \in \Homo_\cc(\un,A)$ such that, for any framed link
$L$, $\tau_\cc(L;\alpha)$ is well-defined and invariant under isotopies and 2-handle slides of $L$. A Kirby element $\al$
of $\cc$ is said to be \emph{normalized} if $\tau_\cc(\bigcirc^{\pm 1};\alpha)$ is invertible in $\Endo_\cc(\un)$, where
$\bigcirc^{\pm 1}$ denotes the unknot with framing~$\pm 1$.
\end{defn}

Note that the unit  $\eta_A:\un \to A$ of the categorical Hopf algebra~$A$ is a normalized Kirby element. The invariant
of framed links associated $\eta_A$ is the trivial one, that is, $\tau_\cc( L;\eta_A)=1$ for any framed link $L$.

\begin{lem}\label{linkdisjunion}
Let $\al$ be a Kirby element of $\cc$. Then $\tau_\cc(L \sqcup L';\al)=\tau_\cc(L;\al)\,\tau_\cc(L';\al)$ for any framed
link $L$ and $L'$, where $L \sqcup L'$ denotes the disjoint union of $L$ and $L'$.
\end{lem}
\begin{proof}
Let $T_L$ (resp.\@ $T_{L'}$) be a special tangle such that $L$ (resp.\@ $L'$) is isotopic $T_L \circ (\cup_- \otimes
\cdots \otimes \cup_-)$ (resp.\@ $T_{L'}\circ (\cup_- \otimes \cdots \otimes \cup_-)$). The tangle $T=T_L \otimes T_{L'}$
is then special and such that the disjoint union $L \sqcup L'$ is isotopic $T \circ (\cup_- \otimes \cdots \otimes
\cup_-)$. Therefore $\phi_T=\phi_{T_L} \otimes \phi_{T_{L'}}$ and so $\tau_\cc(L \sqcup
L';\al)=\tau_\cc(L;\al)\,\tau_\cc(L';\al)$.
\end{proof}

Let $\Theta_\pm : A \to \un$ be the morphism defined by
\begin{equation}\label{tethapmdef}
\Theta_\pm \circ i_X=\ev_X(\id_{X^*} \otimes \theta_X^{\pm 1}),
\end{equation}
where $X$ is any object of $\cc$ (see Figure~5). Note that if $\al$ be a Kirby element of $\cc$, then
$\tau_\cc(\bigcirc^{\pm 1};\al)=\Theta_\pm \al$.
\begin{figure}[h]\label{thetapm}
   \begin{center}
       \psfrag{E}[B][B]{$i_X$}
       \psfrag{X}[Bl][Bl]{$X$}
       \psfrag{T}[B][B]{$\Theta_-$}
       \psfrag{H}[B][B]{$\Theta_+$}
       \psfrag{=}[Bl][Bl]{\textbf{=}}
       \psfrag{A}[Br][Br]{$A$}
       \scalebox{.8}{\includegraphics{thetaplus.eps}}
   \end{center}
     \caption{}
\end{figure}

Recall (see \cite{Li}) that every closed, connected,  and oriented 3-manifold can be obtained from $S^3$ by surgery along
a framed link $L \subset S^3$. In this paper, all considered 3-manifolds are supposed to be closed, connected,  and
oriented. For any framed link $L$ in $S^3$, we will denote by $M_L$ the 3-manifold obtained from $S^3$ by surgery along
$L$, by $ n_L$ the number of components of $L$, and by $b_-(L)$ the number of negative eigenvalues of the linking matrix
of $L$.

Normalized Kirby elements are of special interest due to the following proposition.

\begin{prop}\label{prop3man}
Let $\al$ be a normalized Kirby element of $\cc$. Then
\begin{equation*}
\tau_\cc(M_L;\al)=(\Theta_+ \al)^{b_-(L)-n_L}\,  (\Theta_- \al)^{-b_-(L)} \; \tau_\cc(L;\alpha)
\end{equation*}
is an invariant of 3-manifolds. Moreover $\tau_\cc(M\# M';\al)=\tau_\cc(M;\al)\,\tau_\cc(M';\al)$ for any 3-manifold $M$
and $M'$.
\end{prop}
\begin{proof}
The fact that $\tau_\cc(M_L;\al)$ is an invariant of 3-manifolds follows from Kirby theorem \cite{Ki}. Indeed
$\tau_\cc(L;\al)$, $b_-(L)$ and $n_L$ are invariant under 2-handle slides and  $\tau_\cc(\bigcirc^{\pm 1} \sqcup
L;\alpha)=(\Theta_\pm \al) \, \tau_\cc(L;\alpha)$ by Lemma~\ref{linkdisjunion}, $b_- (\bigcirc^1 \sqcup L)= b_- (L)$,
$b_- (\bigcirc^{-1} \sqcup L)= b_- (L) + 1$, and $n_{\bigcirc^{\pm 1} \sqcup L}= n_L+1$.

Let $L$ and $L'$ be framed links in $S^3$. The disjoint union $L \sqcup L'$ is then a framed link in $S^3$ such that
$M_{L \sqcup L'}\simeq M_L\# M_{L'} $.  Then the multiplicativity of $\tau_\cc(M;\al)$ with respect to the connected sum
of 3-manifolds follows from Lemma~\ref{linkdisjunion} and the equalities $b_-(L \sqcup L')=b_-(L)+b_-(L')$ and $n_{L
\sqcup L'}=n_L+n_{L'}$.
\end{proof}

\begin{remsbold}\
{\bf [a].} For any normalized Kirby element $\al$ of $\cc$, we have that
\begin{equation*}
\tau_\cc(S^3;\al)=1 \quad \text{and} \quad \tau_\cc(S^1 \times S^2;\al)=(\Theta_+ \al)^{-1} \, \varepsilon_A\al.
\end{equation*}

{\bf [b].} The invariant of 3-manifolds associated to the unit  $\eta_A:\un \to A$ of the categorical Hopf algebra~$A$
(which is a normalized Kirby element) is the trivial one, that is, $\tau_\cc( M ;\eta_A)=1$ for any 3-manifold $M$.
\end{remsbold}

In general, determining when a morphism $A^{\otimes n} \to \un$ is of the form $\phi_T$ for some special $n$-tangle $T$
is a quite difficult problem. Hence does so the problem of determining all the (normalized) Kirby elements of $\cc$. In
the next section, we characterize a class of (normalized) Kirby elements of $\cc$ by means of the structural morphisms of
the categorical Hopf algebra $A$. This class will be shown to be sufficiently large to contain the Lyubashenko invariant
(which is a categorical version of the Hennings-Kauffman-Radford invariant) and the Reshetikhin-Turaev invariant
(computed from a semisimple quotient of $\cc$) when these are well-defined.

\subsection{A class of Kirby elements}
Let $\cc$ be a ribbon \kt category such that the functor \eqref{Ffunctor} admits a coend $(A,i)$.

Recall the notion of negligible morphisms (see Section~\ref{sect-negl}). Set:
\begin{align*}
   &\Ki(\cc)=\{ \al \in \Homo_\cc(\un,A) \; | \; S_A\al - \al  \in \Negl_\cc(\un,A) \quad \text{and}\\
   & \phantom{\Ki(\cc)=\{ \al \in \Homo_\cc(\un,A) \; | \;} \Gamma_l(\al \otimes \al)-\al \otimes \al
     \in \Negl_\cc(\un,A \otimes A) \},\\
   &\Ki(\cc)^\norm =\{ \al \in \Ki(\cc) \; | \; \Theta_+ \al \neq 0 \quad \text{and} \quad \Theta_- \al \neq 0 \},
\end{align*}
where $\Gamma_l:A \otimes A \to A \otimes A$ and $\Theta_\pm:A \to \un$ are defined in \eqref{defgaml} and
\eqref{tethapmdef}.

Note that, in the above definitions, the morphism  $\Gamma_l$ can be replaced by the morphism~$\Gamma_r$ which is defined
in~\eqref{defgamr}  (see Corollary~\ref{corgamm}).

Remark that the sets $\Ki(\cc)$ and $\Ki(\cc)^\norm$ always contain a non-zero element, namely the unit $\eta_A:\un \to
A$. Moreover,
\begin{equation}\label{plusnegl}
\kk \Ki(\cc) + \Negl_\cc(\un,A) \subset \Ki(\cc) \text{ and } \kk^*\Ki(\cc)^\norm + \Negl_\cc(\un,A) \subset
\Ki(\cc)^\norm,
\end{equation}
since $\Negl_\cc$ is a two-sided $\otimes$-ideal of $\cc$ and $\Negl_\cc(\un,\un)=0$ (because $\Endo_\cc(\un)=\kk$ is a
field).

\begin{thm}\label{thm3man}
\begin{enumerate}
 \renewcommand{\labelenumi}{{\rm (\alph{enumi})}}
 \item The elements of $\Ki(\cc)$ are Kirby elements of $\cc$.
 \item The elements of $\Ki(\cc)^\norm$ are normalized Kirby elements of $\cc$.
\end{enumerate}
\end{thm}

\begin{remsbold}\label{rmsdeKi}
{\bf [a].} It follows from Proposition~\ref{prop3man} that any $\al \in \Ki(\cc)^\norm$ leads to a 3-manifolds invariant
$\tau_\cc(M;\al)$ with values in $\Endo_\cc(\un)=\kk$. Recall that this invariant is multiplicative with respect to the
connected sum of 3-manifolds and verifies
\begin{equation*}
\tau_\cc(S^3;\al)=1, \quad \tau_\cc(S^1 \times S^2;\al)=(\Theta_+ \al)^{-1} \, \varepsilon_A\al, \quad
\tau_\cc(M;\eta_A)=1,
\end{equation*}
for any 3-manifold $M$.\\[-.8em]

{\bf [b].} Let $\al \in \Ki(\cc)^\norm$, $n \in \Negl_\cc(\un,A)$, and $k \in \kk^*$. Then
\begin{equation*}
k \al +n \in \Ki(\cc)^\norm \quad \text{and} \quad \tau_\cc(M;k\al+n)=\tau_\cc(M;\al)
\end{equation*}
for any 3-manifold $M$. This follows from \eqref{plusnegl} and the choice of the normalization in the definition of
$\tau_\cc(M;\al)$ (see Proposition~\ref{prop3man}).\\[-.8em]

{\bf [c].} Suppose that $\cc$ is 3-modular (see \cite{Lyu2}). Then $A$ admits a (non-zero) two-sided integral $\lambda
\in \Homo_\cc(\un,A)$, that is,
\begin{equation*}
m_A(\lambda \otimes \id_A)=\eta_A \varepsilon_A= m_A(\id_A \otimes \lambda).
\end{equation*}
Moreover $\lambda \in \Ki(\cc)^\norm$ and $\tau_\cc(M;\lambda)$ is the Lyubashenko invariant of 3-manifolds. Note that
when $\cc$ admits split idempotents, then $\lambda$ is unique (up to scalar multiple), see \cite{BKLT}.

The condition that $A$ posses a (non-zero) two-sided integral is quite limitative (for example, when $\cc$ is the
category $\reph$ of representations of a finite-dimensional Hopf algebra $H$, this implies that $H$ must be unimodular).
In Section~\ref{Sect-examples}, we give an example of non-unimodular ribbon Hopf algebra $H$ and of an element $\al \in
\Ki(\reph)$ which is not a two-sided integral and leads to a non-trivial invariant.\\[-.8em]

{\bf [d].} When $\cc$ is finitely semisimple, we show in Section \ref{sectinvarRT} (see Corollary~\ref{corRt}) that there
exists $\al_\cc \in \Ki(\cc)$ corresponding to the Reshetikhin-Turaev invariant defined with $\cc$. Note that $\al_\cc
\in \Ki(\cc)$ is not in general a two-sided integral.

More generally, we show in Section~\ref{sect-lifromsscats} (see Corollary~\ref{invfromsscat}) that $\Ki(\cc)$ contains
elements corresponding to the Reshetikhin-Turaev invariants defined with finitely semisimple full ribbon subcategories of
the ribbon semisimple quotient of $\cc$.
\end{remsbold}

\begin{proof}[Proof of Theorem~\ref{thm3man}]
Let us prove Part (a). Fix $\al \in \Ki(\cc)$. Let $L=L_1 \cup \dots \cup L_n$ be a framed link. Firstly, since $S_A
\alpha-\alpha \in \Negl_\cc(\un,A)$, then $\tau_\cc(L;\alpha)$ does not depend on choice of $T_L$ neither the orientation
of $L$ and is an isotopic invariant of the framed link $L$. Indeed, this is proved in the case $S_A \alpha=\alpha$ in
\cite[Proposition 5.2.1]{Lyu2}. The same arguments work when $S_A \alpha-\alpha \in \Negl_\cc(\un,A)$ since $\Negl_\cc$
is a two-sided $\otimes$-ideal of~$\cc$.

Let us show that $\tau_\cc(L;\alpha)$ is invariant under 2-handle slides. Choose an orientation for $L$. Without loss of
generality, we can suppose that the component $L_1$ slides over $L_2$. Let $L'_2$ be a copy of $L_2$ (following the
framing) and set $L'=(L_1 \# L'_2) \cup L_2 \cup \dots \cup L_n$. We have to show that
$\tau_\cc(L';\alpha)=\tau_\cc(L;\alpha)$. Let $T_L$ be a $n$-special tangle such that $L$ is isotopic $T_L \circ (\cup_-
\otimes \cdots \otimes \cup_-)$,  where the $i$th cup (with clockwise orientation) corresponds to the component $L_i$,
see Figure~\ref{figdem3man3}. Let $\Delta_2(T_L)$ be the $(2n+2,0)$-tangle obtained by copying the 2nd component of $T_L$
(following the framing) in such a way that the endpoints of the new component are between the 2nd and 3th bottom
endpoints and between the $4$th and $5$th bottom endpoints of $T_L$. A $n$-special tangle $T'$ such that $L'$ is isotopic
$T'\circ (\cup_- \otimes \cdots \otimes \cup_-)$ (where the $i$th cup corresponds to the $i$th component of $L'$) can be
constructed from $\Delta_2(T_L)$ as in Figure~\ref{figdem3man4}. For example, if $T_L$ is the 3-special tangle depicted
in Figure~\ref{figdem3man0}, then let $T'$ be the 3-special tangle of Figure~\ref{figdem3man1}.

\begin{figure}[h]
   \begin{center}
   \subfigure[$L \sim T_L \circ (\cup_- \otimes \cdots \otimes \cup_-)$]{\label{figdem3man3}
       \psfrag{E}[B][B]{$T_L$}
       \psfrag{L}[B][B]{$L\sim$}
       \psfrag{X}[Bl][Bl]{$L_1$}
       \psfrag{A}[Bl][Bl]{$L_2$}
       \psfrag{B}[Bl][Bl]{$L_n$}
       \psfrag{-}[B][B]{$\dots$}
       \scalebox{.8}{\; \includegraphics{dem3man3.eps} \quad \;}
       } \qquad
   \subfigure[$T'$]{\label{figdem3man4}
       \psfrag{-}[B][B]{$\dots$}
       \psfrag{I}[B][B]{$\Delta_2(T_L)$}
       \scalebox{.8}{\includegraphics{dem3man4.eps}}
       }\\ \vspace*{.4cm}
   \subfigure[$T_L$]{\label{figdem3man0}
       \scalebox{.8}{\includegraphics{dem3man0.eps}}
       }\qquad \qquad \qquad
   \subfigure[$T'$]{\label{figdem3man1}
       \scalebox{.8}{\includegraphics{dem3man1.eps}}
       }
   \end{center}
     \caption{}
\end{figure}

\noindent By the equalities of Figure~\ref{figdem3man2} where $X_1, \dots, X_n$ are any objects of $\cc$, and by the
uniqueness of the factorisation through a coend, we get that
\begin{equation*}
  \phi_{T'}=\phi_{T_L} \circ (\Gamma_l \otimes \id_{A^{\otimes (n-2)}}).
\end{equation*}
Therefore, since $\Gamma_l(\al \otimes \al)-\al \otimes \al \in \Negl_\cc(\un,A \otimes A)$, we get that
\begin{equation*}
  \tau_\cc(L';\alpha) = \phi_{T'}  \al^{\otimes n}
     = \phi_{T_L}  (\id_{A^{ \otimes 2}} \otimes \al^{\otimes (n-2)})\Gamma_l(\al \otimes \al)
     = \phi_{T_L} \al^{\otimes n}
     = \tau_\cc(L;\alpha).
\end{equation*}

To show Part (b), it suffices to remark that $\tau_\cc(\bigcirc^{\pm 1};\alpha) \in \Endo_\cc(\un)$ is invertible since
$\tau_\cc(\bigcirc^{\pm 1};\alpha)=\Theta_\pm \al\neq 0$ and $\Endo_\cc(\un)=\kk$ is a field.
\end{proof}
\begin{figure}[h]
   \begin{center}
       \psfrag{D}[B][B]{$T'_{(X_1, \dots,X_n)}$}
       \psfrag{Q}[B][B]{$T_{L(X_1, X_2 \otimes X_1, X_3, \dots,X_n)}$}
       \psfrag{X}[Bl][Bl]{$X_1$}
       \psfrag{i}[Bl][Bl]{$\scriptstyle X_1$}
       \psfrag{Y}[Bl][Bl]{$X_2$}
       \psfrag{Z}[Bl][Bl]{$X_3$}
       \psfrag{W}[Bl][Bl]{$X_n$}
       \psfrag{U}[Bl][Bl]{$\scriptstyle X_2 \otimes X_1$}
       \psfrag{I}[B][B]{$\id_{X_2 \otimes X_1}$}
       \psfrag{H}[B][B]{$i_{X_1}$}
       \psfrag{N}[B][B]{$\phi_{T_L}$}
       \psfrag{K}[B][B]{$i_{X_2}$}
       \psfrag{E}[B][B]{$i_{X_3}$}
       \psfrag{T}[B][B]{$i_{X_n}$}
       \psfrag{S}[B][B]{$i_{X_2 \otimes X_1}$}
       \psfrag{G}[B][B]{$\Gamma_l$}
       \psfrag{=}[Bl][Bl]{\textbf{=}}
       \psfrag{-}[B][B]{$\dots$}
       \psfrag{A}[Br][Br]{$\scriptstyle A$}
       \psfrag{B}{$\scriptstyle A$}
       \scalebox{.8}{\includegraphics{dem3man2.eps}}
   \end{center}
   \label{figdem3man2}
   \caption{}
\end{figure}

\subsection{Kirby elements via functors}
Let us see that Kirby elements can be ``pulled back'' via ribbon functors.

Let $\aaa$, $\bb$, and $\cc$ be ribbon \kt categories. Suppose that the functor \eqref{Ffunctor} of $\aaa$ admits a coend
$(A,i)$ and that the functor \eqref{Ffunctor} of $\bb$ admits a coend $(B,j)$. Let $\pi: \aaa \to \cc$ and $\iota:\bb \to
\cc$ be ribbon functors. Since $A$ and $B$ are categorical Hopf algebras and $\pi$ and $\iota$ are ribbon functors, the
objects $\pi(A)$ and $\iota(B)$ are Hopf algebras in $\cc$ (with structure maps induced by $\pi$ and $\iota$
respectively).

\begin{prop}\label{pullback}
Suppose that $\pi$ is surjective, $\iota$ is full and faithful, and that there exists a Hopf algebra morphism
$\varphi:\iota(B) \to \pi(A)$ such that $\pi(i_X)=\varphi \circ \iota(j_Y)$ for all objects $X \in \aaa$ and $Y \in \bb$
with $\pi(X)=\iota(Y)$. Let $\be \in  \Homo_\bb(\un,B)$ and $\al \in \Homo_\aaa(\un,A)$ such that $\pi(\al)=\varphi \circ
\iota(\be)$.
\begin{enumerate}
 \renewcommand{\labelenumi}{{\rm (\alph{enumi})}}
 \item If $\be \in \Ki(\bb)$, then  $\al \in
       \Ki(\aaa)$ and $\tau_\aaa(L;\al)=\tau_\bb(L;\be)$ for any framed link~$L$.
 \item If $\be \in \Ki(\bb)^\norm$, then $\al \in
       \Ki(\aaa)^\norm$ and $\tau_\aaa(M;\al)=\tau_\bb(M;\be)$ for any 3-manifold $M$.
\end{enumerate}
\end{prop}
\begin{proof}
Let us prove Part (a). Suppose that $\be \in \Ki(\bb)$. Since the structure maps of $\pi(A)$ and $\iota(B)$ are induced
by $\pi$ and $\iota$ from those of $A$ and $B$ respectively, and since $\varphi:\iota(B) \to \pi(A)$ is a Hopf algebra
morphism, we have:
\begin{align*}
\pi(S_A \al-\al) & =(S_{\pi(A)}- \id_{\pi(A)})\pi(\al)
   =(S_{\pi(A)}- \id_{\pi(A)})\varphi \iota(\be)\\
 & =\varphi (S_{\iota(B)}-\id_{\iota(B)})\iota(\be)
   =\varphi \iota(S_B\be-\be)
\end{align*}
and
\begin{align*}
\pi(\Gamma^A_l (\al \otimes \al) -\al \otimes \al)
  & =(\Gamma^{\pi(A)}_l-\id_{\pi(A)^{\otimes2}})( \pi(\al) \otimes \pi(\al))\\
  & =(\Gamma^{\pi(A)}_l-\id_{\pi(A)^{\otimes2}})\varphi \iota(\be \otimes \be)\\
  & =\varphi (\Gamma^{\iota(B)}_l-\id_{\iota(B)^{\otimes2}}) \iota(\be \otimes \be)\\
  & =\varphi \iota(\Gamma^B_l(\be \otimes \be)-\be \otimes \be).
\end{align*}
Now, since $S_B\be-\be$ and $\Gamma^B_l(\be \otimes \be)-\be \otimes \be$ are negligible in $\bb$, $\iota$ is full, and
$\tr^\bb=\tr^\cc\circ \iota$, we get that $\pi(S_A \al-\al)$ and $\pi(\Gamma^A_l (\al \otimes \al) -\al \otimes \al)$ are
negligible in $\cc$. Hence, since $\pi$ is surjective and $\tr^\cc\circ \pi=\tr^\aaa$, the morphisms $S_A \al-\al$ and
$\Gamma^A_l (\al \otimes \al) -\al \otimes \al$ are negligible in $\aaa$, that is, $\al \in \Ki(\aaa)$.

Let $L=L_1 \cup \cdots \cup L_n$ be a framed link in $S^3$. Let $T_L$ be a $n$-special tangle such that $L$ is isotopic
$T_L \circ (\cup_- \otimes \cdots \otimes \cup_-)$, where the $i$th cup (with clockwise orientation) corresponds to the
component $L_i$. Let $Y_1, \dots, Y_n$ be any objects of $\bb$. Since $\pi$ is surjective, there exists objets $X_1,
\dots, X_n$ of $\aaa$ such that $\pi(X_k)=\iota(Y_k)$. Recall that, by assumption, $\pi(i_{X_k})=\varphi
\iota(j_{Y_k})$. Since $\iota$ is full and the domain and codomain of the morphism
$\pi(\phi_{T_L}^{\aaa})\varphi^{\otimes n}$ of $\cc$ are $\iota(B^{\otimes n})$ and $\un=\iota(\un)$ respectively, there
exists a morphism $\xi:B^{\otimes n} \to \un$ in $\bb$ such that $\iota(\xi)=\pi(\phi_{T_L}^{\aaa})\varphi^{\otimes n}$.
Then
\begin{align*}
 \iota\bigl(\phi_{T_L}^{\bb} \circ ( j_{Y_1} \otimes \cdots \otimes j_{Y_n})\bigr)
 & = \iota(T_{L(Y_1, \dots , X_n)}^{\bb}) \\
 & = T_{L(\iota(Y_1), \dots , \iota(Y_n))}^{\cc} \\
 & = T_{L(\pi(X_1), \dots , \pi(X_n))}^{\cc} \\
 & = \pi(T_{L(X_1, \dots , X_n)}^\aaa ) \\
 & = \pi(\phi_{T_L}^{\aaa})(\pi(i_{X_1}) \otimes \cdots \otimes \pi(i_{X_n}))\\
 & = \pi(\phi_{T_L}^{\aaa})(\varphi \iota(j_{Y_1}) \otimes \cdots \otimes \varphi \iota(j_{Y_n}))\\
 & = \pi(\phi_{T_L}^{\aaa})\varphi^{\otimes n} \circ\iota(j_{Y_1} \otimes \cdots \otimes j_{Y_n})\\
 & = \iota\bigl(\xi \circ ( j_{Y_1} \otimes \cdots \otimes j_{Y_n})\bigr).
\end{align*}
Therefore, since $\iota$ is faithful, $\phi_{T_L}^{\bb} \circ ( j_{Y_1} \otimes \cdots \otimes j_{Y_n})=\xi \circ (
j_{Y_1} \otimes \cdots \otimes j_{Y_n})$ and so, by the uniqueness of the factorization through a coend, we get that
$\phi_{T_L}^{\bb}=\xi$, that is, $\iota(\phi_{T_L}^{\bb})=\pi(\phi_{T_L}^{\aaa})\varphi^{\otimes n}$. Hence, since the
maps $\Endo_\aaa(\un)=\kk \to \Endo_\cc(\un)=\kk$ and  $\Endo_\bb(\un)=\kk \to \Endo_\cc(\un)=\kk$ induced by $\pi$ and
$\iota$ respectively are the identity of $\kk$, we have that
\begin{align*}
\tau_\cc(L;\al)& = \pi(\tau_\cc(L;\al))   = \pi(\phi_{T_L}^{\aaa} \al^{\otimes n})\\
 & = \pi(\phi_{T_L}^{\aaa})\varphi^{\otimes n}\iota(\be^{\otimes n})\\
 & = \iota(\phi_{T_L}^{\bb} \be^{\otimes n}) = \iota(\tau_{\bb}(L;\be))=\tau_{\bb}(L;\be).
\end{align*}

Let us prove Part (b). Suppose that $\be \in \Ki(\bb)^\norm$. Let $Y$ be any objet of~$\bb$. Since $\pi$ is surjective,
there exists an object $X$ of $\aaa$ such that $\pi(X)=\iota(Y)$. Recall that, by assumption, $\pi(i_X)=\varphi \circ
\iota(j_Y)$. Since $\iota$ is full and the domain and codomain of the morphism $\pi(\Theta^{\aaa}_\pm )\varphi$ of $\cc$
are $\iota(B)$ and $\un=\iota(\un)$ respectively, there exists a morphism $\varsigma_\pm:B \to \un$ in $\bb$ such that
$\iota(\varsigma_\pm)=\pi(\Theta^{\aaa}_\pm )\varphi$. We have that
\begin{align*}
\iota(\varsigma_\pm\circ  \iota(j_Y)) & =\pi(\Theta^{\aaa}_\pm )\varphi \circ  \iota(j_Y)
  =\pi(\Theta^{\aaa}_\pm i_X)
=\pi(\ev^{\aaa}_X(\id_{X^*} \otimes \theta^{\aaa \pm 1}_X))\\
  & =\ev^{\cc}_{\pi(X)}(\id_{\pi(X)^*} \otimes \theta^{\cc \pm 1}_{\pi(X)})
  =\iota(\ev^{\bb}_Y)(\iota(\id_{Y^*}) \otimes \iota(\theta^{\bb}_Y)^{\pm 1})\\
  &=\iota(\ev^{\bb}_Y(\id_{Y^*} \otimes\theta^{\bb\pm 1}_Y)) =\iota(\Theta^{\bb}_\pm \circ j_Y).
\end{align*}
Therefore, since $\iota$ is faithful, $\Theta^{\bb}_\pm \circ j_Y=\varsigma_\pm\circ  \iota(j_Y)$ and so, by the
uniqueness of the factorization through a coend, we get that $\Theta^{\bb}_\pm=\varsigma_\pm$, that is,
$\iota(\Theta^{\bb}_\pm)=\pi(\Theta^{\aaa}_\pm )\varphi$. Then
\begin{equation*}
\Theta^{\cc}_\pm \al=\pi(\Theta^{\aaa}_\pm \al)=\pi(\Theta^{\aaa}_\pm) \varphi \iota(\be)=\iota(\Theta^{\bb}_\pm
\be)=\Theta^{\bb}_\pm \be.
\end{equation*}
Hence, since $\Theta^{\bb}_\pm \be \neq 0$ and by Part (a), we get that $\al \in \Ki(\aaa)^\norm$ and
$\tau_\aaa(M;\al)=\tau_\bb(M;\be)$ for any 3-manifold $M$.
\end{proof}

\section{The case of semisimple ribbon categories}\label{RTinv}

In this section, we focus on the case of semisimple categories $\cc$. We give necessary conditions for being in
$\Ki(\cc)$. In particular, we show that there exist (even in the non-modular case) elements of $\Ki(\cc)$ corresponding
to Reshetikhin-Turaev invariants computed with finitely semisimple ribbon subcategories of $\cc$. Moreover, we study
Kirby elements coming from semisimplification of ribbon categories.

\subsection{Semisimple categories}\label{sect-defsscat}
Recall that a category $\cc$ \emph{admits (finite) direct sums} if, for any finite set of objects $X_1, \dots,X_n$ of
$\cc$, there exists an object $X$ and morphisms $p_i:X \to X_i$ such that, for any object $Y$ and morphisms $f_i:Y \to
X_i$, there is a unique morphism $f:Y \to X$ with $p_i \circ f=f_i$ for all $i$. The object $X$ is then unique up to
isomorphism. We write $X=\oplus_i X_i$ and $f=\oplus_i f_i$.

A \kt category is \emph{abelian} if it admits (finite) direct sums, every morphism has a kernel and a cokernel, every
monomorphism is the kernel of its cokernel, every epimorphism is the cokernel of its kernel, and every morphism is
expressible as the composite of an epimorphism followed by a monomorphism. In particular, an abelian category admits a
null object (which is unique up to isomorphism). Note that a morphism of an abelian \kt category which is both a
monomorphism and an epimorphism is an isomorphism.

Let $\cc$ be an abelian \kt category. A non-null object $U$ of $\cc$ is \emph{simple} if every non-zero monomorphism $V
\to U$ is an isomorphism, and every non-zero epimorphism $U \to V$ is an isomorphism. Any non-zero morphism between
simple objects is an isomorphism. An object $U$ of $\cc$ is \emph{scalar} if $\Endo_\cc(V)=\kk$. Note that if $\kk$ is
algebraically closed, then every simple object is scalar. An object of $\cc$ is \emph{indecomposable} if it cannot be
written as a direct sum of two non-null objects. Note that every scalar or simple object is indecomposable.

By a \emph{semisimple \kt category}, we shall mean an abelian \kt category for which every object is a (finite) direct
sum of simple objects. By a \emph{finitely semisimple \kt category}, we shall mean a semisimple \kt category which has
finitely many isomorphism classes of simple objects. Note that in a semisimple \kt category, every scalar or
indecomposable object is simple.

In a semisimple ribbon \kt category $\cc$, any negligible morphism is null. Therefore, for every object $X,Y$ of $\cc$,
the pairing
\begin{equation}\label{pairss}
\Homo_\cc(X,Y) \otimes \Homo_\cc(Y,X) \to \kk \quad f \otimes g \mapsto \tr(gf)
\end{equation}
is non-degenerate. Note that this implies that the quantum dimension of a scalar object of $\cc$ is invertible.

\begin{lem}\label{lemidsum}
Let $\cc$ be a finitely semisimple ribbon \kt category whose simple objects are scalar. Let $\Lambda$ be a (finite) set
of representatives of isomorphism classes of simple objects of $\cc$. Fix an object $X$ of $\cc$. For any $\lambda \in
\Lambda$, set
\begin{equation*}
n_\lambda=\kdim\Homo_\cc(\lambda,X)=\kdim\Homo_\cc(X,\lambda)
\end{equation*}
and let $\{f_i^\lambda \, | \, 1 \leq i \leq n_\lambda\}$ be a basis of $\Homo_\cc(\lambda,X)$ and $\{g_i^\lambda\, | \,
1 \leq i \leq n_\lambda\}$ be a basis of $\Homo_\cc(X,\lambda)$ such that $ g_i^\lambda f_j^\lambda= \delta_{i,j} \,
\id_{\lambda}$ for all $1 \leq i,j \leq n_\lambda$ (such basis exist since the pairing \eqref{pairss} is non-degenerate).
Then
\begin{equation*}
\id_X=\sum_{\lambda \in \Lambda}\, \sum_{1 \leq i \leq n_\lambda} f_i^\lambda g_i^\lambda.
\end{equation*}
\end{lem}
\begin{proof}
Note that since $\cc$ is semisimple with scalar simple objects, the $\Homo_\cc$'s \kt spaces are finite-dimensional.
Since the category $\cc$ is semisimple, the composition induces a \kt linear isomorphism $\oplus_{\lambda \in \Lambda}
\Homo_\cc(X, \lambda) \otimes \Homo_\cc(\lambda,X) \to \Endo_\cc(X)$. Therefore, for all $\lambda \in \Lambda$ and $1\leq
i,j \leq n_\lambda$, there exist $a_{\lambda,i,j}\in\kk$ such that $\id_X=\sum_{\lambda \in \Lambda}\sum_{1 \leq i,j \leq
n_\lambda}a_{\lambda,i,j} \, f_i^\lambda g_j^\lambda$. Let $\lambda \in \Lambda$ and $1\leq i,j \leq n_\lambda$. Then
\begin{align*}
\delta_{i,j} \, \id_\lambda & = g_i^\lambda f_j^\lambda = g_i^\lambda\, \id_X \, f_j^\lambda
 = \sum_{\mu \in \Lambda}\sum_{1 \leq k,l \leq n_\lambda} \!\!\! a_{\mu,k,l} \, g_i^\lambda f_k^\mu g_l^\mu f_j^\lambda \\
 &= \sum_{\mu \in \Lambda}\sum_{1 \leq k,l \leq n_\lambda} \!\!\! a_{\mu,k,l} \, \delta_{\lambda,\mu} \,
    \delta_{i,k} \, \delta_{j,l} \, \id_\mu
 =  a_{\lambda,i,j} \, \id_\lambda,
\end{align*}
and so $a_{\lambda,i,j}=\delta_{i,j}$. Hence $\id_X=\sum_{\lambda \in \Lambda}\sum_{1 \leq i \leq n_\lambda} f_i^\lambda
g_i^\lambda$.
\end{proof}

\subsection{Kirby elements of finitely semisimple ribbon categories}\label{sect-Kiss}
Let $\bb$ be a finitely semisimple ribbon \kt category whose simple objects are scalar.

Note that under these assumptions, the \kt space $\Homo_\bb(X,Y)$ is finite-dimensional for any objects $X,Y$ of $\bb$.

Denote by $\Lambda$ a (finite) set of representatives of isomorphism classes of simple objects of $\bb$. We can suppose
that $\un \in \Lambda$. For any $\lambda \in \Lambda$, there exists a unique $\lambda^\vee \in \Lambda$ such that
$\lambda^*\simeq \lambda^\vee$.  The map $^\vee : \Lambda \to \Lambda$ is an involution and $\un^\vee=\un$. Recall that
$\qdim(\lambda) \neq 0$ for any $\lambda \in \Lambda$.

Set
\begin{equation}
   B= \bigoplus_{\lambda \in \Lambda} \lambda^* \otimes \lambda \in \bb.
\end{equation}
In particular, there exist morphisms $p_\lambda:B \to \lambda^* \otimes \lambda$ and $q_\lambda: \lambda^* \otimes
\lambda \to B$ such that $\id_B= \sum_{\lambda \in \Lambda} q_\lambda p_\lambda$ and $ p_\lambda q_\mu =
\delta_{\lambda,\mu}\, \id_{\lambda^* \otimes \lambda}$. Let $X$ be an object of $\bb$. Since $\bb$ is semisimple, we can
write $X=\oplus_{i \in I} \lambda_i$, where $I$ is a finite set and $\lambda_i \in \Lambda$. We set
\begin{equation}\label{coendinss}
  j_X = \sum_{i \in I} q_{\lambda_i} \circ (q_i^* \otimes p_i)
  : X^* \otimes X \to B,
\end{equation}
where $p_i : X \to \lambda_i$ and $q_i: \lambda_i \to X$ are morphisms in $\bb$ such that $\id_X= \sum_{i \in I} q_i
 p_i$ and $p_i q_j =\delta_{i,j} \, \id_{\lambda_i}$. Note that $j_X$ does not depend on the choice of such a family
of morphisms. Remark that $j_\lambda=q_\lambda$ for any $\lambda \in \Lambda$. It is not difficult to verify that $(B,
j)$ is a coend of the functor \eqref{Ffunctor} of $\bb$. By Section \ref{coendHA}, the object $B$ is a Hopf algebra in
$\bb$.

For any $\lambda \in \Lambda$, set
\begin{equation}
e_\lambda=j_\lambda \tcoev_\lambda : \un \to B \quad \text{and} \quad f_\lambda=\ev_\lambda p_\lambda : B \to \un.
\end{equation}
Note that $f_\lambda e_\mu=\delta_{\lambda,\mu} \qdim (\lambda)$ for any $\lambda,\mu \in \Lambda$.

\begin{lem}\label{lemavRT}
\begin{enumerate}
\renewcommand{\labelenumi}{{\rm (\alph{enumi})}}
\item $(e_\lambda)_{\lambda \in \Lambda}$ is a basis of the \kt space $\Homo_\bb(\un,B)$.
\item $(f_\lambda)_{\lambda \in \Lambda}$ is a basis of the \kt space $\Homo_\bb(B,\un)$.
\item For any $\lambda,\mu \in \Lambda$,
\begin{align*}
& S_B(e_\lambda)=e_{\lambda^\vee}, \qquad \eta_B=e_\un,  \qquad \varepsilon_B(e_\lambda)=\qdim(\lambda),\\
& (\id_B \otimes f_\mu)\Delta_B(e_\lambda)=\delta_{\lambda,\mu}e_\lambda,  \qquad m_B(e_\lambda \otimes e_\mu)=\sum_{\nu
\in \Lambda} N_{\lambda,\mu}^\nu e_\nu,
\end{align*}
where $ N_{\lambda,\mu}^\nu=\kdim \Homo_\bb(\lambda \otimes \mu,\nu)$.
\end{enumerate}
\end{lem}
\begin{proof}
Let us prove Part (a). For any $\lambda \in \Lambda$, the \kt space $\Homo_\bb(\un, \lambda^* \otimes \lambda)$ is
one-dimensional (since $\lambda$ is scalar) with basis $\tcoev_\lambda$. Therefore there exists $x_\lambda \in \kk$ such
that $p_\lambda f=x_\lambda \tcoev_\lambda$, and so $f=\id_B f=\sum_{\lambda \in \Lambda} q_\lambda p_\lambda
f=\sum_{\lambda \in \Lambda} x_\lambda e_\lambda$. Hence the family $(e_\lambda)_{\lambda \in \Lambda}$ generates
$\Homo_\bb(\un,B)$. To show that it is free, suppose that $\sum_{\lambda \in \Lambda} x_\lambda e_\lambda=0$. Then, for
any $\mu \in \Lambda$, $0=\sum_{\lambda \in \Lambda} x_\lambda f_\mu e_\lambda=\qdim(\mu) x_\mu$ and so $x_\mu=0$ since
$\qdim(\mu) \neq 0$.

Part (b) can be shown similarly. Let us prove Part (c). Let $\lambda,\mu \in \Lambda$. By definition,
$\eta_B=j_\un=j_\un\tcoev_\un=e_\un$ (since $\tcoev_\un=\id_\un$) and $\varepsilon_B(e_\lambda)=\varepsilon_B j_\lambda
\tcoev_\lambda=\ev_\lambda \tcoev_\lambda=\qdim(\lambda)$. The equalities $S_B(e_\lambda)=e_{\lambda^\vee}$ and $(\id_B
\otimes f_\mu)\Delta_B(e_\lambda)=\delta_{\lambda,\mu}e_\lambda$ are shown in Figure~\ref{demRT1fig} and \ref{demRT2fig}
respectively. Write $\lambda \otimes \mu=\oplus_{i \in I} \lambda_i$. In particular, there exits morphisms $p_i:\lambda
\otimes \mu \to \lambda_i$ and $q_i:\lambda_i \to \lambda \otimes \mu$ such that $\id_{\lambda\otimes \mu}=\sum_{i \in I}
q_ip_i$ and $p_i q_j=\delta_{i,j}\id_{\lambda_i}$. Recall that $j_{\lambda \otimes \mu}=\sum_{i \in I}
j_{\lambda_i}(q_i^* \otimes p_i)$. For any $\nu \in \Lambda$, since $\Homo_\bb(\lambda \otimes \mu,\nu)\cong \oplus_{i
\in I} \Homo_\bb(\lambda_i,\nu)$, we have that
\begin{equation*}
N_{\lambda,\mu}^\nu=\kdim \Homo_\bb(\lambda \otimes \mu,\nu)=\sum_{i \in I} \delta_{\lambda_i,\nu}.
\end{equation*}
Then the equality $m_B(e_\lambda \otimes e_\mu)=\sum_{\nu \in \Lambda} N_{\lambda,\mu}^\nu e_\nu$ is shown in
Figure~\ref{demRT3fig}.
\end{proof}
\begin{figure}[t]
        \subfigure[$S_B(e_\lambda)=e_{\lambda^\vee}$]{
                                       \label{demRT1fig}
                                       \psfrag{A}[Br][Br]{$B$}
                                       \psfrag{V}[Bl][Bl]{$\lambda^*$}
                                       \psfrag{U}[Bl][Bl]{$\lambda$}
                                       \psfrag{I}[B][B]{$\id_{\lambda^*}$}
                                       \psfrag{H}[B][B]{$S_B$}
                                       \psfrag{S}[B][B]{$e_\lambda$}
                                       \psfrag{X}[B][B]{$e_{\lambda^\vee}$}
                                       \psfrag{E}[B][B]{$j_{\lambda^*}$}
                                       \psfrag{=}[Bl][Bl]{\textbf{=}}
                                       \scalebox{.8}{\includegraphics{demlemRT1.eps}}}
        \subfigure[$(\id_B \otimes f_\mu)\Delta_B(e_\lambda)=\delta_{\lambda,\mu}e_\lambda$]{
                                       \label{demRT2fig}
                                       \psfrag{A}[Br][Br]{$B$}
                                       \psfrag{V}[Bl][Bl]{$\mu$}
                                       \psfrag{U}[Bl][Bl]{$\lambda$}
                                       \psfrag{I}[B][B]{$p_\mu$}
                                       \psfrag{d}[Bl][Bl]{$=\, \delta_{\lambda,\mu}$}
                                       \psfrag{P}[B][B]{$e_\lambda$}
                                       \psfrag{D}[B][B]{$\Delta_B$}
                                       \psfrag{N}[B][B]{$f_\mu$}
                                       \psfrag{E}[B][B]{$j_\lambda$}
                                       \psfrag{X}[B][B]{$q_\lambda$}
                                       \psfrag{=}[Bl][Bl]{\textbf{=}}
                                       \scalebox{.8}{\includegraphics{demlemRT2.eps}}}
        \caption{}
\end{figure}
\begin{figure}[t]
                                       \psfrag{A}[Br][Br]{$B$}
                                       \psfrag{D}[B][B]{$m_B$}
                                       \psfrag{F}[B][B]{$e_\lambda$}
                                       \psfrag{K}[B][B]{$e_\mu$}
                                       \psfrag{M}[B][B]{$e_\nu$}
                                       \psfrag{c}[B][B]{$j_{\lambda_i}$}
                                       \psfrag{g}[B][B]{$q_i^*$}
                                       \psfrag{p}[B][B]{$j_\nu$}
                                       \psfrag{r}[B][B]{$q_i$}
                                       \psfrag{o}[B][B]{$p_i$}
                                       \psfrag{E}[B][B]{$j_{\mu \otimes \lambda}$}
                                       \psfrag{I}[B][B]{$\id_{\mu \otimes \lambda}$}
                                       \psfrag{J}[B][B]{$\id_{\lambda \otimes \mu}$}
                                       \psfrag{H}[B][B]{$j_{\lambda \otimes \mu}$}
                                       \psfrag{V}[Bl][Bl]{$\mu$}
                                       \psfrag{w}[Bl][Bl]{$\nu$}
                                       \psfrag{X}[Bl][Bl]{$\lambda$}
                                       \psfrag{U}[Bl][Bl]{$\lambda_i$}
                                       \psfrag{T}[Bl][Bl]{$\mu \otimes \lambda$}
                                       \psfrag{S}[Bl][Bl]{$\lambda \otimes \mu$}
                                       \psfrag{L}[Bl][Bl]{$\ds=\, \sum_{i\in I}$}
                                       \psfrag{B}[Bl][Bl]{$\ds=\, \sum_{\nu \in \Lambda}\Bigl (
                                                                   \sum_{i\in I}\delta_{\lambda_i,\nu} \Bigr )$}
                                       \psfrag{y}[Bl][Bl]{$\ds=\, \sum_{\nu\in \Lambda} N_{\lambda,\mu}^\nu$}
                                       \psfrag{=}[Bl][Bl]{\textbf{=}}
                                       \scalebox{.8}{\includegraphics{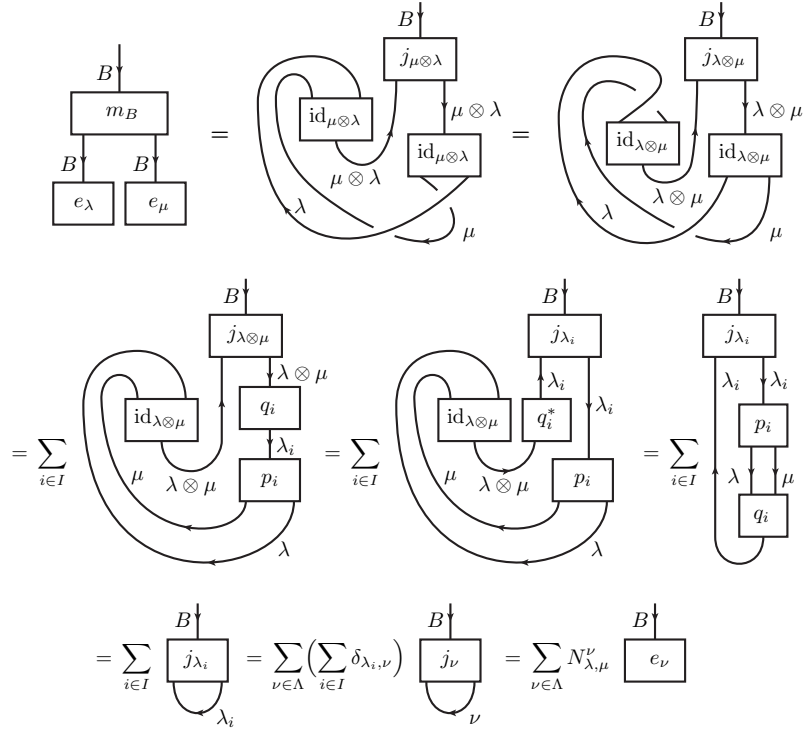}}
        \caption{$\Delta_B(e_\lambda \otimes e_{\mu})=\sum_{\nu\in \Lambda} N_{\lambda,\mu}^\nu e_\nu$}
        \label{demRT3fig}
\end{figure}

Since $\bb$ is semisimple and so its negligible morphisms are null, we have that:
\begin{equation*}
\Ki(\bb)=\{ \al \in \Homo_\bb(\un,A) \; | \; S_A\al=\al \, \text{ and } \, \Gamma_l(\al \otimes \al)=\al \otimes \al \}.
\end{equation*}

\begin{lem}\label{aleqal1dim}
Let $\al=\sum_{\lambda \in \Lambda} \al_\lambda e_\lambda \in \Homo_\bb(\un,B)$, where $\al_\lambda \in \kk$. Suppose
that $\al \in \Ki(\bb)$. Set $\Lambda_\al=\{ \lambda \in \Lambda \, | \, \al_\lambda \neq 0 \}$. Then $\al = \al_\un
\sum_{\lambda \in  \Lambda_\al} \qdim(\lambda) e_\lambda$. Moreover, the set $\Lambda_\al$ is such that
$\Lambda_\al^\vee=\Lambda_\al$ and $m_B(\al \otimes e_\lambda)=\qdim(\lambda) \, \al=m_B(e_\lambda \otimes \al)$ for all
$\lambda \in \Lambda_\al$.
\end{lem}
\begin{proof}
By Lemma~\ref{lemavRT}(c), since  $S_B(\al)=\al$, we have that $\al_\lambda=\al_{\lambda^\vee}$ for all $\lambda \in
\Lambda$ and so $\Lambda_\al^\vee=\Lambda_\al$. Let $\mu,\nu \in \Lambda$. Since $\Gamma_r(\al \otimes \al)=\al \otimes
\al$ and by using Lemma~\ref{lemavRT}(c), we have that
\begin{align*}
\qdim(\mu) \qdim(\nu)  \al_\nu \al_\mu & = (f_\nu \otimes f_\lambda)(\al \otimes \al)
= (f_\nu \otimes f_\mu)\Gamma_r(\al \otimes \al)\\
&=\sum_{\lambda,\omega \in \Lambda}\al_\lambda \al_\omega \,(f_\nu \otimes f_\mu)\Gamma_r(e_\lambda \otimes e_\omega) \\
&=\sum_{\lambda,\omega \in \Lambda}\al_\lambda \al_\omega \,
  f_\nu m_B(e_\lambda \otimes (\id_B \otimes f_\mu) \Delta_B(e_\omega)) \\
&=\sum_{\lambda \in \Lambda}\al_\lambda \al_\mu \, f_\nu m_B(e_\lambda \otimes e_\mu) \\
&=\sum_{\lambda,\omega \in \Lambda}\al_\lambda \al_\mu N_{\lambda,\mu}^\omega \, f_\nu e_\omega\\
&=\sum_{\lambda \in \Lambda}\al_\lambda \al_\mu N_{\lambda,\mu}^\nu \qdim(\nu),
\end{align*}
and so
\begin{equation}\label{equaquadrdeal}
\qdim(\mu) \al_\mu \al_\nu= \al_\mu\sum_{\lambda \in \Lambda} N_{\lambda,\mu}^\nu \al_\lambda.
\end{equation}
Note that $N_{\lambda,\mu}^\un=\kdim \Homo_\bb(\lambda \otimes \mu,\un)=\kdim \Homo_\bb(\lambda,
\mu^*)=\delta_{\lambda,\mu^\vee}$ for all $\lambda,\mu \in \Lambda$. Hence \eqref{equaquadrdeal} gives that $\qdim(\mu)
\al_\mu \al_\un=\al_\mu \al_{\mu^\vee}= \al_\mu^2$ and so $\al_\mu=\al_\un\qdim(\mu)$ whenever $\al_\mu \neq 0$, that is,
$\al = \al_\un \sum_{\lambda \in \Lambda_\al} \qdim(\lambda) e_\lambda$. Finally, for any $\mu \in \Lambda_\al$,
\begin{align*}
m_B(\al \otimes e_\mu)&=\sum_{\lambda \in \Lambda} \al_\lambda \, m_B(e_\lambda \otimes e_\mu) \\
&=\sum_{\lambda,\nu \in \Lambda}  N_{\lambda,\mu}^\nu \al_\lambda \, e_\nu\\
&=\sum_{\nu \in \Lambda} \qdim(\mu) \al_\nu e_\nu \quad \text{by \eqref{equaquadrdeal} since $\al_\mu \neq 0$}\\
&=\qdim(\mu) \al.
\end{align*}
Likewise, using the fact that $\Gamma_l(\al \otimes \al)=\al \otimes \al$, on gets that $m_B(e_\mu \otimes
\al)=\qdim(\mu) \al$.
\end{proof}

By Lemma~\ref{aleqal1dim}, determining $\Ki(\bb)$ resumes to find the subsets $E$ of $\Lambda$ for which $\sum_{\lambda
\in E} \qdim(\lambda) e_\lambda$ belongs to $\Ki(\bb)$. In the next theorem, we show that among these subsets, there are
those corresponding to monoidal subcategories of $\bb$.

\begin{thm}\label{RTsubmonoidal}
Let $\dd$ be a full ribbon and abelian subcategory of $\bb$. Let $\Lambda_\dd$ be a (finite) set of representatives of
isomorphism classes of simple objects of $\dd$. We can suppose that  $\Lambda_\dd \subset \Lambda$. Then $\sum_{\lambda
\in \Lambda_\dd} \qdim(\lambda) e_\lambda \in \Ki(\bb)$.
\end{thm}
\begin{remsbold}\label{rmsthsubmonoRT}
{\bf [a].} We do not know if every element of $\Ki(\bb)$ is of this form.\\[-.8em]

{\bf [b].} In Section~\ref{sectinvarRT}, we verify that $\sum_{\lambda \in \Lambda_\dd} \qdim(\lambda) e_\lambda$ leads
to
the Reshetikhin-Turaev invariant defined with $\dd$.\\[-.8em]

{\bf [c].} When $\bb$ is modular in the sense that the pairing $\omega_B:B \otimes B \to \un$ defined in \eqref{parBB} is
non-degenerate or, equivalently, that the $S$-matrix $[\tr(c_{\mu,\lambda}\circ c_{\lambda,\mu})]_{\lambda,\mu \in
\Lambda}$ is inversible, then $\sum_{\lambda \in \Lambda} \qdim(\lambda) e_\lambda$ is a two-sided integral of $B$ (see
\cite{Ker1}) and so belongs to $\Ki(\bb)$. Nevertheless, $\sum_{\lambda \in \Lambda} \qdim(\lambda) e_\lambda$ is not in
general a two-sided integral of $B$.
\end{remsbold}
\begin{proof}[Proof of Theorem~\ref{RTsubmonoidal}]
Firstly, since $\Lambda_\dd^\vee=\Lambda_\dd$ and by Lemma~\ref{lemavRT}(c),
\begin{equation*}
S_B(\al_\dd)=\sum_{\lambda \in \Lambda_\dd} \qdim(\lambda) S_B(e_\lambda)=\sum_{\lambda \in \Lambda_\dd}
\qdim(\lambda^\vee) e_{\lambda^\vee}=\al_\dd.
\end{equation*}
Secondly, to show that $\Gamma_r(\al_\dd \otimes \al_\dd)=\al_\dd \otimes \al_\dd$, it suffices to show that, for any
$\lambda \in \Lambda_\dd$,
\begin{equation}\label{ademRT}
\Gamma_r(\al_\dd \otimes e_\lambda)=\al_\dd \otimes e_\lambda.
\end{equation}
Fix $\lambda \in \Lambda_\dd$. Let $\mu,\nu \in \Lambda_\dd$ and set $n_{\mu,\nu}=\kdim \Homo_\bb(\nu \otimes
\lambda,\mu)$. Since $\dd$ is a full subcategory of $\bb$ and the pairing
\begin{equation*}
g \otimes f \in \Homo_\bb(\nu \otimes \lambda, \mu) \otimes \Homo_\bb(\mu, \nu \otimes \lambda) \mapsto \tr(gf) \in \kk
\end{equation*}
is non-degenerate, there exist basis $\{f_i^{\mu,\nu} \, | \, 1 \leq i \leq n_{\mu,\nu}\}$ of $\Homo_\bb(\mu,\nu \otimes
\lambda)$ and $\{g_i^{\mu,\nu} \, | \, 1 \leq i \leq n_{\mu,\nu}\}$ of $\Homo_\bb(\nu \otimes \lambda,\mu)$ such that
$g_j^{\mu,\nu} f_i^{\mu,\nu}= \delta_{i,j} \, \id_\mu$ for all $1 \leq i,j \leq n_{\mu,\nu}$. By Lemma~\ref{lemidsum}, we
have that
\begin{equation}\label{ideqdem1}
\sum_{\mu \in \Lambda} \sum_{1 \leq i \leq n_{\mu,\nu}} f_i^{\mu,\nu} g_i^{\mu,\nu} = \id_{\nu \otimes \lambda}.
\end{equation}
Let $\mu,\nu \in \Lambda_\dd$. For any $1 \leq i \leq n_{\mu,\nu}$, set
\begin{align*}
 & F_i^{\nu,\mu}=(g_i^{\mu,\nu} \otimes \id_\lambda)(\id_\nu \otimes \tcoev_\lambda) : \nu \to \mu \otimes \lambda^*,\\
 & G_i^{\nu,\mu}=(\id_\nu \otimes \tev_\lambda)(f_i^{\mu,\nu} \otimes \id_\lambda) : \mu \otimes \lambda^* \to \nu.
\end{align*}
Since $(\lambda^*, \tev_\lambda,\tcoev_\lambda)$ is a right duality for $\lambda$, we have that $\{F_i^{\nu,\mu}\, | \,
1\leq i \leq n_{\mu,\nu}\}$ is a basis for $\Homo_\bb(\nu,\mu\otimes \lambda^*)$ and $\{G_i^{\nu,\mu} \, | \, 1\leq i
\leq n_{\mu,\nu}\}$ is a basis for $\Homo_\bb(\mu\otimes \lambda^*,\nu)$. For any $1 \leq i,j \leq n_{\mu,\nu}$, since
$G_j^{\nu,\mu} F_i^{\nu,\mu}\in \Endo_\bb(\nu)$ and $\nu$ is scalar, we have that
\begin{align*}
\qdim(\nu) \, G_j^{\nu,\mu} F_i^{\nu,\mu} & = \tr(G_j^{\nu,\mu} F_i^{\nu,\mu}) \, \id_\nu
 = \tr( F_i^{\nu,\mu} G_j^{\nu,\mu}) \, \id_\nu\\
 & = \tr((g_i^{\mu,\nu} \otimes \id_\lambda)(\id_\nu \otimes \tcoev_\lambda)
                        (\id_\nu \otimes \tev_\lambda)(f_j^{\mu,\nu} \otimes \id_\lambda) ) \, \id_\nu \\
 & =  \tr(g_i^{\mu,\nu} f_j^{\mu,\nu}) \, \id_\nu
 = \tr(\delta_{i,j}\, \id_\mu) \, \id_\nu \\
 & =  \delta_{i,j}\, \qdim(\mu) \, \id_\nu.
\end{align*}
Therefore, by Lemma~\ref{lemidsum},
\begin{equation}\label{ideqdem2}
\sum_{\nu \in \Lambda} \sum_{1 \leq i \leq n_{\mu,\nu}} \qdim(\nu) \, F_i^{\nu,\mu} G_i^{\nu,\mu} = \qdim(\mu) \,
\id_{\mu \otimes \lambda^*}.
\end{equation}
Finally one gets \eqref{ademRT} by the equalities depicted in Figure~\ref{longuefigRT} which follow from the dinaturality
of $j$, the definition of $m_B$ and $\Delta_B$, and equalities \eqref{ideqdem1} and \eqref{ideqdem2}.
\end{proof}
\begin{figure}[t]
   \begin{center}
       \psfrag{A}[Br][Br]{$B$}
       \psfrag{V}[Bl][Bl]{$\nu$}
       \psfrag{X}[Bl][Bl]{$\lambda$}
       \psfrag{T}[Bl][Bl]{$\lambda\otimes \nu$}
       \psfrag{S}[Bl][Bl]{$\nu\otimes \lambda$}
       \psfrag{U}[Bl][Bl]{$\mu$}
       \psfrag{I}[B][B]{$\id_{\lambda \otimes \nu}$}
       \psfrag{J}[B][B]{$\id_{\nu \otimes \lambda}$}
       \psfrag{H}[B][B]{$S_B$}
       \psfrag{E}[B][B]{$j_{\lambda \otimes \nu}$}
       \psfrag{H}[B][B]{$j_{\nu \otimes \lambda}$}
       \psfrag{N}[B][B]{$j_\mu$}
       \psfrag{K}[B][B]{$j_\lambda$}
       \psfrag{B}[B][B]{$e_\lambda$}
       \psfrag{D}[B][B]{$\Gamma_r$}
       \psfrag{F}[B][B]{$\al_\dd$}
       \psfrag{f}[B][B]{$f_i^{\mu,\nu}$}
       \psfrag{g}[B][B]{$g_i^{\mu,\nu}$}
       \psfrag{O}[B][B]{$F_i^{\nu,\mu}$}
       \psfrag{P}[B][B]{$G_i^{\nu,\mu}$}
       \psfrag{x}[Bl][Bl]{$\ds=\sum_{\mu\in \Lambda_\dd} \qdim(\mu)$}
       \psfrag{u}[Bl][Bl]{$\ds\Gamma_r(\al_\dd \otimes e_\lambda)=\sum_{\nu\in \Lambda} \qdim(\nu)$}
       \psfrag{v}[Bl][Bl]{$\ds\!\!=\sum_{\nu\in \Lambda_\dd} \qdim(\nu)$}
       \psfrag{w}[Bl][Bl]{$\ds\;=\!\!\!\!\!\sum_{\scriptstyle \mu,\nu\in \Lambda_\dd \atop \scriptstyle 1 \leq i \leq
                                                                                  n_{\mu,\nu}}\!\!\!\!\!\qdim(\nu)$}
       \psfrag{c}[Bl][Bl]{$\ds=\sum_{\mu\in \Lambda_\dd}\sum_{\scriptstyle \nu\in \Lambda_\dd \atop \scriptstyle
                                                                   1 \leq i \leq n_{\mu,\nu}} \!\!\!\!\!\qdim(\nu)$}
       \psfrag{=}[Bl][Bl]{\textbf{=}}
       \scalebox{.8}{\includegraphics{demRT2.eps}}
   \end{center}
     \caption{$\Gamma_r(\al_\dd \otimes e_\lambda)=\al_\dd \otimes e_\lambda$}
     \label{longuefigRT}
\end{figure}

Recall that an object $X$ of $\bb$ is invertible if $X^* \otimes X$ is isomorphic to $\un$.
\begin{cor}
Suppose that every simple object of $\bb$ is invertible or, equivalently, that $(\Lambda,\otimes,\un)$ is a group. Then
$\Ki(\bb)$ is composed by the scalar multiples of $\al_G=\sum_{\lambda \in G} \qdim(\lambda) e_\lambda$, where $G$ is any
subgroup of $\Lambda$.
\end{cor}
\begin{proof}
Let $G$ be a subgroup of $\Lambda$. The full abelian subcategory $\dd$ of $\bb$ generated by $G$ is ribbon. Therefore, by
Theorem~\ref{RTsubmonoidal}, the element $\al_G=\sum_{\lambda \in G} \qdim(\lambda) e_\lambda$ (and so its scalar
multiples) belongs to $\Ki(\bb)$.

Conversely, let $\al=\sum_{\lambda \in \Lambda} \al_\lambda e_\lambda \in \Ki(\bb)$. If $\al=0$, then $\al=0 \cdot
\al_\Lambda$. Suppose that $\al \neq 0$. By Lemma~\ref{aleqal1dim},  $\al = \al_\un \sum_{\lambda \in \Lambda_\al}
\qdim(\lambda) e_\lambda$ where $\Lambda_\al=\{\lambda \in \Lambda \, | \, \al_\lambda \neq 0 \}$. Then $\un
\in\Lambda_\al$. Since every element of $\Lambda$ is invertible, we have that $N_{\lambda,\mu}^\nu=\delta_{\nu,\lambda
\otimes \mu}$ for any $\lambda,\mu,\nu \in \Lambda$. Let $\mu,\nu \in \Lambda_\al$. By using \eqref{equaquadrdeal}, $0
\neq \qdim(\mu) \al_\mu \al_\nu= \al_\mu \al_{\nu \otimes \mu^\vee}$ and so $\al_{\nu \otimes \mu^\vee}\neq 0$, that is,
$\nu \otimes \mu^\vee \in \Lambda_\al$. Hence $\Lambda_\al$ is a subgroup of $\Lambda$.
\end{proof}

\subsection{On the Reshetikhin-Turaev invariant}\label{sectinvarRT}
Let $\bb$ be a finitely semisimple ribbon \kt category whose simple objects are scalar. Let $\Lambda$ be a (finite) set
of representatives of isomorphism classes of simple objects containing $\un$.

Set $\Delta_{\pm}=\sum_{\lambda \in \Lambda} v_\lambda^{\pm 1} \qdim(\lambda)^2 \in \kk$, where $v_\lambda \in \kk$ is
the (invertible) scalar defined by $\theta_\lambda=v_\lambda \, \id_\lambda$. Recall (see \cite{Tur2,Brug1}) that the
Reshetikhin-Turaev invariant of 3-manifolds is well-defined when $\Delta_+ \neq 0 \neq \Delta_-$. Moreover, if $L$ is a
framed link in $S^3$, it is given by:
\begin{equation*}
 \RT_\bb(M_L)= \Delta_+^{b_-(L)-n_L}\,  \Delta_-^{-b_-(L)}
  \sum_{c \in \mathrm{Col}(L)} \Bigl ( \prod_{j=1}^n \qdim(c(L_j))
 \Bigr ) F(L,c) \in \kk.
\end{equation*}
Here $\mathrm{Col}(L)$ is the set of functions $c:\{L_1, \dots, L_n\} \to \Lambda$, where $L_1, \dots, L_n$ are the
components of $L$, and $F(L,c) \in \Endo_\bb(\un)=\kk$ is the morphism represented by a plane diagram of $L$ where the
component $L_j$ is colored with the object $c(L_j)$.

By Theorem~\ref{RTsubmonoidal}, $\al_\bb=\sum_{\lambda \in \Lambda} \qdim(\lambda) e_\lambda \in \Ki(\bb)$, where
$e_\lambda$ is defined as in Section~\ref{sect-Kiss}.

\begin{cor}\label{corRt}
Suppose that $\Delta_\pm \neq 0$. Then $\al_\bb \in \Ki(\bb)^\norm$ and $\tau_\bb(M;\al_\bb)=\RT_\bb(M)$ for any
3-manifold $M$.
\end{cor}

\begin{proof}
We have that
\begin{equation*}
  \Theta_{\pm} \al_\bb
   = \sum_{\lambda \in \Lambda} \qdim(\lambda) \, \ev_\lambda (\id_{\lambda^*} \otimes \theta_\lambda^{\pm 1}) \tcoev_\lambda
   =  \sum_{\lambda \in \Lambda} v_\lambda^{\pm 1} \qdim(\lambda)^{2}
   = \Delta_{\pm}  \neq 0.
\end{equation*}
Therefore, since $\al_\bb \in \Ki(\bb)$, one gets that $\al_\bb \in \Ki(\bb)^\norm$.

Let $L=L_1 \cup \cdots \cup L_n$ be a framed link in $S^3$. Let $T_L$ be a $n$-special tangle such that $L$ is isotopic
$T_L \circ (\cup_- \otimes \cdots \otimes \cup_-)$, where the $i$th cup (with clockwise orientation) corresponds to the
component $L_i$. Then, for any $c \in \mathrm{Col}(L)$,
\begin{equation*}
F(L,c)=\phi_{T_L} \circ  (j_{c(L_1)}\tcoev_{c(L_1)} \otimes \cdots \otimes j_{c(L_n)}\tcoev_{c(L_n)}) =\phi_{T_L} \circ
(e_{c(L_1)} \otimes \cdots \otimes e_{c(L_n)}).
\end{equation*}
Therefore
\begin{eqnarray*}
\lefteqn{\sum_{c \in \mathrm{Col}(L)} \Bigl ( \prod_{j=1}^n \qdim(c(L_j))
       \Bigr ) F(L,\lambda)}\\
 & = & \sum_{c \in \mathrm{Col}(L)}
        \phi_{T_L} \circ \Bigl (\qdim(c(L_1))\, e_{c(L_1)} \otimes \cdots \otimes
       \qdim(c(L_n))e_{c(L_n)}\Bigr)\\
 & = &  \phi_{T_L} \circ \Bigl (\sum_{\lambda_1 \in \Lambda}\qdim(\lambda_1)e_{\lambda_1} \otimes \cdots \otimes
       \sum_{\lambda_n \in \Lambda} \qdim(\lambda_n)e_{\lambda_n}\Bigr )\\
 & = & \phi_{T_L} \circ (\al_\bb \otimes \cdots \otimes \al_\bb ).
\end{eqnarray*}
Hence
\begin{eqnarray*}
 \RT_\bb(M_L)
 & = & \Delta_+^{b_-(L)-n_L}\,  \Delta_-^{-b_-(L)} \!\!\sum_{\lambda \in \mathrm{Col}(L)}
        \Bigl ( \prod_{j=1}^n \qdim(\lambda(L_j)) \Bigr ) F(L,\lambda),\\
 & = & (\Theta_+ \al_\bb)^{b_-(L)-n_L}\,  (\Theta_- \al_\bb)^{-b_-(L)} \;
        \phi_{T_L} \circ (\al_\bb \otimes \cdots \otimes \al_\bb )\\
 & = & \tau_\bb(M_L;\al_\bb).
\end{eqnarray*}
\end{proof}

\subsection{Semisimplification of ribbon categories}\label{sect-semisimplific}
Let $\cc$ be a ribbon \kt category. For any objects $X,Y \in \cc$, recall that $\Negl_\cc(X,Y)$ denotes the \kt space of
negligible morphisms of $\cc$ from $X$ to $Y$ (see Section~\ref{sect-negl}). Let $\cc^\sss$ be the category whose objects
are the same as in $\cc$, and whose morphisms are
\begin{equation}
\Homo_{\cc^\sss}(X,Y)=\Homo_{\cc}(X,Y)/\Negl_\cc(X,Y)
\end{equation}
for any objets $X,Y \in \cc^\sss$. The composition, monoidal structure, braiding, twist, and duality of $\cc^\sss$ are
induced by those of $\cc$.

When $\cc$ has finite-dimensional $\Homo$'s \kt spaces, the category $\cc^\sss$ is a semisimple ribbon \kt category,
called the \emph{semisimplification} of $\cc$, and the simple objects of $\cc^\sss$ are the indecomposable objects of
$\cc$ with non-zero quantum dimension, see~\cite{Brug2}.

Let $\pi: \cc \to \cc^\sss$ be the functor defined by $\pi(X)=X$ and $\pi(f)=f+\Negl_\cc(X,Y)$ for any object $X$ and any
morphism $f:X \to Y$ in $\cc$. This is a surjective ribbon functor. Note that $\pi$ is bijective on the objects.

\subsection{Kirby elements from semisimplification}\label{sect-lifromsscats}
Let $\cc$ be a ribbon \kt category which admits a coend $(A,i)$ for the functor \eqref{Ffunctor} and whose $\Homo$'s
spaces are finite-dimensional. Denote by $\cc^\sss$ the semisimplification of $\cc$ and let $\pi:\cc \to \cc^\sss$ be its
associated surjective ribbon functor (see Section~\ref{sect-semisimplific}).

Let $\bb$ be a full ribbon and abelian subcategory of $\cc^\sss$ which admits finitely many isomorphism classes of simple
objects, and whose simple objects are scalar. Let $\Lambda$ be a (finite) set of representatives of isomorphism classes
of simple objects of $\bb$ containing $\un$. For any object $X$ of $\cc^\sss$, we denote by $\pi^{-1}(X)$ the (unique)
object of $\cc$ such that $\pi(\pi^{-1}(X))=X$.

Let $B=\oplus_{\lambda \in \Lambda} \lambda^* \otimes \lambda$. In particular, there exist morphisms $p_\lambda:B \to
\lambda^* \otimes \lambda$ and $q_\lambda: \lambda^* \otimes \lambda \to B$ of $\bb$ such that $\id_B= \sum_{\lambda \in
\Lambda} q_\lambda p_\lambda$ and $p_\lambda q_\mu = \delta_{\lambda,\mu}\, \id_{\lambda^* \otimes \lambda}$. For any
object $X$ of $\bb$, we let $j_X:X^* \otimes X \to B$ as in \eqref{coendinss}. Recall that $(B,j)$ is a coend for the
functor $\eqref{Ffunctor}$ of $\bb$ and that $j_\lambda=q_\lambda$ for any $\lambda \in \Lambda$ (see
Section~\ref{sect-Kiss}).

Since $\bb$ is a full ribbon subcategory of $\cc^\sss$ and $\pi$ is a ribbon functor, the objects $B$ and $\pi(A)$ are
Hopf algebra in $\cc^\sss$. Set
\begin{equation}\label{phibb}
\varphi=\sum_{\lambda \in \Lambda} \pi(i_{\pi^{-1}(\lambda)})p_\lambda \in \Homo_{\cc^\sss}(B,\pi(A)).
\end{equation}

\begin{lem}\label{lemmseisipli}
$\varphi:B \to \pi(A)$ is a Hopf algebra morphism such that $\pi(i_X)=\varphi j_{\pi(X)}$ for all object $X$ of $\cc$
with $\pi(X) \in \bb$.
\end{lem}
\begin{proof}
Let $X$ be an object of $\cc$ such that $\pi(X) \in \bb$. Since $\bb$ is semisimple, we can write $\pi(X)=\oplus_{k \in
K} \lambda_k$, where $K$ is a finite set and $\lambda_k \in \Lambda$. Recall that $j_{\pi(X)} = \sum_{k \in K}
q_{\lambda_k} (q_k^* \otimes p_k)$, where $p_k : \pi(X) \to \lambda_k$ and $q_k: \lambda_k \to \pi(X)$ are morphisms in
$\bb$ such that $\id_{\pi(X)}= \sum_{k \in K} q_k p_k$ and $p_k q_l =\delta_{k,l} \, \id_{\lambda_k}$. For any $k \in K$,
since $\pi$ surjective, there exist morphisms $P_k : X \to \pi^{-1}(\lambda_k)$ and $Q_k: \pi^{-1}(\lambda_k) \to X$ in
$\cc$ such that $\pi(P_k)=p_k$ and $\pi(Q_k)=q_k$. Then, using the dinaturality of $i$ and since the functor $\pi$ is
ribbon,
\begin{align*}
\varphi j_{\pi(X)}
&=\sum_{\lambda \in \Lambda}\sum_{k \in K} \pi(i_{\pi^{-1}(\lambda)})p_\lambda q_{\lambda_k}  (q_k^*\otimes p_k)\\
&=\sum_{\lambda \in \Lambda}\sum_{\scriptstyle k \in K \atop \scriptstyle \lambda_k=\lambda}
\pi(i_{\pi^{-1}(\lambda)}) (q_k^*\otimes p_k)\\
&=\sum_{\lambda \in \Lambda}\sum_{\scriptstyle k \in K \atop \scriptstyle \lambda_k=\lambda}
\pi\bigl (i_{\pi^{-1}(\lambda)} (Q_k^*\otimes P_k) \bigr )\\
&=\sum_{\lambda \in \Lambda}\sum_{\scriptstyle k \in K \atop \scriptstyle \lambda_k=\lambda}
\pi\bigl (i_X (\id_{X^*}\otimes Q_kP_k) \bigr )\\
&=\pi(i_X) (\id_{\pi(X)^*}\otimes \sum_{k \in K }q_k p_k)\\
&=\pi(i_X) (\id_{\pi(X)^*}\otimes \id_{\pi(X)})=\pi(i_X).
\end{align*}

Let us verify that $\varphi$ is a Hopf algebra morphism. For any objects $X,Y$ of a ribbon category, set
\begin{equation*}
\be_X=(\tev_X \otimes\id_{X^{**}})(\id_X \otimes \coev_{X^*}):X \stackrel{\simeq}{\longrightarrow} X^{**}
\end{equation*}
and
\begin{align*}
\gamma_{X,Y}=&(\ev_{X^*} \otimes \id_{(Y \otimes X)^*})(\id_{X^*} \otimes \ev_{Y^*} \otimes \id_{(Y \otimes X)^*})\\
&(\id_{X^* \otimes Y^*} \otimes \coev_{Y \otimes X}):X^* \otimes Y^* \stackrel{\simeq}{\longrightarrow} (Y \otimes X)^*.
\end{align*}
Let $\lambda,\mu \in \Lambda$. Set $U=\pi^{-1}(\lambda)$ and $V=\pi^{-1}(\mu)$. Then
\begin{equation*}
\varepsilon_{\pi(A)}\varphi j_\lambda=\pi(\varepsilon_A i_U)=\pi(\ev_U) =\ev_\lambda=\varepsilon_B j_\lambda,
\end{equation*}
\begin{align*}
\Delta_{\pi(A)}\varphi j_\lambda &=\pi(\Delta_A i_U)\\
&=\pi ( (i_U \otimes i_U )
  (\id_{U^*} \otimes \coev_U \otimes \id_U) ) \\
&= (\varphi j_\lambda \otimes \varphi j_\lambda)
  (\id_{\lambda^*} \otimes \coev_\lambda \otimes \id_\lambda) = (\varphi \otimes \varphi) \Delta_B j_\lambda,
\end{align*}
\begin{align*}
m_{\pi(A)}(\varphi j_\lambda \otimes \varphi j_\mu) &=\pi(m_A (i_U \otimes i_V))\\
&=\pi(i_{V \otimes U} (\gamma_{U,V} \otimes \id_{V \otimes U})(\id_{U^*} \otimes c_{U, V^* \otimes V})) \\
&=\varphi j_{\mu \otimes \lambda} (\gamma_{\lambda,\mu} \otimes \id_{\mu \otimes \lambda})(\id_{\lambda^*} \otimes
c_{\lambda, \mu^* \otimes \mu})) \\
&=\varphi m_B(j_\lambda \otimes j_\mu),
\end{align*}
and
\begin{align*}
S_{\pi(A)}\varphi j_\lambda &=\pi(S_A i_U)=\pi ( i_{U^*} (\beta_U \otimes \theta_{U^*}^{-1})
c_{U^*,U})\\
&=\varphi j_{\lambda^*}(\beta_\lambda \otimes \theta_{\lambda^*}^{-1}) c_{\lambda^*,\lambda}=\varphi S_B j_\lambda.
\end{align*}
Therefore, since $\id_B= \sum_{\lambda \in \Lambda} j_\lambda p_\lambda$, we get that
$\varepsilon_{\pi(A)}\varphi=\varepsilon_B$, $\Delta_{\pi(A)}\varphi=(\varphi \otimes \varphi) \Delta_B$,
$m_{\pi(A)}(\varphi \otimes \varphi)=\varphi m_B$, and $S_{\pi(A)}\varphi=\varphi S_B$. Finally, we conclude by remarking
that $\varphi\eta_B=\varphi j_\un=\pi(i_\un)=\pi(\eta_A)=\eta_{\pi(A)}$.
\end{proof}
\begin{cor}\label{invfromsscat}
Let $\be \in  \Homo_\bb(\un,B)$ and $\al \in \Homo_\cc(\un,A)$ such that $\pi(\al)=\varphi \be$.
\begin{enumerate}
 \renewcommand{\labelenumi}{{\rm (\alph{enumi})}}
 \item If $\be \in \Ki(\bb)$, then  $\al \in
       \Ki(\cc)$ and $\tau_\cc(L;\al)=\tau_\bb(L;\be)$ for any framed link~$L$.
 \item If $\be \in \Ki(\bb)^\norm$, then $\al \in
       \Ki(\cc)^\norm$ and $\tau_\cc(M;\al)=\tau_\bb(M;\be)$ for any 3-manifold $M$.
\end{enumerate}
\end{cor}
\begin{proof}
This is an immediate consequence of Lemma~\ref{lemmseisipli} and Proposition~\ref{pullback}.
\end{proof}

\begin{rembold}
By Corollary~\ref{invfromsscat}, we have that
\begin{equation*}
\bigcup_\bb\pi^{-1}\bigl(\varphi_\bb(\Ki(\bb))\bigr) \subset \Ki(\cc) \quad \text{and} \quad \bigcup_\bb
\pi^{-1}\bigl(\varphi_\bb(\Ki(\bb)^\norm)\bigr) \subset \Ki(\cc)^\norm,
\end{equation*}
where $\bb$ runs over (equivalence classes of) finitely semisimple full ribbon and abelian subcategories of $\cc^\sss$
whose simple objects are scalar, and $\varphi_\bb$ is the morphism \eqref{phibb} corresponding to $\bb$. We will see in
Section~\ref{HKRsect} that these inclusions may be strict (see Remark~\ref{inclustrict}). This means that the
semisimplification process ``lacks'' some invariants.
\end{rembold}

\section{The case of categories of representations}\label{HKRsect}
In this section, we focus on the case of  the category $\reph$ of representations of a finite-dimensional ribbon Hopf
algebra $H$. In particular, we describe $\Ki(\reph)$ in purely algebraic terms. One of the interest of such a description
is to avoid the representation theory of $H$, which may be of wild type. Moreover, we show that the 3-manifolds
invariants obtained with these Kirby elements can be computed by using the Kauffman-Radford algorithm (even in the
non-unimodular case).

\subsection{Finite-dimensional Hopf algebras}\label{subsecfindimHA}
All considered algebras are supposed to be over the field $\kk$.

Let $H$ be a finite-dimensional Hopf algebra. Recall that a \emph{left} (resp.\@ \emph{right}) \emph{integral} for $H$ is
an element $\Lambda \in H$ such that $x\Lambda=\varepsilon(x) \Lambda$ (resp.\@ $\Lambda x =\varepsilon(x) \Lambda$) for
all $x \in H$. A left (resp.\@ right) integral for $H^*$ is then an element $\lambda \in H^*$ such that
$x_{(1)}\lambda(x_{(2)})=\lambda(x) \,1$ (resp.\@ $\lambda(x_{(1)})x_{(2)}=\lambda(x) \,1$) for all $x \in H$. Since $H$
is finite-dimensional, the space of left (resp.\@ right) integrals for $H$ is one-dimensional, and there always exist
non-zero right integral $\lambda$ for $H^*$ and a non-zero left integral $\Lambda$ for $H$ such that
$\lambda(\Lambda)=\lambda(S(\Lambda))=1$, see \cite[Proposition 3]{Rad5}.

By \cite[Corollary 2]{Rad1}, the space $H^*$ endowed with the right $H$-action defined, for any $f \in H^*$ and $h,x \in
H$, by
\begin{equation}\label{actionpoint}
\langle f \cdot h, x \rangle = \langle f,h x \rangle,
\end{equation}
is a free $H$-module of rank 1 with basis every non-zero right integral $\lambda$ for $H^*$. Likewise, $H$ endowed with
the right $H^*$-action $\leftharpoonup$ defined, for any $f \in H^*$ and $x \in H$, by
\begin{equation}\label{actharpoon}
 x \leftharpoonup f= f(x_{(1)}) x_{(2)},
\end{equation}
is a free $H^*$-module of rank 1 with basis $S(\Lambda)$, where $\Lambda$ is a non-zero left integral for $H$.

\begin{lem}\label{lemdufa}
Let $\lambda$ be a right integral for $H^*$ and $\Lambda$ be a left integral for $H$ such that
$\lambda(\Lambda)=\lambda(S(\Lambda))=1$. Let $a \in H$ and $f \in H^*$. Then
\begin{enumerate}
\renewcommand{\labelenumi}{{\rm (\alph{enumi})}}
\item $a=\lambda(a \Lambda_{(1)}) \, S(\Lambda_{(2)}) = \lambda(S(\Lambda_{(2)}) a) \, \Lambda_{(1)}$.
\item $f=\lambda \cdot a\,$ if and only if $\,a=(f \otimes S)\Delta (\Lambda)$.
\end{enumerate}
\end{lem}
\begin{proof}
Let us prove Part (a). Since $\lambda$ is a right integral for $H^*$ and $\lambda$ is a left integral for $H$ such that
$\lambda(\Lambda)=1$, we have that
\begin{eqnarray*}
\lambda(a \Lambda_{(1)}) \, S(\Lambda_{(2)})
 & = & \lambda(a_{(1)} \Lambda_{(1)}) \, a_{(2)} \Lambda_{(2)} \, S(\Lambda_{(3)}) \\
 & = & \lambda(a_{(1)} \Lambda_{(1)}) \, a_{(2)} \varepsilon(\Lambda_{(2)})\\
 & = & \lambda(a_{(1)} \Lambda) \, a_{(2)}\\
 & = & \lambda(\Lambda) \, \varepsilon(a_{(1)})a_{(2)} \\
 & = & a.
\end{eqnarray*}
The second equality of Part (a) can be proved similarly (by using the facts that $\lambda$ is a right integral for $H^*$
and $S(\Lambda)$ is a right integral for $H$ such that $\lambda(S(\Lambda))=1$).

Let us prove Part (b). If $f=\lambda \cdot a$ then, by Part (a),
\begin{equation*}
a=\lambda(a \Lambda_{(1)}) \, S(\Lambda_{(2)}) = f(\Lambda_{(1)}) \, S(\Lambda_{(2)}) = (f \otimes S)\Delta (\Lambda).
\end{equation*}
Conversely, if $a=(f \otimes S)\Delta (\Lambda)$ then, by Part (a),
\begin{equation*}
\lambda(a x)=\lambda(f(\Lambda_{(1)}) \, S(\Lambda_{(2)}) x) = f(\lambda(S(\Lambda_{(2)}) x)\Lambda_{(1)})=f(x)
\end{equation*}
for all $x \in H$, and so $f=\lambda \cdot a$.
\end{proof}
Recall that an element $h \in H$ is \emph{grouplike} if $\Delta(g)=g \otimes g$ and $\varepsilon(g)=1$. We denote by
$G(H)$ the space of grouplike elements of $H$. Recall (see \cite{Abe}) that there exist a unique grouplike element $g$ of
$H$ such that $x_{(1)}\lambda(x_{(2)})=\lambda(x) \,g$ for any $x \in H$ and any right integral $\lambda \in H^*$, and a
unique grouplike element $\nu \in G(H^*)=\mathrm{Alg}(H,\kk)$ such that $\Lambda x= \nu(x) \Lambda$ for any $x \in H$ and
any left integral $\Lambda \in H$. The element $g \in H$ (resp.\@ $\nu \in H^*$) is called the \emph{distinguished
grouplike element} of $H$ (resp.\@ of~$H^*$). The Hopf algebra $H$ is said to be \emph{unimodular} if its integrals are
to sided, that is, if $\nu=\varepsilon$.

By \cite[Theorem 3]{Rad1} and \cite[Proposition 3]{Rad1}, right integrals for $H^*$ and distinguished  grouplike elements
of $H$ and $H^*$ are related by:
\begin{align}
  & \lambda(xy)=\lambda(S^2(y\leftharpoonup \nu)x)=\lambda((S^2(y)\leftharpoonup \nu)x)
    \quad \text{for all $x,y \in H$}, \label{lambxy}\\
  & \lambda(S(x))=\lambda(gx) \quad \text{for all $x \in H$.}\label{Slambax}
\end{align}

\subsection{Quasitriangular Hopf algebras}  Following \cite{Drin}, a Hopf algebra $H$ is \emph{quasitriangular} if it
is endowed with an invertible element $R \in H \otimes H$ (the \emph{\R matrix}) such that, for any $x \in H$,
\begin{align}
  & R \Delta(x)= \sigma \Delta(x) R,\label{Rmat1}\\
  & (\id_H \otimes \Delta)(R)=R_{13} R_{12},\label{Rmat2}\\
  & (\Delta \otimes \id_H)(R)=R_{13} R_{23}. \label{Rmat3}
\end{align}
where $\sigma: H \otimes H \to H \otimes H$ denotes the usual flip map. The \R matrix verifies:
\begin{align}
 & (\varepsilon \otimes \id_H )(R)=(\id_H \otimes  \varepsilon) (R) = 1 , \label{pptRmat1}\\
 & (S \otimes \id_H)(R)= R^{-1}=(\id_H \otimes S^{-1})(R),\label{pptRmat2}\\
 & R_{23} R_{13} R_{12} = R_{12} R_{13} R_{23}.\label{pptRmat3}
\end{align}

The \emph{Drinfeld element} $u$ associated to $R$ is $u=m(S \otimes \id_H) \sigma(R) \in H$. It is invertible and
verifies:
\begin{align}
 & u^{-1}=m (\id_H \otimes S^2)(R_{21}), \label{pptu1}\\
 & S^2(x)=uxu^{-1} \quad \text{for all $x \in H$},\\
 & \Delta (u)=(R_{21}R)^{-1}(u \otimes u),\\
 & \varepsilon (u)=1.
\end{align}

Let $H$ be a finite-dimensional quasitriangular Hopf algebra and let $\nu \in H^*$ be the distinguished grouplike element
of $H^*$. Set $h_\nu=(\id_H \otimes \nu)(R) \in H$. By \eqref{Rmat3} and \eqref{pptRmat1}, $h_\nu$ is grouplike.

\subsection{Ribbon Hopf algebras}
Following \cite{RT2}, a \emph{ribbon} Hopf algebra is a quasitriangular Hopf algebra $H$ endowed with an invertible
element $\theta\in H$ (the \emph{twist}) such that:
\begin{align}
 & \text{$\theta$ is central,} \\
 & S(\theta)=\theta,\\
 & \Delta (\theta)=(\theta \otimes \theta) R_{21}R.
\end{align}

The twist verifies:
\begin{align}
 & \varepsilon (\theta_1) = 1, \\
 & \theta^{-2}=uS(u)=S(u)u.
\end{align}

Set $G=\theta u \in H$. Then $G$ is grouplike and verifies:
\begin{align}
 & S(u)=G^{-1} u G^{-1}, \label{SuGuG}\\
 & S^2(x)=GxG^{-1}. \label{ScarreG}
\end{align}
The element $G$ is called the \emph{special grouplike element} of $H$.

By \cite[Theorem 2]{Rad2} and \eqref{SuGuG}, the special grouplike element $G$ of a finite-dimensional ribbon Hopf
algebra $H$ is related to the distinguished grouplike element $g$ of $H$ and to $h_\nu=(\id_H \otimes \nu)(R)\in G(H)$ by
\begin{equation}\label{gG2h}
g=G^2 h_\nu.
\end{equation}
Relations \eqref{Slambax} and \eqref{gG2h} imply that
\begin{align}\label{lambdaSribb}
  \lambda(S(x))=\lambda(G^2 h_\nu x) \quad \text{for all $x \in H$.}
\end{align}

\subsection{Category of representations of a ribbon Hopf algebra}\label{sect-hmodu}
Let $H$ be a ribbon Hopf algebra with \R matrix $R$ and twist $\theta$. Denote by $\reph$ the \kt category of
finite-dimensional left $H$-modules and $H$-linear homomorphisms. The category $\reph$ is a monoidal \kt category with
tensor product and unit object defined in the usual way using the comultiplication and counit of $H$. The category
$\reph$ posses a left duality: for any module $M\in \reph$, set $M^*=\Homo_\kk(M,\kk)$, where $h\in H$ acts as the
transpose of $x\in M \mapsto S(h) x \in M$. The duality morphism $\ev_M: M^*\otimes M \to \un=\kk$ is the evaluation
pairing and, if $(e_k)_k$ is a basis of $M$ with dual basis $(e_k^*)_k$, then $\coev(1_\kk)=\sum_k e_k \otimes e_k^*$.
The category $\reph$ is braided: for modules $M, N \in \reph$, the braiding $c_{M,N}: M \otimes N \to N \otimes M$ is the
composition of multiplication by $R$ and the flip map $\sigma_{M,N}: M \otimes N \to N \otimes M$. The category $\reph$
is ribbon: for any module $M\in \reph$, the twist $\theta_M: M \to M$ is the multiplication by $\theta$. Recall, see
\eqref{catrightdual1} and \eqref{catrightdual2}, that $\reph$ posses a right duality $M \in \reph \mapsto (M^*, \tev_M,
\tcoev_M)$. Finally, the \kt category $\reph$ is pure, that is, $\Endo_{\reph}(\kk)=\kk$. Hence $\reph$ is a ribbon \kt
category in the sense of Section~\ref{ribkkcats}.

\begin{lem}\label{tevinreplem}
Let $G$ be the special grouplike element of $H$ and $M$ be a finite-dimensional left $H$-module. Then,
\begin{enumerate}
\renewcommand{\labelenumi}{{\rm (\alph{enumi})}}
\item $\tev_M(m \otimes f)=f(G m)$ for any $f \in M^*$ and $m \in M$.
\item $\tcoev_M=(\id_{M^*} \otimes G^{-1}\id_M) \sigma_{M,M^*} \coev_M$.
\end{enumerate}
\end{lem}
\begin{proof}
Write $R=\sum_i a_i \otimes b_i$. Recall that $u= \sum_i S(b_i) a_i$. Then
\begin{align*}
  \tev_M(m \otimes f)
     & = \ev_{M}\, c_{M,M^*} (\theta_M \otimes \id_{M^*})(m \otimes f)\\
     & = \ptsum_i \ev_{M} (b_i f \otimes a_i \theta m)
     = \ptsum_i f(S(b_i) a_i \theta m)\\
     &= f(u \theta  m)
     = f(G  m).
\end{align*}
Let $(e_k)_k$ be a basis of $M$ and $(e_k^*)_k$ be its dual basis. Note that if $g$ is any \kt linear endomorphism of
$M$, then $\sum_k g^*(e_k^*) \otimes e_k=\sum_k e_k^* \otimes g(e_k)$. For any $h \in H$, denote by $\rho(h)$ the \kt
linear endomorphism of $M$ defined by $m \in M \mapsto h m \in M$. By \eqref{pptRmat2} and \eqref{pptu1}, we have that
$R^{-1}=\sum_i S(a_i) \otimes b_i$ and $u^{-1}=\sum_i b_i S^2(a_i)$. Then
\begin{align*}
  \tcoev_M(1_\kk)
     & = (\id_{M^*} \otimes \theta_M^{-1}) (c_{M^*,M})^{-1} \coev_M(1_\kk)\\
     & = \ptsum_{k,i} S(a_i) e_k^* \otimes \theta^{-1}b_i \cdot e_k\\
     & = \ptsum_i \bigl (\id_{M^*} \otimes \rho(\theta^{-1}b_i) \bigr )
         \bigl ( \ptsum_k \rho(S^2(a_i))^*(e_k^*) \otimes e_k \bigr )\\
     & = \ptsum_i \bigl (\id_{M^*} \otimes \rho(\theta^{-1}b_i) \bigr )
         \bigl ( \ptsum_k e_k^* \otimes \rho(S^2(a_i))(e_k) \bigr ) \\
     & = \bigl(\id_{M^*} \otimes \rho (\theta^{-1}\ptsum_i b_iS^2(a_i) )\bigr ) \bigl
         ( \ptsum_k e_k^* \otimes e_k \bigr ) \\
     &=  \bigl(\id_{M^*} \otimes \rho(\theta^{-1}u^{-1})\bigr ) \sigma_{M,M^*} \coev_M(1_k) \\
     &=  \bigl(\id_{M^*} \otimes \rho(G^{-1})\bigr) \sigma_{M,M^*} \coev_M(1_k) \\
     &=  (\id_{M^*} \otimes G^{-1} \id_M) \sigma_{M,M^*} \coev_M(1_k).
\end{align*}
This completes the proof of the lemma.
\end{proof}
We immediately deduce from Lemma~\ref{tevinreplem} that, in the category $\reph$, we have
\begin{align*}
\tr(f)&=\Tr(G f)=\Tr(G^{-1} f),\\
\qdim(M)&=\Tr(G \,\id_{M})=\Tr(G^{-1} \id_{M}),
\end{align*}
for all module $M \in \reph$ and all $H$-linear endomorphism $f$ of $M$, where $\Tr$ denotes the usual trace of \kt
linear endomorphisms.

\subsection{Braided Hopf algebra associated to ribbon Hopf algebras}\label{sectbraidHA}
Let $H$ be a finite-dimensional ribbon Hopf algebra. The ribbon \kt category $\reph$ of finite-dimensional left
$H$-modules posses a coend $(A,i)$ for the functor \eqref{Ffunctor}. More precisely, $A=H^*=\Homo_\kk(H,\kk)$ as a \kt
space and is endowed with the coadjoint left $H$-action $\triangleright$ given, for any $f \in A$ and $h,x \in H$, by
\begin{equation}
\langle h \triangleright f, x \rangle = \langle f,S(h_{(1)}) x h_{(2)} \rangle,
\end{equation}
where $\langle, \rangle$ denotes the usual pairing between a \kt space and its dual. Given a module $M \in \reph$, the
map $i_M: M^* \otimes M \to A$ is given  by
\begin{equation}\label{imdefit}
  \langle i_M(l \otimes m),x \rangle=\langle l,x  m \rangle,
\end{equation}
for all $l \in M^*$, $m \in M$, and $x \in H$.

\begin{lem}\label{repcoendlemma}
If $\xi$ is a dinatural transformation from the functor \eqref{Ffunctor} of $\reph$ to a module $Z \in \reph$, then the
(unique) morphism $r: A \to Z$ such that $\xi_M=r \circ i_M$ for all $M \in \reph$ is given by $f \in A=H^* \mapsto
r(f)=\xi_{H} (f \otimes 1_H)$.
\end{lem}

Recall (see Section~\ref{coendHA}) that $A$ is a Hopf algebra in $\reph$. Using Lemma~\ref{repcoendlemma}, the structural
morphisms of $A$ can be explicitly described in terms of the structure maps of the Hopf algebra $H$. Nevertheless, it is
more convenient to write down its pre-dual structural morphisms: for example, since $A=H^*$ as a \kt space and $H$ is
finite-dimensional, the pre-dual of the multiplication $m_A:A \otimes A \to A$ of $A$ is a morphism $\Delta^\bd:H \to H
\otimes H$ such that $(\Delta^\bd)^*=m_A$. That yields a \kt space $H^\bd=H$ endowed with a comultiplication
$\Delta^\bd:H^\bd \to H^\bd \otimes H^\bd$, a counit $\varepsilon^\bd:H^\bd_1 \to \kk$, a unit $\eta^\bd: \kk \to H^\bd$,
a multiplication $m^\bd:H^\bd \otimes H^\bd \to H^\bd$, and an antipode $S^\bd: H^\bd \to H^\bd$. These structure maps,
described in the following lemma, verify the same axioms as those of a Hopf algebra except that the usual flip maps are
replaced by the braiding of $\reph$. The space $H^\bd$ is called the \emph{braided Hopf algebra associated to $H$}, see
\cite{Maj1,Lyu1}.
\begin{lem}[cf. \cite{Lyu2}] \label{braidedHopfalg}
The braided Hopf algebra $H^\bd$ associated to $H$ can be described as follows :
\begin{itemize}
  \item $H^\bd=H$ as an algebra;
  \item $\Delta^\bd(x)
             = \sum_i x_{(2)} a_i \otimes S(b_{i(1)}) x_{(1)} b_{i(2)}
             = \sum_i S(a_{i(1)}) x_{(1)} a_{i(2)} \otimes S(b_i) x_{(2)}$,
  \item $\varepsilon^\bd=\varepsilon$;
  \item $S^\bd_\al(x) = \sum_i S(a_i) \theta^2 S(x) u b_i
                        = \sum_i S(a_i) S(x) S(u)^{-1} b_i$,
\end{itemize}
for any $x \in H$, where $R=\sum_i a_i \otimes b_i$ is the \R matrix, $u$ the Drinfeld element, and $\theta$ the twist of
$H$.
\end{lem}

\subsection{Kirby elements of ribbon Hopf algebras}\label{subsect-kirbelemnH}
Throughout this section, $H$ will denote a finite-dimensional ribbon Hopf algebra with \R matrix $R \in H \otimes H$ and
twist $\theta \in H$. Let $u\in H$ be the Drinfeld element of $H$, $G=\theta u$ be the special grouplike element of $H$,
$\lambda \in H^*$ be a non-zero right integral for $H^*$, $\Lambda \in H$ be a non-zero left integral for $H$ such that
$\lambda(\Lambda)=\lambda(S(\Lambda))=1$, $g \in G(H)$ be the distinguished grouplike element of $H$, $\nu \in
G(H^*)=\mathrm{Alg}(H,\kk)$ be the distinguished grouplike element of $H^*$, and $h_\nu=(\id_H \otimes \nu)(R) \in G(H)$.
Let $(A,i)$ be the coend for the functor \eqref{Ffunctor} of $\reph$, and let $H^\bd$ be the braided Hopf algebra
associated to $H$. Recall that $A=H^*$ endowed with the coadjoint left $H$-action $\triangleright$ and that $A$ is a Hopf
algebra in the category $\reph$ whose structure maps are dual to those of~$H^\bd$. Recall that $\cdot$ denotes the right
action of $H$ on $H^*$ defined in \eqref{actionpoint} and $\leftharpoonup$ denotes the right
$H^*$-action on $H$ defined in \eqref{actharpoon}.\\

Set:
\begin{equation*}
\phi : \begin{cases} H & \to \Homo_\kk(\kk, H^*) \\ z & \mapsto \phi_z =(\al \mapsto \al \lambda \cdot z)
\end{cases},
\qquad T : \begin{cases} H & \to H \\ z & \mapsto T(z)=(S(z)\leftharpoonup \nu)h_\nu
\end{cases},
\end{equation*}
and
\begin{equation*}
\psi : \begin{cases} H^{\otimes n} & \to \Homo_\kk(H^{* \otimes n},\kk) \\ X & \mapsto \psi_X = \bigl (f_1\otimes \cdots
\otimes f_n \mapsto \langle f_1 \otimes \cdots \otimes f_n, X \rangle \bigr )
\end{cases}.
\end{equation*}
Denote by $\triangleleft$ be the right action of $H$ on $H^{\otimes n}$ given by
\begin{equation*}
(x_1\otimes \cdots \otimes x_n) \triangleleft h= S(h_{(1)}) x_1 h_{(2)}\otimes \cdots \otimes S(h_{(2n-1)}) x_n h_{(2n)}
\end{equation*}
for any $h \in H$ and $x_1, \dots, x_n \in H$. Set
\begin{align*}
 & L(H)=\{ z \in H \, | \, (x \leftharpoonup \nu) z = z x \text{ for any } x \in H \},\\
 & N(H)=\{ z \in L(H) \, | \, \lambda(za) = 0 \text{ for any } a \in Z(H) \},\\
 & V_n(H)=\{ X \in H^{\otimes n} \, | \, X \triangleleft h = \varepsilon(h) X \text{ for any } h \in H \}.
\end{align*}

Note that if $H$ is unimodular, then $L(H)=Z(H)$ and $T(z)=S(z)$ for all $z \in H$.

By the definition of $L(H)$ and by \eqref{lambxy}, we have that
\begin{equation}\label{eqphihlin}
\lambda(zxy)=\lambda(zS^2(y)x)
\end{equation}
for all $z \in L(H)$ and $x,y \in H$.

\begin{lem}\label{lemdespiptspaces}
\begin{enumerate}
\renewcommand{\labelenumi}{{\rm (\alph{enumi})}}
\item The map $\phi$ is an \kt isomorphism with inverse given by:
$$\al \in \Homo_\kk(\kk, H^*) \mapsto \phi^{-1}(\al)=(\al(1_\kk) \otimes S) \Delta (\Lambda)\in H.$$
\item The map $\phi$ induces a \kt isomorphism between $L(H)$ and $\Homo_{\reph}(\kk, A)$.
\item The map $\phi$ induces a \kt isomorphism between $N(H)$ and $\Negl_{\reph}(\kk, A)$.
\item The map $\psi$ induces a \kt isomorphism between $V_n(H)$ and $\Homo_{\reph}(A^{\otimes n},\kk)$.
\item $V_1(H)=Z(H)$.
\end{enumerate}
\end{lem}

\begin{rembold}
In Appendix~\ref{sect-app}, we use the space $L(H)$ and the morphism $T$, defined from topological considerations, to
parameterize all the traces of a finite-dimensional ribbon Hopf algebra $H$.
\end{rembold}

\begin{proof}
Let us prove Part (a). Since $(H^*,\cdot)$ is a free right $H$-module of rank 1 with basis $\lambda$, $\phi$ is an
isomorphism. The expression of $\phi^{-1}$ follows from Lemma~\ref{lemdufa}(b).

Let us prove Part (b). Let $z \in H$. We have to show that $z \in L(H)$ if and only if $\phi_z$ is $H$-linear. Suppose
that $\phi_z$ is $H$-linear. For all $a,h \in H$,
\begin{align*}
\varepsilon(h) \, \lambda(za) & = \langle \varepsilon(S^{-1}(h)) \,\phi_z(1_\kk), a \rangle = \langle S^{-1}(h)
\triangleright \phi_z(1_\kk), a \rangle \\
 & = \lambda(zh_{(2)}aS^{-1}(h_{(1)})) = \lambda(S^2(S^{-1}(h_{(1)})\leftharpoonup \nu)zh_{(2)}a) \quad \text{by
 \eqref{lambxy}},
\end{align*}
and so, since $(H^*,\cdot)$ is a free right $H$-module of rank 1 with basis $\lambda$,
\begin{equation}\label{pepe}
\varepsilon(h) z = S^2(S^{-1}(h_{(1)})\leftharpoonup \nu)zh_{(2)}= \nu^{-1}(h_{(2)})\, S(h_{(1)})zh_{(3)}
\end{equation}
for all $h \in H$. Hence, for any $x \in H$,
\begin{align*}
(x\leftharpoonup \nu)z & =  \nu(x_{(1)}) \, x_{(2)} z = \nu(x_{(1)}) \, x_{(2)} \varepsilon(x_{(3)}) z \\
 & = \nu(x_{(1)}) \, \nu^{-1}(x_{(4)})\, x_{(2)}  S(x_{(3)})zx_{(5)} \quad \text{by \eqref{pepe}}\\
 & = \nu(x_{(1)}) \, \nu^{-1}(x_{(3)})\, \varepsilon(x_{(2)}) \, zx_{(4)}
   = \nu(x_{(1)}) \, \nu^{-1}(x_{(2)})\, zx_{(3)} \\
 & = \varepsilon(x_{(1)}) \,  zx_{(2)} = z x,
\end{align*}
and so $z \in L(H)$. Conversely, suppose that $z \in L(H)$. Then, for any $x,h \in H$,
\begin{align*}
\lambda(zS(h_{(1)})xh_{(2)}) & = \lambda(zS^2(h_{(2)})S(h_{(1)})x) \quad \text{by \eqref{eqphihlin}}\\
       & = \lambda(zS(h_{(1)}S(h_{(2)}))x) = \varepsilon(h) \, \lambda(zx),
\end{align*}
and so $\phi_z$ is $H$-linear.

Let us prove Part (d). Note that $\psi$ is an isomorphism since $H$ is finite-dimensional. Let $X \in H^{\otimes n}$. For
all $h \in H$ and $f_1, \dots, f_n \in H^*$,
\begin{equation*}
\psi_X\bigl(h \triangleright (f_1 \otimes \cdots \otimes f_n)\bigr) = \langle f_1 \otimes \cdots \otimes f_n, X
\triangleleft h \rangle
\end{equation*}
and
\begin{equation*}
\varepsilon(h) \, \psi_X(f_1 \otimes \cdots \otimes f_n) = \langle f_1 \otimes \cdots \otimes f_n, \varepsilon(h) \, X
\rangle.
\end{equation*}
Therefore $\psi_X$ is $H$-linear if and only if $X \in V_n(H)$.

Let us prove Part (e). Let $a \in Z(H)$. For all $h \in H$,
\begin{equation*}
a \triangleleft h= S(h_{(1)})a h_{(2)} = S(h_{(1)})h_{(2)} a =\varepsilon(h) \, a
\end{equation*}
and so $a \in V_1(H)$. Conversely, let $a \in V_1(H)$. For all $x \in H$,
\begin{equation*}
xa= x_{(1)}\varepsilon(x_{(2)})\, a = x_{(1)} (a \triangleleft x_{(2)})= x_{(1)} S(x_{(2)})a x_{(3)}
=\varepsilon(x_{(1)}) \, a x_{(2)} = ax,
\end{equation*}
and so $a \in Z(H)$.

Finally, let us prove Part (c). Let $z \in L(H)$. Since $\Endo_{\reph}(\kk)=\kk$, we have that $\phi_z$ is negligible if
and only if $\psi_a \phi_z=\lambda(za)=0$ for all $a \in V_1(H)=Z(H)$, that is, if and only if $z \in N(H)$.
\end{proof}

\begin{lem}\label{lemLH}
\begin{enumerate}
\renewcommand{\labelenumi}{{\rm (\alph{enumi})}}
\item $L(H)$ is a commutative algebra with product $*$ defined by
\begin{align*}
x*z & =\lambda(x S(z_{(2)}) )\, z_{(1)}=\lambda(z S(x_{(2)}) )\, x_{(1)}\\
 & =\lambda(z_{(1)}S^{-1}(x))\, z_{(2)}=\lambda(x_{(1)}S^{-1}(z))\, x_{(2)}
\end{align*}
for any $x,z \in L(H)$, and with $S(\Lambda)$ as unit element.

\item For any $z \in L(H)$, $S_A \circ \phi_z=\phi_{T(z)}$, where $S_A$ denotes the antipode of the categorical Hopf algebra $A$.

\item $T$ induces on $L(H)$ an involutory algebra-automorphism, that is,
\begin{equation*}
T(x * z)=T(x)*T(z), \quad T(S(\Lambda))=S(\Lambda), \quad \text{and} \quad T^2(x)=x
\end{equation*}
for all $x,z \in L(H)$.
\end{enumerate}
\end{lem}
\begin{proof}
Let us prove Part (a). Since $A$ is a Hopf algebra in $\reph$, the space $\Homo_{\reph}(\kk,A)$ is an algebra for the
convolution product $\al * \beta=m_A (\al \otimes \be)$ and with unit element $\eta_A$. This algebra structure transports
to $L(H)$ via the \kt isomorphism $\phi:L(H) \to \Homo_{\reph}(\kk,A)$. Let $x,z \in L(H)$. Then
\begin{align*}
 x * z & = \phi^{-1}(\phi_x * \phi_z) = \phi^{-1}(m_A(\phi_x \otimes \phi_z)) \\
  & = \langle m_A(\phi_x \otimes \phi_z)(1_\kk), \Lambda_{(1)} \rangle \, S(\Lambda_{(2)}) \quad
    \text{by Lemma~\ref{lemdespiptspaces}(a)} \\
  & = \langle \lambda \cdot x \otimes \lambda \cdot z, \Delta^\bd(\Lambda_{(1)}) \rangle \, S(\Lambda_{(2)})
\end{align*}
Write $R=\sum_i a_i \otimes b_i$. By using Lemma~\ref{braidedHopfalg}, \eqref{pptRmat1} and \eqref{eqphihlin}, we have
that
\begin{align*}
 x * z & = \ptsum_i \lambda(x\Lambda_{(2)}a_i) \, \lambda(zS(b_{i(1)})\Lambda_{(1)}b_{i(2)}) \, S(\Lambda_{(3)}) \\
   & = \ptsum_i \lambda(x\Lambda_{(2)}a_i) \, \lambda(zS^2(b_{i(2)})S(b_{i(1)})\Lambda_{(1)}) \, S(\Lambda_{(3)}) \\
   & = \ptsum_i \lambda(x\Lambda_{(2)}a_i \varepsilon(b_i)) \, \lambda(z\Lambda_{(1)}) \, S(\Lambda_{(3)}),
\end{align*}
that is,
\begin{equation}\label{eqcomalg1}
   x * z = \lambda(x\Lambda_{(2)}) \, \lambda(z\Lambda_{(1)}) \, S(\Lambda_{(3)}).
\end{equation}
Likewise
\begin{align*}
 x * z & = \ptsum_i \lambda(xS(a_{i(1)})\Lambda_{(1)}a_{i(2)}) \, \lambda(zS(b_i)\Lambda_{(2)}) \, S(\Lambda_{(3)}) \\
   & = \ptsum_i \lambda(xS^2(a_{i(2)})S(a_{i(1)})\Lambda_{(1)}) \, \lambda(zS(b_i)\Lambda_{(2)}) \, S(\Lambda_{(3)}) \\
   & = \ptsum_i \lambda(x\Lambda_{(1)}) \, \lambda(zS(\varepsilon(a_i)\,b_i)\Lambda_{(2)}) \, S(\Lambda_{(3)}),
\end{align*}
that is,
\begin{equation}\label{eqcomalg2}
   x * z = \lambda(x\Lambda_{(1)}) \, \lambda(z\Lambda_{(2)}) \, S(\Lambda_{(3)}).
\end{equation}
Now, by Lemma~\ref{lemdufa}(a),
\begin{align}
  &z =\lambda(z \Lambda_{(1)}) \, S(\Lambda_{(2)}), \label{eqcomalg3}\\
  &x = \lambda(x \Lambda_{(1)}) \,S(\Lambda_{(2)}), \label{eqcomalg4}
\end{align}
so that
\begin{align}
  & z_{(1)} \otimes S^{-1}(z_{(2)}) =\lambda(z \Lambda_{(1)}) \, S(\Lambda_{(3)}) \otimes \Lambda_{(2)} \label{eqcomalg5}\\
  & x_{(1)} \otimes S^{-1}(x_{(2)}) =\lambda(x \Lambda_{(1)}) \, S(\Lambda_{(3)}) \otimes \Lambda_{(2)}. \label{eqcomalg6}
\end{align}
Hence
\begin{align*}
 x * z & = \lambda(x\Lambda_{(2)}) \, \lambda(z\Lambda_{(1)}) \, S(\Lambda_{(3)}) \quad \text{by \eqref{eqcomalg1}}\\
       & = \lambda(x S^{-1}(z_{(2)})) \, z_{(1)} \quad \text{by \eqref{eqcomalg5}}\\
       & = \lambda(x S(z_{(2)}) )\, z_{(1)} \quad \text{by \eqref{eqphihlin},}
\end{align*}
\begin{align*}
 x * z & = \lambda(z\Lambda_{(1)}) \, \lambda(x\Lambda_{(2)}) \, S(\Lambda_{(3)}) \quad \text{by \eqref{eqcomalg1}}\\
       & = \lambda(z\Lambda_{(1)}) \, \lambda(x_{(1)}\Lambda_{(2)}) \, x_{(2)}\Lambda_{(3)} S(\Lambda_{(4)})\\
       & = \lambda(x_{(1)}S^{-1}(z))\, x_{(2)} \quad \text{by \eqref{eqcomalg3},}
\end{align*}
\begin{align*}
 x * z & = \lambda(x\Lambda_{(1)}) \, \lambda(z\Lambda_{(2)}) \, S(\Lambda_{(3)}) \quad \text{by \eqref{eqcomalg2}}\\
       & = \lambda(z S^{-1}(x_{(2)})) \, x_{(1)} \quad \text{by \eqref{eqcomalg6}}\\
       & = \lambda(z S(x_{(2)}) )\, x_{(1)} \quad \text{by \eqref{eqphihlin},}
\end{align*}
and
\begin{align*}
 x * z & = \lambda(x\Lambda_{(1)}) \, \lambda(z\Lambda_{(2)}) \, S(\Lambda_{(3)}) \quad \text{by \eqref{eqcomalg2}}\\
       & = \lambda(x\Lambda_{(1)}) \, \lambda(z_{(1)}\Lambda_{(2)}) \, z_{(2)}\Lambda_{(3)} S(\Lambda_{(4)})\\
       & = \lambda(z_{(1)}S^{-1}(x))\, z_{(2)} \quad \text{by \eqref{eqcomalg4}.}
\end{align*}
Note that these expressions of the product of $L(H)$ show that $L(H)$ is commutative. Moreover, by
Lemma~\ref{lemdespiptspaces}(a), the unit element of $L(H)$ is
\begin{equation*}
\phi^{-1}(\eta_A)=(\eta_A(1) \otimes S)\Delta(\Lambda)=(\varepsilon \otimes S)\Delta(\Lambda)=S(\Lambda).
\end{equation*}

Let us prove Part (b). Let $z \in L(H)$. Write $R=\sum_i a_i \otimes b_i$. For any $x \in H$,
\begin{align*}
\langle S_A \phi_z(1_\kk),x \rangle & =\langle \lambda \cdot z, S^\bd(x) \rangle \\
 & = \ptsum_i \lambda(z S(a_i) \theta^2 S(x) u b_i) \quad \text{by Lemma~\ref{braidedHopfalg}}\\
 & = \ptsum_i \lambda(z S^2(u)S^2(b_i)S(a_i) \theta^2 S(x)) \quad \text{by \eqref{eqphihlin}}\\
 & = \lambda(z u^2 \theta^2 S(x)) \quad \text{since $S^2(u)=u$ and $(S \otimes S)(R)=R$}\\
 & = \lambda(z G^2S(x)) = \lambda((G^2 \leftharpoonup \nu)z S(x))=\nu(G^2)\,\lambda(G^2z S(x)) \\
 & = \nu(G^2)\,\lambda(G^2 h_\nu x S^{-1}(z) G^{-2}) \quad \text{by \eqref{lambdaSribb}}\\
 & = \nu(G^2)\,\lambda((G^{-2}\leftharpoonup \nu)G^2 h_\nu x S^{-1}(x) ) \quad \text{by \eqref{lambxy}}\\
 & = \nu(G^2)\,\nu(G^{-2})\, \lambda(h_\nu x S^{-1}(z) ) \\
 & = \lambda((S(z)\leftharpoonup \nu)h_\nu x) \quad \text{by \eqref{lambxy}}\\
 & = \lambda(T(z) x) = \langle \phi_{T(z)}(1_\kk),x \rangle,
\end{align*}
that is $S_A \phi_z=\phi_{T(z)}$.

Let us prove Part (c). Let $z \in L(H)$. Firstly $T(z) \in L(H)$ since $\phi_{T(z)}=S_A \phi_z$ is $H$-linear. Moreover,
since $S_A^2=\theta_A$ and the twist of $\reph$ is natural and verify $\theta_\kk=\id_\kk$, we have that
\begin{equation*}
 T^2(z)=\phi^{-1}(\phi_{T^2(z)})=\phi^{-1}(S_A^2 \phi_z)=\phi^{-1}(\theta_A \phi_z)=\phi^{-1}(\phi_z\theta_\kk)
 =\phi^{-1}(\phi_z)=z.
\end{equation*}
For any $x,z \in L(H)$,
\begin{align*}
\phi_{T(x * z)} & = S_A \phi_{x*z} = S_A(\phi_x * \phi_z) = S_A m_A (\phi_x \otimes \phi_z) \\
 & = m_A (S_A \otimes S_A) c_{A,A} (\phi_x \otimes \phi_z) = m_A (S_A \otimes S_A) (\phi_z \otimes \phi_x)c_{\kk,\kk}\\
 & = m_A (S_A\phi_z \otimes S_A\phi_x) = \phi_{T(z)}*\phi_{T(x)}=\phi_{T(z)*T(x)}
\end{align*}
and so $T(x * z)=T(z)*T(x)=T(x)*T(z)$. Finally $T(S(\Lambda))=S(\Lambda)$ since $\phi_{T(S(\Lambda))}=S_A
\phi_{S(\Lambda)}=S_A \eta_A=\eta_A=\phi_{S(\Lambda)}$.
\end{proof}

In the next theorem, we describe the sets $\Ki(\reph)$ and $\Ki(\reph)^\norm$ in algebraic terms. Set
\begin{equation}
  \Ki(H) =\phi^{-1}(\Ki(\reph)) \quad \text{and} \quad
  \Ki(H)^\norm =\phi^{-1}(\Ki(\reph)^\norm),
\end{equation}
where $\phi:L(H) \to \Homo_{\reph}(k,A)$ is as in Lemma~\ref{lemdespiptspaces}. Note that $\Ki(H)^\norm \subset \Ki(H)$.
It follows from \eqref{plusnegl} that
\begin{equation*}
\kk \Ki(H)+N(H) \subset \Ki(H) \quad \text{and} \quad  \kk^* \Ki(H)^\norm+N(H) \subset \Ki(H)^\norm.
\end{equation*}
Since $\eta_A \in \Ki(\reph)^\norm$ and by Lemma~\ref{lemLH}, the right integral $S(\Lambda)$ for $H$ belongs to $\Ki(H)$
and $\Ki(H)^\norm$.
\begin{thm}\label{KpourH}
$\Ki(H)$ is constituted by the element $z \in L(H)$ verifying
\begin{enumerate}
\renewcommand{\labelenumi}{{\rm (\alph{enumi})}}
  \item $T(z)-z \in N(H)$, that is, $\lambda(T(z)a)=\lambda(za)$ for all $a \in Z(H)$;
  \item $\ptsum_i \lambda(zx_{i(1)})\lambda(zx_{i(2)}y_i)=\ptsum_i \lambda(zx_i)\lambda(zy_i)$ for all
        $X=\sum_i x_i \otimes y_i \in V_2(H)$.
\end{enumerate}
In particular, $\Ki(H)$ contains the elements $z \in L(H)$ satisfying $T(z)=z$ and
$\lambda(zx_{(1)})zx_{(2)}=\lambda(zx)z$ for all $x \in H$. Moreover, an element $z \in \Ki(H)$ belongs to $\Ki(H)^\norm$
if and only if $\lambda(z\theta) \neq 0 \neq \lambda(z\theta^{-1})$.
\end{thm}
Note that the sets $\Ki(H)$ and $\Ki(H)^\norm$ do not depend on the choice of the non-zero right integral $\lambda$ for
$H^*$. In Section~\ref{Sect-examples}, we give an example of determination of these sets for a family of non-unimodular
ribbon Hopf algebras.

\begin{proof}
Let $z \in L(H)$. By Lemma~\ref{lemLH}(b), we have that $S_A \phi_z=\phi_{T(z)}$. Therefore, using
Lemma~\ref{lemdespiptspaces}(c), we get that $S_A \phi_z - \phi_z \in \Negl_{\reph}(\kk, A)$ if and only if $T(z)-z \in
N(H)$ , that is, if and only if $\lambda(T(z)a)=\lambda(za)$ for all $a \in Z(H)$. Note that this last property is in
particular satisfied when $T(z)=z$.

Since $\Endo_{\reph}(\kk)=\kk$, the morphism $\Gamma_r(\phi_z \otimes \phi_z)-\phi_z \otimes \phi_z: \kk \to A \otimes A$
is negligible if and only if $G \circ \bigl ( \Gamma_r(\phi_z \otimes \phi_z)-\phi_z \otimes \phi_z \bigr )=0$ for any $G
\in \Homo_{\reph}(A \otimes A,\kk)$. By Lemma~\ref{lemdespiptspaces}(d), this is equivalent to $\psi_X (\Gamma_r(\phi_z
\otimes \phi_z)-\phi_z \otimes \phi_z)=0$ for all $X \in V_2(H)$. Now, writing $R=\sum_i a_i \otimes b_i$ and using
Lemma~\ref{braidedHopfalg}, we have that for any $x,y \in H$,
\begin{eqnarray*}
  \lefteqn{\langle \Gamma_r(\phi_z \otimes \phi_z)(1_\kk), x \otimes y \rangle}\\
   & = &\langle (m_A \otimes \id_A)(\id_A \otimes \Delta_A)(\lambda \cdot z \otimes \lambda \cdot z), x \otimes y
   \rangle\\
   & = &\langle (\id_A \otimes \Delta_A)(\lambda \cdot z \otimes \lambda \cdot z), \Delta^\bd(x) \otimes y \rangle\\
   & = &\ptsum_i \langle (\id_A \otimes \Delta_A)(\lambda \cdot z \otimes \lambda \cdot z),
       S(a_{i(1)}) x_{(1)} a_{i(2)} \otimes S(b_i) x_{(2)} \otimes y \rangle\\
   & = &\ptsum_i \langle \lambda \cdot z \otimes \lambda \cdot z,
       S(a_{i(1)}) x_{(1)} a_{i(2)} \otimes S(b_i) x_{(2)} y \rangle \\
   & = &\ptsum_i \lambda (zS(a_{i(1)}) x_{(1)} a_{i(2)}) \, \lambda(zS(b_i) x_{(2)} y ) \\
   & = &\ptsum_i \lambda (zS^2(a_{i(2)})S(a_{i(1)}) x_{(1)}) \, \lambda(zS(b_i) x_{(2)} y ) \quad \text{by \eqref{eqphihlin}}\\
   & = &\ptsum_i \lambda (z x_{(1)}) \, \lambda(zS(\varepsilon(a_i)b_i) x_{(2)} y ) \\
   & = &\lambda (z x_{(1)}) \, \lambda(z x_{(2)} y ) \quad \text{by \eqref{pptRmat1}.}
\end{eqnarray*}
Therefore the morphism $\Gamma_r(\phi_z \otimes \phi_z)-\phi_z \otimes \phi_z: \kk \to A \otimes A$ is negligible if and
only if $\lambda (z x_{i(1)}) \, \lambda(z x_{i(2)} y_i )= \lambda (z x_i) \, \lambda(z y_i )$ for all $X=\sum x_i
\otimes y_i \in V_2(H)$. Note that this last property is in particular satisfied when $\lambda (z x_{(1)}) \, \lambda(z
x_{(2)} y )=\lambda (z x) \, \lambda(z y )$ for all $x,y \in H$, that is, since $(H^*,\cdot)$ is a free right $H$-module
of rank 1 with basis $\lambda$, when $\lambda (z x_{(1)}) \, z x_{(2)}=\lambda (z x) \, z$ for all $x \in H$.

Finally, by using Lemma~\ref{repcoendlemma}, we have that:
\begin{equation*}
\Theta_\pm \phi_z=\ev_H(\id_{H^*} \otimes \theta_H^{\pm 1})(\phi_z(1_\kk) \otimes 1_H)
 =\ev_H (\lambda \cdot z \otimes \theta^{\pm 1})=\lambda(z \theta^{\pm 1}).
\end{equation*}
This completes the proof of the theorem.
\end{proof}

\begin{cor}\label{corKiKR}
$1 \in \Ki(H)$ if and only if $H$ is unimodular.
\end{cor}
\begin{proof}
Suppose that $H$ is unimodular. Therefore $L(H)=Z(H)$ and $T(z)=S(z)$ for all $z \in Z(H)$. In particular $1 \in L(H)$
and $T(1)=1$. Moreover, $\lambda(x_{(1)})x_{(2)}=\lambda(x)1$ for all $x \in H$ (since $\lambda$ is a right integral for
$H^*$). Hence $1 \in \Ki(H)$ by Theorem~\ref{KpourH}.

Conversely, suppose that $1 \in \Ki(H)$. In particular $1 \in L(H)$ and so $z\leftharpoonup \nu=z$ for all $z \in H$.
Therefore $\varepsilon(z)=\varepsilon(z\leftharpoonup \nu)= \nu(z_{(1)})\varepsilon(z_{(2)})=\nu(z)$ for all $z \in H$,
that is, $H$ is unimodular.
\end{proof}

\subsection{Elements elements from semisimplification}
Let $H$ be a finite-dimen\-sional ribbon Hopf algebra. Let $(A,i)$ be the coend of the functor \eqref{Ffunctor} of
$\reph$ (as in Section~\ref{sectbraidHA}). Denote by $\reph^\sss$ the semisimplification of $\reph$ and by $\pi$ its
associated surjective ribbon functor $\reph \to \reph^\sss$ (see Section~\ref{sect-semisimplific}). Let $\phi:L(H) \to
\Homo_{\reph}(\kk,A)$ be as is Section~\ref{subsect-kirbelemnH}. Set
\begin{equation}
\Ki(H)^\sss=\bigcup_\bb\phi^{-1} \Bigl (\pi^{-1}\bigl(\varphi_\bb(\Ki(\bb))\bigr) \Bigr ),
\end{equation}
where $\bb$ runs over (equivalence classes of) finitely semisimple full ribbon and abelian subcategories of $\reph^\sss$
whose simple objects are scalar, and $\varphi_\bb$ is the morphism \eqref{phibb} corresponding to $\bb$. By
Corollary~\ref{invfromsscat}, we have that
\begin{equation*}
\Ki(H)^\sss \subset \Ki(H).
\end{equation*}
Note that this inclusion may be strict (see Remark~\ref{inclustrict}).

Let $\vv$ be a set of representatives of isomorphism classes of indecomposable finite-dimensional left $H$-modules with
non-zero quantum dimension. Note that $\pi(\vv)=\{\pi(V) \, | \, V \in \vv\}$ is a set of representatives of isomorphism
classes of simple objects of $\reph^\sss$. Let $\lambda$ be a non-zero right integral for $H^*$. Since $H^*$ is a free
right $H$-module with basis $\lambda$ (see Section~\ref{subsecfindimHA}), there exists a (unique) element $z_V \in H$
such that
\begin{equation}\label{zzv}
\lambda(z_V x)=\Tr(G^{-1}x \,\id_V)
\end{equation}
for all $x \in H$, where $G$ is the special grouplike element of $H$. Recall that $\qdim(V)=\Tr(G^{-1} \id_V)$ denotes
the quantum dimension of $V$ (see Lemma~\ref{tevinreplem}).
\begin{cor}\label{corhsss}
\begin{enumerate}
  \renewcommand{\labelenumi}{{\rm (\alph{enumi})}}
  \item If $z \in \Ki(H)^\sss$, then $z=k\sum_{V \in \ww}\qdim(V) \, z_V$ for some finite subset $\ww$ of~$\vv$ and some scalar
  $k \in \kk$.
  \item Let $\ww$ be a (finite) set of representatives of isomorphism classes of simple objects of a full ribbon and
        abelian subcategory of $\reph^\sss$. We can suppose that
        $\ww \subset \pi(\vv)$. If the objects of $\ww$ are scalar, then
        \begin{equation*}
        \sum_{V \in \pi^{-1}(\ww)}\qdim(V) \,
        z_V \in \Ki(H)^\sss.
        \end{equation*}
\end{enumerate}
\end{cor}
\begin{proof}
Let $\bb$ be finitely semisimple full ribbon and abelian subcategory of $\reph^\sss$ whose simple objects are scalar, and
let $(B,j)$ be the coend of the functor \eqref{Ffunctor} of $\bb$ (as in Section~\ref{sect-Kiss}). We can suppose that
there exists a (finite) subset $\ww$ of $\vv$ such that $\pi(\ww)$ is a set of representatives of isomorphism classes of
simple objects of $\bb$. Recall that $B=\oplus_{V \in \ww} \, \pi(V)^* \otimes \pi(V)$. In particular, there exist
morphisms $p_V:B \to \pi(V)^* \otimes \pi(V)$ and $q_V: \pi(V)^* \otimes \pi(V) \to B$ of $\bb$ such that $\id_B= \sum_{V
\in \ww} q_V p_V$ and $p_V q_W = \delta_{V,W}\, \id_{\pi(V)^* \otimes \pi(V)}$. Recall that $j_V=q_V$ for any $V \in
\ww$. Let $\phi:L(H) \to \Homo_{\reph}(\kk,A)$ be as is Section~\ref{subsect-kirbelemnH}. As in \eqref{phibb}, we set
\begin{equation*}
\varphi_\bb=\sum_{V \in \ww} \pi(i_V)p_V \in \Homo_{\reph^\sss}(B,\pi(A)).
\end{equation*}
Let $V\in \ww$. Let $(e_i)_i$ be a basis of $V$ with dual basis $(e_i^*)_i$. By Lemma~\ref{tevinreplem}(b) and
\eqref{imdefit} we have that
\begin{align*}
 \langle i_V \tcoev_V (1_\kk),x \rangle
  & = \langle i_V (\id_{V^*} \otimes G^{-1}  \id_V) \, \sigma_{V,V^*} \, \coev_V (1_\kk),x \rangle\\
  & = \ptsum_i \langle i_V (e_i^* \otimes G^{-1}  e_i) ,x \rangle
  = \ptsum_i \langle e_i^*, xG^{-1}  e_i \rangle\\
  & = \Tr( xG^{-1}  \id_V)
  = \Tr( G^{-1}x \, \id_V)\\
  & =\lambda(z_V x) = \langle \phi_{z_V}(1_\kk),x \rangle,
\end{align*}
for any $x \in H$, that is, $i_V\tcoev_V=\phi_{z_\vv}$. Moreover,
\begin{align*}
\varphi_\bb j_{\pi(V)}\tcoev_{\pi(V)} &=\sum_{W \in \ww} \pi(i_W)p_W j_{\pi(V)}\tcoev_{\pi(V)}\\
 &=  \sum_{W \in \ww} \delta_{V,W}\,  \pi(i_W) \pi(\tcoev_V) \\
 &= \pi (i_V \tcoev_V)=\pi(\phi_{z_V}).
\end{align*}
Hence Part (a) follows from Lemma~\ref{aleqal1dim} and Corollary~\ref{invfromsscat}(a), and Part (b) follows from
Theorem~\ref{RTsubmonoidal} and Corollary~\ref{invfromsscat}(a).
\end{proof}

\begin{lem}\label{KiHnonsslem}
If $H$ is not semisimple, then $\varepsilon(z)=0$ for any $z \in \Ki(H)^\sss$.
\end{lem}
\begin{rembold}\label{inclustrict}
When $H$ is not semisimple, it is possible that $\Ki(H)^\sss \varsubsetneq \Ki(H)$. For example, if $H$ is unimodular but
not semisimple, then $1 \in \Ki(H)$ (by Corollary~\ref{corKiKR}) and $1 \not \in \Ki(H)^\sss$ (by
Lemma~\ref{KiHnonsslem}, since $\varepsilon(1)=1$).
\end{rembold}
\begin{proof}
Let $\Lambda$ be a left integral for $H$ such that $\lambda(\Lambda)=1$. Since $H$ is not semisimple, we have
$\varepsilon (\Lambda)=0$ (by \cite[Theorem 3.3.2]{Abe}) and $\Lambda^2=\varepsilon (\Lambda)\Lambda=0$. Now, if $M$ is a
finite-dimensional left $H$-module, then $(\Lambda\, \id_M)^2=\Lambda^2 \id_M=0$ and so $\Tr(\Lambda \,\id_M)=0$. Let $z
\in \Ki(H)^\sss$. By Corollary~\ref{corhsss}(a), there exist $k \in \kk$ and a finite subset $\ww$ of~$\vv$ such that
$z=k\sum_{V \in \ww}\qdim(V) \, z_V$. Then
\begin{equation*}
\lambda(z \Lambda)=k\sum_{V \in \ww} \qdim(V) \,\Tr(G^{-1} \Lambda\, \id_V) = k\sum_{V \in \ww} \qdim(V)
\varepsilon(G^{-1}) \,\Tr( \Lambda \,\id_V)=0.
\end{equation*}
Hence $\varepsilon(z)=\varepsilon(z)\lambda(\Lambda)=\lambda(\varepsilon(z)\Lambda)=\lambda(z \Lambda)=0$.
\end{proof}

Recall (see \cite{CR}) that $\kk$ is a \emph{splitting field} for a \kt algebra $A$ if every simple finite dimensional
left $A$-module is scalar. Note that this is always the case if $\kk$ is algebraically closed.

\begin{cor}
If $H$ is semisimple and $\kk$ is a splitting field for $H$, then $\Ki(H)^\sss=\Ki(H)$ and this set is composed by
elements $z \in Z(H)$ satisfying $S(z)=z$ and $\lambda(zx_{(1)})zx_{(2)}=\lambda(zx)z$ for all $x \in H$.
\end{cor}
\begin{proof}
Since $H$ is semisimple and finite-dimensional, we have that $\reph^\sss=\reph$ and that $\vv$ finite. Moreover, since
$\kk$ is a splitting field for a $H$, every $V \in \vv$ is scalar. Then $\Ki(H)=\phi^{-1}(\Ki(\reph)) \subset
\Ki(H)^\sss$ and so $\Ki(H)^\sss=\Ki(H)$. Moreover, since $H$ is unimodular (because it is semisimple), we have that
$L(H)=Z(H)$, $N(H)=0$, and $T(x)=S(x)$ for all $x \in H$ (see Section~\ref{subsect-kirbelemnH}). Therefore, by using
Theorem~\ref{KpourH}, we get that $z \in \Ki(H)$ if and only if $z \in Z(H)$, $S(z)=z$, and
$\lambda(zx_{(1)})zx_{(2)}=\lambda(zx)z$  for all $x \in H$.
\end{proof}

\begin{prop}\label{hsscossksfc0}
Suppose that $H$ is semisimple, that $\kk$ is of characteristic 0, and that $\kk$ is a splitting field for $H$. Then
$Z(H)=\bigoplus_{V \in \vv} \, \kk z_V$ and $\sum_{V \in \vv}  \qdim(V)\, z_V \in \kk^* 1$.
\end{prop}
\begin{proof}
Note that $H$ is cosemisimple since any finite-dimensional semisimple Hopf algebra over a field of characteristic 0 is
cosemisimple (see \cite[Theorem 3.3]{LR1}). Then $S^2=\id_H$ by \cite[Theorem 4]{LR2} and  $\lambda(1)\neq 0$
by~\cite[Theorem 3.3.2]{Abe}.

Note that $\vv$ is finite (since $H$ is finite-dimensional). By \cite[Theorem 25.10]{CR}, any simple left ideal of $H$ is
isomorphic (as a left $H$-module) to a (unique) element of $\vv$. For any $V \in \vv$, let $H_V \subset H$ be the sum of
all the simple left ideals of $H$ which are isomorphic to $V$. By \cite[Theorem 25.15]{CR}, $H_V$ is a two-sided ideal of
$H$, $H_V$ is simple \kt algebra (the operations being those induced by $H$), $B_VB_W=0$ for $V \neq W \in \vv$,
$H=\bigoplus_{V \in \vv} H_V$, and $H$ is isomorphic (as an algebra) to $\prod_{V \in \vv} H_V$. Moreover, if $e_V$
denotes the unit of $H_V$, then $1=\sum_{V \in \vv} e_V$, $H_V=He_V$, and $e_Ve_W=\delta_{V,W} \, e_V$ for all $V,W \in
\vv$. By \cite[Theorem 26.4]{CR}, since $H_V$ is a simple \kt algebra, $V$ is a simple left $H_V$-module, and
$\Endo_{H_V}(V)=\kk$ (because $V$ is a scalar $H$-module), we have that $H_V$ is isomorphic (as an algebra) to
$\Endo_\kk(V)$ and that $\kdim(V)$ is the number of simple left ideals appearing in a direct sum decomposition of $H_V$
as such a sum. Then $Z(H_V)=\kk e_V$ (since $Z(\Endo_\kk(V))=\kk \id_V$) and so $Z(H)=\bigoplus_{V \in \vv} \kk e_V$.
Moreover, for any $x \in H$, we have that
\begin{equation}\label{eqHtrss2}
\Tr(x \,\id_H)=\sum_{V \in \vv} \Tr(x \,\id_{H_V})=\sum_{V \in \vv} \kdim(V)\, \Tr(x \,\id_V).
\end{equation}

The map $x \in H \mapsto \Tr((x\,\id_H)\circ S^2) \in \kk$ is a right integral for $H^*$ (by \cite[Proposition
2(b)]{Rad1} applied to $H^\mathrm{op}$). Therefore, since $S^2=\id_H$ and  by the uniqueness of integrals, there exists
$k \in \kk$ such that $\Tr(x \,\id_H)=k \, \lambda(x)$ for all $x \in H$. Evaluating with $x=1$, we get that
$k=\kdim(H)/\lambda(1)$. Hence, by \eqref{eqHtrss2},
\begin{equation}\label{eqHtrss1}
k \, \lambda(x)=\sum_{V \in \vv} \kdim(V)\, \Tr(x \,\id_V)
\end{equation}
for all $x \in H$.

Let $V \in \vv$. By \eqref{ScarreG} and since $S^2=\id_H$, the special grouplike element $G$ of $H$ is central and so
$G^{-1} \id_V$ is $H$-linear. Therefore, since $V$ is scalar and $G$ is invertible, there exists a (unique) $\gamma_V \in
\kk^*$ such that $G^{-1} \id_V=\gamma_V \id_V$. Since $H$ and $H^*$ are semisimple and so unimodular, their special
grouplike elements are trivial. Then $G^2=1$ by \eqref{gG2h} and so $\gamma_V^2=1$. Hence, for all $x \in H$,
\begin{equation}\label{eqHtrss3}
\begin{split}
\qdim(V)\, \Tr(G^{-1} x \,\id_V)&=\Tr(G^{-1} \id_V)\, \Tr(xG^{-1}\id_V)\\
&=\gamma_V^2\Tr(\id_V)\, \Tr(x \,\id_V)=\kdim(V) \,\Tr(x \,\id_V).
\end{split}
\end{equation}

For any $x \in H$, we have that:
\begin{align*}
\lambda(\qdim(V)\, z_V x) &= \qdim(V)\,\Tr(G^{-1}x \,\id_V) \quad \text{by \eqref{zzv}}\\
& = \kdim(V)\,\Tr(x \,\id_V) \quad \text{by \eqref{eqHtrss3}}\\
& = \ptsum_{W \in \vv} \kdim(W)\,\Tr(x e_V\,\id_W) \quad \text{since $e_V \,\id_W=\delta_{V,W} \,\id_V$}\\
& = \lambda(k\, e_Vx) \quad \text{by \eqref{eqHtrss1},}
\end{align*}
and so $\qdim(V)\, z_V=k\, e_V$ (since $H^*$ is a free right $H$-module with basis~$\lambda$).

Finally, since $k=\kdim(H)/\lambda(1)\neq0$ (because the characteristic of $\kk$ is 0) and $\qdim(V) \neq 0$ (because
$\reph$ is semisimple, see Section~\ref{sect-defsscat}), we get that $Z(H)=\bigoplus_{V \in \vv} \kk e_V=\bigoplus_{V \in
\vv} \, \kk z_V$ and $\sum_{V \in \vv} \qdim(V) \, z_V=k \sum_{V \in \vv} e_V=k \,1 \in \kk^*1$.
\end{proof}

\subsection{HKR-type invariants}\label{subsect-HKRtyp}
Let $H$ be a finite-dimensional ribbon Hopf algebra. We use notations of Section~\ref{subsect-kirbelemnH}.

By Theorem~\ref{thm3man} and Proposition~\ref{prop3man}, for any $z \in \Ki(H)^\norm$,
\begin{equation}\label{deftauHz}
 \tau_{(H,z)}(M)=\tau_{\reph}(M;\phi_z) \in \kk
\end{equation}
is an invariant of 3-manifolds. Note that the choice of the normalisation in the definition of $\tau_{\reph}(M;\phi_z)$
(see Proposition~\ref{prop3man}) implies that $\tau_{(H,z)}(M)$ does not depend on the choice of the non-zero right
integral $\lambda$ for $H^*$ used to define $\tau_{\reph}(M;\phi_z)$.

By Remarks~\ref{rmsdeKi}, for any $z \in \Ki(H)^\norm$, we have that
\begin{align*}
&\tau_{(H,z)}(S^3)=1,  &&\tau_{(H,z)}(M\#M')=\tau_{(H,z)}(M)\, \tau_{(H,z)}(M'),
\end{align*}
for all 3-manifolds $M$, $M'$. Moreover, if $z \in \Ki(H)^\norm$, $n \in N(H)$, and $k \in \kk^*$, then $kz+n \in
\Ki(H)^\norm$ and, for all 3-manifold $M$,
\begin{equation}\label{taukzeqtauz}
\tau_{(H,k z+n)}(M)=\tau_{(H,z)}(M).
\end{equation}

\begin{defn}
An invariant of closed 3-manifolds $I$ with values in $\kk$ is said to be of \emph{HKR-type} if there exist a
finite-dimensional ribbon Hopf algebra $H$ (over~$\kk$) and an element $z \in \Ki(H)^\norm$  such that
$I(M)=\tau_{(H,z)}(M)$ for all 3-manifolds $M$.
\end{defn}
In Proposition~\ref{corRTvsHKR}, we show that the Reshetikhin-Turaev invariants defined with premodular Hopf
algebras are of HKR-type.\\

Let us show that any HKR-type invariant can be computed by using the Kauffman-Radford algorithm (which is given in
\cite{KR1} for the case $H$  unimodular and for $z=1$). Fix a finite-dimensional ribbon Hopf algebra $H$, a non-zero
right integral $\lambda$ for $H^*$, and an element $z \in \Ki(H)^\norm$. Let $M$ be a 3-manifold and $L=L_1 \cup \cdots
\cup L_n$ be a framed link in $S^3$ such that $M \simeq M_L$. Let us recall the Kauffman-Radford algorithm\footnote{This
corresponds to \cite{KR1} applied with the opposite ribbon Hopf algebra $H^\mathrm{op}$ to $H$.}:

(A). Consider a diagram $D$ of $L$ (with blackboard framing). Each crossing of $D$ is decorated with the \R matrix
$R=\sum_i a_i \otimes b_i$ as in Figure~\ref{algoHKRfigur1}. The diagram obtained after this step is called the
\emph{flat diagram of $D$}. Note that the flat diagram of $D$ is composed by $n$ closed plane curves, each of them
arising from a component of $L$.
\begin{figure}[h]
    \begin{center}
      \psfrag{x}{$a_i$}
      \psfrag{y}{$b_i$}
      \psfrag{S}{$\ds \sum_i$}
      \psfrag{u}[Br][Br]{$S(a_i)$}
      \scalebox{.8}{\includegraphics{HKRalg1.eps}}
    \end{center}
      \caption{}
      \label{algoHKRfigur1}
\end{figure}

(B). On each component of the flat diagram of $D$, the algebraic decoration is concentrated in an arbitrary point (other
than extrema and crossings) according to the rules of Figure~\ref{algoHKRfigur2}, where $a,b \in H$.
\begin{figure}[h]
    \begin{center}
      \psfrag{z}{$a$}
      \psfrag{x}{$a$}
      \psfrag{y}{$S(a)$}
      \psfrag{z}{$a$}
      \psfrag{r}[Br][Br]{$a$}
      \psfrag{s}[Br][Br]{$b$}
      \psfrag{c}{$ab$}
      \psfrag{i}[Br][Br]{$1$}
      \psfrag{t}[Br][Br]{$S(a)$}
      \psfrag{d}{$G^{-1}$}
      \psfrag{e}{$G$}
      \scalebox{.8}{\includegraphics{HKRalg2.eps}}
    \end{center}
      \caption{}
      \label{algoHKRfigur2}
\end{figure}

\noindent In that way we get an element $\sum_k v_1^k \otimes \cdots \otimes v_n^k \in H^{\otimes n}$, where $v_i^k$
corresponds to the component of the flat diagram of $D$ arising from $L_i$, see Figure~\ref{algoHKRfigur3}.
\begin{figure}[h]
    \begin{center}
      \psfrag{a}[Br][Br]{$v_1^k$}
      \psfrag{b}{$v_2^k$}
      \psfrag{d}{$v_n^k$}
      \psfrag{z}{$a$}
      \psfrag{r}[Br][Br]{$L_1$}
      \psfrag{s}{$L_2$}
      \psfrag{c}{$L_n$}
      \psfrag{-}{$\dots$}
      \psfrag{S}{$\ds \sum_k$}
      \scalebox{.8}{\includegraphics{HKRalg3.eps}}
    \end{center}
      \caption{}
      \label{algoHKRfigur3}
\end{figure}

\noindent For $1 \leq i \leq n$, let $d_i$ be the Whitney degree of the flat diagram of $L_i$ obtained by traversing it
upwards from the point where the algebraic decoration have been concentrated. The Whitney degree is the total turn of the
tangent vector to the curve when one traverses it in the given direction, see Figure~\ref{algoHKRfigur4}.
\begin{figure}[h]
    \begin{center}
      \psfrag{k}{$d=1$}
      \psfrag{s}{$d=-1$}
      \psfrag{c}{$d_i=-2$}
      \psfrag{u}[Br][Br]{$v_i^k$}
      \scalebox{.8}{\includegraphics{withdeg.eps}}
    \end{center}
      \caption{}
      \label{algoHKRfigur4}
\end{figure}

\begin{prop}\label{calcbyHKRalgo}
\begin{equation*}
\tau_{(H,z)}(M)=\lambda(z\theta)^{b_-(L)-n_L} \, \lambda(z\theta^{-1})^{-b_-(L)}\, \sum_k \lambda(zG^{d_1+1}v_1^k) \cdots
\lambda(zG^{d_n+1}v_n^k).
\end{equation*}
\end{prop}
\begin{proof}
Choose an orientation for $L$. Let $T$ be a $n$-special tangle such that $L$ is isotopic $T \circ (\cup_- \otimes \cdots
\otimes \cup_-)$,  where the $i$th cup (with clockwise orientation) corresponds to the component $L_i$. Let $D_T$ be a
diagram of $T$. Note that $D=D_T \circ (\cup \otimes \cdots \otimes \cup)$ is a diagram of $L$. Apply steps (A) and (B)
to $D_T$, see Figure~\ref{algoHKRfigur5}. Note that, in this case, $d_i=-1$.
\begin{figure}[h]
    \begin{center}
      \psfrag{u}{$v_1^k$}
      \psfrag{v}{$v_n^k$}
      \psfrag{=}{$\dots$}
      \psfrag{T}{$T$}
      \psfrag{S}{$\ds \sum_k$}
      \scalebox{.8}{\includegraphics{HKRalg4.eps}}
    \end{center}
      \caption{}
      \label{algoHKRfigur5}
\end{figure}

From the definition of the monoidal structure, duality, braiding and twist of $\reph$ (see Section~\ref{sect-hmodu}), it
is not difficult to verify that, for any finite-dimensional left $H$-modules $M_1, \dots, M_n$,
\begin{equation*}
T_{(M_1, \dots, M_n)}=\sum_k \ev_{M_1}(\id_{M_1^*} \otimes v_k^1 \,\id_{M_1}) \otimes \cdots \otimes
     \ev_{M_n}(\id_{M_n^*} \otimes v_k^n \,\id_{M_n}).
\end{equation*}
Then, by Lemma~\ref{repcoendlemma},
\begin{align*}
\tau_{\reph}(L;\phi_z) & =\phi_T \circ \phi_z^{\otimes n } = \sum_k \ev_{H}(\phi_z(1_\kk) \otimes v_k^1 ) \otimes \cdots
\otimes
     \ev_{H}(\phi_z(1_\kk) \otimes v_k^n ) \\
     & =\sum_k \lambda(zv_1^k) \cdots \lambda(zv_n^k)=
     \sum_k \lambda(zG^{d_1+1}v_1^k) \cdots \lambda(zG^{d_n+1}v_n^k).
\end{align*}
Hence we get the results since $\Theta_\pm \phi_z=\lambda(z\theta^{\pm 1})$.
\end{proof}
\begin{cor}\label{invHKRpouHop}
Suppose that $H$ is unimodular and $\lambda(\theta) \neq 0 \neq \lambda(\theta^{-1})$. Then  $1 \in \Ki(H)^\norm$ and
$\tau_{(H,1)}(M)$ is the Hennings-Kauffman-Radford invariant of 3-manifolds defined with the opposite ribbon Hopf algebra
$H^\mathrm{op}$ to $H$.
\end{cor}
\begin{proof}
This is an immediate consequence of Corollary~\ref{corKiKR}, Proposition~\ref{calcbyHKRalgo}, and the definition of the
Hennings-Kauffman-Radford invariant given, e.g., in \cite{KR1}.
\end{proof}

\subsection{Reshetikhin-Turaev from premodular Hopf algebras}

Let $(H,\vv)$ be a finite-dimensional \emph{premodular} Hopf algebra. This means (see \cite{Tur2}) that $H$ is a
finite-dimensional ribbon Hopf algebra and $\vv$ is a finite set of finite-dimensional pairwise non-isomorphic left
$H$-modules such that:
\begin{itemize}
\item each $V \in \vv$ is non-negligible scalar;
\item the trivial left $H$-module $\kk$ belongs to $\vv$;
\item for any $V \in \vv$, there exists $W \in \vv$ such that $V^*\simeq W$;
\item for any $V,W \in \vv$, $V \otimes W$ splits as a (finite) direct sum of certain modules of $\vv$ (possibly with
multiplicities) and a negligible $H$-module.
\end{itemize}
By \emph{negligible} $H$-module, we mean a finite-dimensional left $H$-module $N$ such that $\tr(f)=0$ for any $f \in
\Endo_{\reph}(N)$ or, equivalently, such that $\qdim(N)=0$.

Consider the semisimplification $\reph^\sss$ of $\reph$ (see Section~\ref{sect-semisimplific}) and let $\pi : \reph \to
\reph^\sss$ be the ribbon functor associated to this semisimplification.

Let $\bb_\vv$ be the full subcategory of $\reph^\sss$ whose objects are finite directs sums of objects of
$\pi(\vv)=\{\pi(V) \, | \, V \in \vv\}$. By the definition of a premodular Hopf algebra, we have that $\bb_\vv$ is a full
ribbon subcategory of $\reph^\sss$. Note that $\bb_\vv$ is a finitely semisimple ribbon \kt category whose simples
objects are scalar and has $\pi(\vv)$ as a (finite) set of representatives of isomorphism classes of simple objects.
Recall that the Reshetikhin-Turaev invariant $\RT_{\bb_\vv}(M)$ of 3-manifolds is well-defined when $\Delta_\pm^{\bb_\vv}
\neq 0$ (see Section~\ref{sectinvarRT}).

Let $\lambda$ be a non-zero right integral for $H^*$. For any $V \in \vv$, as in \eqref{zzv}, we let $z_V \in H$ such
that $\lambda(z_V x)=\Tr(G^{-1}x \,\id_V)$ for all $x \in H$. Set
\begin{equation}\label{zzvv}
z_\vv=\sum_{V \in \vv} \qdim(V) \,z_V,
\end{equation}
where $\qdim(V)=\Tr(G^{-1}\id_V)$ is the quantum dimension of $V$. By Corollary~\ref{corhsss}(b), we have that $z_\vv \in
\Ki(H)$.
\begin{prop}\label{corRTvsHKR}
If $\Delta_\pm^{\bb_\vv} \neq 0$, then $z_\vv \in \Ki(H)^\norm$ and $\tau_{(H,z_\vv)}(M)=\RT_{\bb_\vv}(M)$ for all
3-manifold $M$.
\end{prop}
Note that Proposition~\ref{corRTvsHKR} says that the Reshetikhin-Turaev invariant defined from a premodular Hopf algebra
is of HKR-type.
\begin{proof}
Let $(A,i)$ be the coend of the functor \eqref{Ffunctor} of $\reph$ (as in Section~\ref{sectbraidHA}), let $(B,j)$ be the
coend of the functor \eqref{Ffunctor} of $\bb_\vv$ (as in Section~\ref{sect-Kiss}). Set $\al_{\bb_\vv}=\sum_{V \in \vv}
\qdim(V) j_{\pi(V)}\tcoev_{\pi(V)}$. Suppose that $\Delta_\pm^{\bb_\vv} \neq 0$. By Corollary~\ref{corRt}, we have that
$\al_{\bb_\vv} \in \Ki(\bb_\vv)^\norm$ and $\RT_{\bb_\vv}(M)=\tau_{\bb_\vv}(M;\al_{\bb_\vv})$ for all 3-manifold $M$. Set
$\varphi_{\bb_\vv}:B \to \pi(A)$ as in \eqref{phibb}. As in the proof of Corollary~\ref{corhsss}, we have that
$\pi(\phi_{z_V})=\varphi_{\bb_\vv} j_{\pi(V)}\tcoev_{\pi(V)}$. Then $\pi(\phi_{z_\vv})=\varphi_{\bb_\vv} \al_{\bb_\vv}$.
Since $\al_{\bb_\vv} \in \Ki(\bb_\vv)^\norm$ and $\pi(\phi_{z_\vv})=\varphi \al_{\bb_\vv}$,
Corollary~\ref{invfromsscat}(b) gives that $\phi_{z_\vv} \in \Ki(\reph)^\norm$ and
$\tau_{\bb_\vv}(M;\al_{\bb_\vv})=\tau_{\reph}(M;\phi_{z_\vv})$ for all 3-manifold $M$. Hence $z_\vv \in \Ki(H)^\norm$ and
$\tau_{(H,z_\vv)}(M)=\RT_{\bb_\vv}(M)$ for all 3-manifold $M$.
\end{proof}

Note that if $H$ is semisimple finite-dimensional ribbon Hopf algebra, $\kk$ is a splitting field for $H$, and $\vv$ is a
set of representatives of isomorphism classes of simple left $H$-modules, then $(H,\vv)$ is a premodular Hopf algebra and
$\bb_\vv=\reph$.

\begin{cor}\label{corHssRTeqHKR}
Let $H$ be a finite-dimensional semisimple ribbon Hopf algebra. Suppose that the base field $\kk$ is of characteristic
$0$ and is a splitting field for $H$. Then the Hennings-Kauffman-Radford invariant of 3-manifolds computed with
$H^\mathrm{op}$ and the Reshetikhin-Turaev invariant of 3-manifolds computed with $\reph$ are simultaneously well-defined
(that is, $\Delta_\pm^{\reph} \neq 0$ if and only if $1 \in \Ki(H)^\norm$). Moreover, if they are, then they coincide,
that is, $\tau_{(H,1)}(M)=\RT_{\reph}(M)$ for all 3-manifold $M$.
\end{cor}
\begin{rembold}
Conclusions of Corollary~\ref{corHssRTeqHKR} may be no more true when $H$ is not semisimple (see
Remark~\ref{inclustrict}). Moreover, in the modular case (in the sense of Remark~\ref{rmsthsubmonoRT}),
Corollary~\ref{corHssRTeqHKR} was first shown in \cite{Ker1}.
\end{rembold}
\begin{proof}
By Proposition~\ref{corRTvsHKR}, the Reshetikhin-Turaev invariant of 3-manifolds computed from $\reph$ is well-defined if
$z_\vv=\sum_{V \in \vv} \qdim(V) \, z_V \in \Ki(H)^\norm$ and is equal to $\tau_{(H,z_\vv)}$. By
Corollary~\ref{invHKRpouHop}, the Hennings-Kauffman-Radford invariant of 3-manifolds computed with $H^\mathrm{op}$ is
well-defined if $1 \in \Ki(H)^\norm$ and is equal to $\tau_{(H,1)}$. Now, by Proposition~\ref{hsscossksfc0}, $\sum_{V \in
\vv} \qdim(V)\, z_V=k\, 1$ for some $k \in \kk^*$. We conclude by using \eqref{taukzeqtauz}.
\end{proof}

\section{A non-unimodular example}\label{Sect-examples}
Let us examine the case of family of non-unimodular ribbon Hopf algebras, defined by Radford \cite{Rad4}, which includes
Sweedler's Hopf algebra.

Let $n$ be an odd positive integer and $\kk$ be a field whose characteristic does not divide $2n$. Let $H_n$ be the \kt
algebra generated by $a$ and $x$ with the following relations
\begin{equation*}
a^{2n}=1, \quad x^2=0, \quad ax=-xa.
\end{equation*}
The algebra $H_n$ is a Hopf algebra for the following structure maps:
\begin{align*}
& \Delta(a)=a \otimes a, && \Delta(x)= x \otimes a^n+  1 \otimes x,\\
& \varepsilon(a)=1,  && \varepsilon(x)=0, \\
& S(a)=a^{-1},  &&S(x)=a^nx.
\end{align*}
The set $\mathcal{B}=\{a^lx^m \, | \, 0 \leq l < 2n, \; 0 \leq m \leq 1 \}$ is a basis for $H_n$. The dual basis of
$\mathcal{B}$ is $\{\overline{a^kx^r} \, | \, 0 \leq k < 2n, \; 0 \leq r \leq 1 \}$, where
$\overline{a^kx^n}(a^lx^r)=\delta_{l,k} \, \delta_{m,r}$. Set
\begin{equation*}
 \Lambda=(1+a+a^2+ \cdots + a^{2n-1})x \quad \text{and} \quad \lambda=\overline{a^nx}.
\end{equation*}
Then $\Lambda$ is a left integral for $H_n$ and $\lambda$ is a right integral for $H_n^*$ such that
$\lambda(\Lambda)=\lambda(S(\Lambda))=1$. The distinguished grouplike element of $H_n$ is $g=a^n\in G(H_n)$ and the
distinguished grouplike element $\nu\in G(H_n^*)=\mathrm{Alg}(H_n,\kk)$ of $H_n^*$ is given by $\nu(a)=-1$ and
$\nu(x)=0$.

Suppose that $\kk$ has a primitive $2n$-root of unity $\omega$. Let an odd integer $1 \leq s < 2n$ and  $\be \in \kk$.
Then
\begin{equation*}
R_{\omega,s,\be}=\frac{1}{2n} \sum_{0\leq i,l < 2n} w^{-il} \, a^i \otimes a^{sl} + \frac{\be}{2n} \sum_{0\leq i,l < 2n}
w^{-il} \, a^ix \otimes a^{sl+n}x
\end{equation*}
is a \R matrix for $H_n$ and $h_\nu=(\id_{H_n} \otimes \nu)(R_{\omega,s,\be})=a^n$.

The set $\{e_l=\frac{1}{2n}\sum_{0 \leq i < 2n} \omega^{-il} a^i \, | \, 0\leq l < 2n \}$ is a basis of the algebra
$\kk[a]$ and verifies:
\begin{align*}
& e_ie_j=\delta_{i,j}\, e_i, && 1=e_0+ \cdots + e_{2n-1}, && a^i e_j=\omega^{ij}e_j, \\
& \varepsilon(e_l)=\delta_{l,0}, &&\Delta(e_l)=\sum_{0 \leq i < 2n} e_i \otimes e_{l-i}, && S(e_l)=e_{-l}.
\end{align*}
Denote by $\chi: \kk^* \to \kk[a]$ the algebra homomorphism defined by
\begin{equation*}
\al \in \kk^* \mapsto \chi(\al)=\sum_{0 \leq l < 2n} \al^{l^2} \, e_l.
\end{equation*}
The quasitriangular Hopf algebra $(H_n,R_{\omega,s,\be})$ is ribbon with twist $\theta = a^n \chi(\omega^s)$. Note that
$\theta^p=a^{pn} \chi(\omega^{ps})$ for any integer $p$. The special grouplike element of $H_n$ is then $G=\sum_{0 \leq l
< 2n} (-1)^l \, e_l=a^n$.

Let $T:H_n \to H_n$, $L(H_n)$, $N(H_n)$, $V_2(H_n)$, $\Ki(H_n)$, and $\Ki(H_n)^\norm$ be as in
Section~\ref{subsect-kirbelemnH}. It is not difficult to verify that
\begin{align*}
  & T(a^k)=(-1)^k a^{n-k} \quad \text{and} \quad T(a^kx)=a^{-k}x \quad \text{for all $0 \leq k <2n$,}\\
  & L(H_n)=\kk x \oplus \kk ax \oplus \cdots \oplus \kk a^{2n-1}x, \\
  & Z(H_n)=\kk 1 \oplus \kk a^2 \oplus \cdots \oplus \kk a^{2n-2}, \\
  & N(H_n)=\kk x \oplus \kk a^2x \oplus \cdots \oplus \kk a^{2n-2}x=Z(H_n)x, \\
  & V_2(H_n)=\underset{0\leq p,q <n}{\oplus} \kk (a^{2p} \otimes a^{2q}) \oplus
           \underset{0 \leq k,l<2n}{\oplus} \kk (a^kx \otimes a^lx).
\end{align*}

Let $z \in L(H_n)$. Since $a^{2n}=1$, we can write $z=\sum_{k \in \Z / 2n \Z} \al_k a^k x$ for some function $\al:\Z / 2n
\Z \to \kk$. Using Theorem~\ref{KpourH}, we get that $z \in \Ki(H_n)$ if and only if
\begin{equation*}
\al_{-k}=\al_{k} \quad \text{and} \quad \al_k \al_{k+l-n}=\al_k\al_l.
\end{equation*}
for all $k,l$ odd. Let $\gamma:\al:\Z / n \Z \to \kk$ defined by $\gamma_k=\al_{2k+n}$. Since $n$ is odd, we have that
$z-\sum_{k \in \Z / n \Z} \gamma_k a^{2k+n} x \in N(H_n)$. Then $z \in \Ki(H_n)$ if and only if
\begin{equation}\label{eqgamma}
\gamma_{-k}=\gamma_{k} \quad \text{and} \quad \gamma_k \gamma_{k+l}=\gamma_k\gamma_l.
\end{equation}
for all $k,l$.

Suppose that $z \in \Ki(H_n)$ and $z \neq 0$. Using \eqref{eqgamma}, we get that $\gamma_k\gamma_0=\gamma_k
\gamma_{-k}=\gamma_k^2$ for all $k$. In particular, $\gamma_0 \neq 0$ and $\gamma_k=\gamma_0$ whenever $\gamma_k \neq 0$.
Set
\begin{equation*}
d=\min \{ 1 \leq k \leq n \, | \, \gamma_k \neq 0 \} \geq 1.
\end{equation*}
Note that $\gamma_k=0$ for all $1 \leq k < d$ and, by $\eqref{eqgamma}$, that $\gamma_{k+d}=\gamma_k$ for all $k$. The
integer $d$ divides $n$. Indeed, let $r$ such that $rd \leq n <rd+d$. Then $0 \leq n-rd \leq n$ and
$\gamma_{n-rd}=\gamma_n=\gamma_0 \neq 0$. Therefore, by definition of $d$, we get that $n-rd=0$ and so $d | n$. Hence
$z=\gamma_0 z_d +w$, where $w \in N(H_n)$ and
\begin{equation*}
z_d=\sum_{k=0}^{\frac{n}{d}-1} a^{2dk+n}x.
\end{equation*}

Conversely, it is not difficult to verify  that $z_d \in \Ki(H_n)$ whenever $d$ divides $n$. Hence
\begin{equation*}
 \Ki(H_n)=\bigcup_{d |n} \{ \al z_d +w \, | \, \al \in \kk \text{ and } w \in N(H_n) \}.
\end{equation*}

Let $d$ be an positive integer dividing $n$. For any $\al \in \kk$ and $w \in N(H_n)$, we have that
\begin{equation*}
\lambda((\al z_d +w)\theta^{\pm 1})=\frac{\al}{2d}\,  \sum_{k=0}^{2d-1} \, \bigl (\omega^\frac{n}{d} \bigr )^{\pm s
\frac{n}{d} k^2 +n k}.
\end{equation*}
The sum of the right-hand sum of this equality is a Gauss sum which is non-zero if and only if the enhancement $k \in
\Z/2d\Z \mapsto \psi(k)=\pm s \frac{n}{d} k^2 +n k \in \Z/2d\Z$ is tame, that is, $\psi(x)=0$ for any $x \in \Z/2d\Z$
such that $\psi(x+y)=\psi(x)+\psi(y)$ for all $y \in \Z/2d\Z$, see \cite{Tay}. Since $n$ and $s$ are odd, it is not
difficult to verify that $\psi$ is tame. Therefore $\lambda((\al z_d +w)\theta^{\pm 1})\neq 0$ if and only if $\al \neq
0$. Hence
\begin{equation*}
 \Ki(H_n)^\norm=\bigcup_{d|n} \{ \al z_d +w \, | \, \al \in \kk^* \text{ and } w \in N(H_n) \}.
\end{equation*}

In conclusion, the ribbon Hopf algebra $H_n$ leads to $D(n)$ HKR-type invariants of 3-manifolds, where $D(n)$ denotes the
number of positive divisors of $n$, which are $\tau_{(H_n,z_d)}$ where $1 \leq d \leq n$ and $d|n$.

Note that $H_n$ is not unimodular (and so neither semisimple) since $\nu\neq \varepsilon$. Therefore $1 \not\in \Ki(H_n)$
(by Corollary~\ref{corKiKR}), that is, the Hennings-Kauffman-Radford invariant is not defined for $H_n$. Moreover, the
categorical Hopf algebra $A=H_n^*$ of $\rephn$ does not possess any non-zero two-sided integral (since $H_n$ is not
unimodular), and so the Lyubashenko invariant of 3-manifolds is not defined for $\rephn$.

\appendix \section{Traces on ribbon Hopf algebras}\label{sect-app}
Recall that a \emph{trace} on a Hopf algebra $H$ is a form $t \in H^*$ such that $t(xy)=t(yx)$ and $t(S(x))=t(x)$ for all
$x,y \in H$.

Let $H$ be a finite-dimensional ribbon Hopf algebra. Let $\lambda \in H^*$ be a non-zero right integral for $H^*$, $\nu
\in G(H^*)=\mathrm{Alg}(H,\kk)$ be the distinguished grouplike element of $H^*$ and $G$ be the special grouplike element
of $H$. Recall that $\cdot$ denotes the right action of $H$ on $H^*$ defined in \eqref{actionpoint}, that
$\leftharpoonup$ denotes the right $H^*$-action on $H$ defined in \eqref{actharpoon}, that $L(H)$ denotes the \kt
subspace of $H$ constituted by the elements $z \in H$ verifying $(x \leftharpoonup \nu) z = z x$ for all $x \in H$, and
that $T$ denotes the \kt endomorphism of $H$ defined by $z \mapsto T(z)=(S(z)\leftharpoonup \nu)h_\nu$, where
$h_\nu=(\id_H \otimes\nu)(R) \in G(H)$.\\

In the next proposition, we give an algebraic description of the space of traces on $H$.
\begin{prop}
The space $\{ z \in L(H) \, | \, T(z)=z \}$ is \kt isomorphic to the space of traces on $H$ via the map $z \mapsto
\lambda \cdot (zG)$.
\end{prop}
If $H$ is unimodular, then $L(H)=Z(H)$ and $T=S$, and so we recover the parameterization of traces on $H$ given in
\cite{Rad1,He2}.
\begin{proof}
Let $z \in L(H)$ such that $T(z)=z$. Set $t=\lambda \cdot (zG) \in H^*$. By using \eqref{eqphihlin} and \eqref{ScarreG},
we have that for any $x,y \in H$,
\begin{equation*}
t(xy)  = \lambda(zGxy)
   = \lambda(z S^2(y)Gx)
  = \lambda(z Gyx) = t(yx).
\end{equation*}
Moreover, for any $x \in H$,
\begin{align*}
t(S(x))&= \lambda( zGS(x)) = \lambda(G^2h_\nu x G^{-1} S^{-1}(z)) \quad \text{by \eqref{lambdaSribb}}\\
  & = \lambda((S(z)\leftharpoonup \nu)G^2h_\nu x G^{-1} ) \quad \text{by \eqref{lambxy}}\\
  & = \lambda(z h_\nu^{-1}G^2h_\nu x G^{-1} ) \quad \text{since $T(z)=z$}\\
  & = \lambda(z  S^2(G^{-1}) h_\nu^{-1}G^2h_\nu x ) \quad \text{by \eqref{eqphihlin}}\\
  & = \lambda(z  G^{-1} h_\nu^{-1}G^2h_\nu x ) \\
  & = \lambda(z  G S^{-4}(h_\nu^{-1})h_\nu x) \quad \text{by \eqref{ScarreG}}\\
  & = \lambda(z  G h_\nu^{-1}h_\nu x) = \lambda(z  G  x) = t(x).
\end{align*}

Conversely, let $t \in H^*$ be a trace on $H$. Since $(H^*,\cdot)$ is a free right $H$-module of rank 1 with basis
$\lambda$ and $G$ is invertible, there exists a (unique) $z \in H$ such that $t=\lambda \cdot (zG)$. Let $x \in H$. For
all $a \in H$,
\begin{align*}
\lambda(zxGy) &= \lambda(z  G S^{-2}(x)y) \quad \text{by \eqref{ScarreG}}\\
  & = t( S^{-2}(x)y) = t(yS^{-2}(x))\\
  & = \lambda(zG yS^{-2}(x))\\
  & = \lambda((x\leftharpoonup \nu)zG y) \quad \text{by \eqref{lambxy}}.
\end{align*}
Therefore, since $(H^*,\cdot)$ is a free right $H$-module of rank 1 with basis $\lambda$, we have that $(x\leftharpoonup
\nu)zG=zxG$ and so $(x\leftharpoonup \nu)z=zx$. Hence $z \in L(H)$. Moreover, for all $x \in H$,
\begin{align*}
\lambda(zGx) &= t(x)=t(S(x))=\lambda(zGS(x))\\
  & = \lambda(G^2h_\nu x G^{-1} S^{-1}(z)) \quad \text{by \eqref{lambdaSribb}} \\
  & = \lambda(G^2h_\nu x S^{-1}(zG)) \\
  & = \lambda(G^2h_\nu x S^{-1}((G\leftharpoonup \nu)z)) \quad \text{since $z \in L(H)$}\\
  & = \lambda((S((G\leftharpoonup \nu)z)\leftharpoonup \nu)G^2h_\nu x ) \quad \text{by \eqref{lambxy}}.
\end{align*}
Therefore, since $(H^*,\cdot)$ is a free right $H$-module of rank 1 with basis $\lambda$,
\begin{align*}
 zG &= (S((G\leftharpoonup \nu)z)\leftharpoonup \nu)G^2h_\nu\\
  & = (S(z)\leftharpoonup \nu)(S(G\leftharpoonup \nu)\leftharpoonup \nu)G S^2(h_\nu)G\quad \text{by \eqref{ScarreG}}\\
  & = (S(z)\leftharpoonup \nu) \nu(G) \nu(G^{-1}) G^{-1} G h_\nu G\\
  & = (S(z)\leftharpoonup \nu) h_\nu G = T(z)G
\end{align*}
and so $T(z)=z$.
\end{proof}

\bibliographystyle{amsalpha}
\bibliography{virel2}
\end{document}